\documentstyle[amssymb,11pt]{book}

\begin{document}

\begin{titlepage}
\begin{center}
{\large
Department of Theoretical Physics}
\vspace*{25mm}

{\huge\rm\bf

The Whittaker model of the center of the quantum group
and Hecke algebras

}
\vspace*{10mm}

{\large
BY

\vspace*{2mm}
ALEXEI SEVOSTYANOV \\
}

\vfill
\vspace*{10mm}

{\large
UPPSALA UNIVERSITY  1999
}

\end{center}

\end{titlepage}

\newpage

\thispagestyle{plain}

{\flushleft Dissertation for the Degree of Doctor of 
Philosophy in Theoretical Physics 
presented at Uppsala University in 1999}
\vspace*{10mm}

%{{{ Abstract

\noindent
{\rm ABSTRACT}
 
\vspace*{2ex}
\noindent
Sevostyanov, A. 1999. The Whittaker model of the center of the quantum group
and Hecke algebras. 
77 pp. Uppsala. ISBN 91-506-1342-1.

\vspace*{2ex}
\noindent
{\small In 1978 Kostant suggested the {\em Whittaker model}
of the center of the universal enveloping algebra $U({\frak g})$
of a complex simple  Lie algebra ${\frak g}$. An essential role
in this construction is played by a non--singular character $\chi$
of the maximal nilpotent subalgebra ${\frak n}_+ \subset \frak{g}$.
The main result is that the center of $U({\frak g})$ is isomorphic
to a commutative subalgebra in $U({\frak b}_-)$, where 
${\frak b}_- \subset {\frak g}$ is the opposite Borel subalgebra. 
This observation is used in the theory of principal series
representations of the corresponding Lie group $G$ and in 
the proof of complete integrability of the quantum Toda lattice.

We show that the Whittaker model introduced by Kostant has a natural
homological interpretation in terms of Hecke algebras. Moreover,
we introduce a general definition of a Hecke algebra $Hk^*(A,B,\chi)$
associated to the triple of an associative algebra $A$, 
a subalgebra $B \subset A$ and a character $\chi$ of $B$. 
In particular, the Whittaker model of the center of $U({\frak g})$ is identified 
with
$Hk^0(U({\frak g}), U({\frak n}_+), \chi)^{opp}$.

The goal of this thesis is to generalize the Kostant's construction
to quantum groups. An obvious obstruction is the fact that the
subalgebra in $U_h({\frak g})$ 
generated by positive root generators (subject to the quantum Serre relations)
does not have non--singular characters. In order to overcome
this difficulty we introduce a family of new realizations of
quantum groups, one for each Coxeter element of the corresponding
Weyl group. The modified quantum Serre relations allow for non--singular
characters, and we are able to construct the Whittaker model
of the center of $U_h({\frak g})$.

The new Whittaker model is applied to the deformed quantum Toda
lattice recently studied by Etingof. We give new proofs of his results
which resemble the original Kostant's proofs for the quantum Toda lattice.

Finally, we study the ``quasi-classical'' limit of the Whittaker model
for $U_h({\frak g})$. A remarkable new result is a cross-section theorem
for the action of a complex simple  Lie group on itself by conjugations. 
We are able to prove this theorem for all such Lie groups except for
the case of $E_6$! Using the cross--section theorem we establish a relation 
between the
Whittaker model and the set of conjugacy classes of regular elements in
the corresponding Lie group $G$.}

\vspace*{0.5ex}
{\flushleft \em Alexei Sevostyanov, Department of Theoretical Physics, 
Uppsala University,
Box 803, S-751\,08 Uppsala, Sweden}

\newpage

\renewcommand{\theequation}{\thesection.\arabic{equation}}

\newtheorem{theorem}{Theorem}{}
\newtheorem{lemma}[theorem]{Lemma}{}
\newtheorem{corollary}[theorem]{Corollary}{}
\newtheorem{conjecture}[theorem]{Conjecture}{}
\newtheorem{proposition}[theorem]{Proposition}{}
\newtheorem{axiom}{Axiom}{}
\newtheorem{remark}{Remark}{}
\newtheorem{example}{Example}{}
\newtheorem{exercise}{Exercise}{}
\newtheorem{definition}{Definition}{}

\renewcommand{\thetheorem}{\thesection.\arabic{theorem}}

\renewcommand{\thelemma}{\thesection.\arabic{lemma}}

\renewcommand{\theproposition}{\thesection.\arabic{proposition}}

\renewcommand{\thecorollary}{\thesection.\arabic{corollary}}

\renewcommand{\theremark}{\thesection.\arabic{remark}}

\renewcommand{\thedefinition}{\thesection.\arabic{definition}}

\setcounter{equation}{0}
\setcounter{theorem}{0}

\thispagestyle{plain}

\begin{center}
{\LARGE \bf Introduction}
\end{center}

In 1978 Kostant suggested the {\em Whittaker model}
of the center of the universal enveloping algebra $U({\frak g})$
of a complex simple  Lie algebra ${\frak g}$. An essential role
in this construction is played by a non--singular character $\chi$
of the maximal nilpotent subalgebra ${{\frak n}_+} \subset \frak{g}$.
The main result is that the center of $U({\frak g})$ is isomorphic
to a commutative subalgebra in $U({\frak b}_-)$, where 
${\frak b}_- \subset {\frak g}$ is the opposite Borel subalgebra. 
This observation is used in the theory of principal series
representations of the corresponding Lie group $G$ and in 
the proof of complete integrability of the quantum Toda lattice.

We show that the Whittaker model introduced by Kostant has a natural
homological interpretation in terms of Hecke algebras. Moreover, following 
\cite{S2}
we introduce a general definition of a Hecke algebra $Hk^*(A,B,\chi)$
associated to the triple of an associative algebra $A$, 
a subalgebra $B \subset A$ and a character $\chi$ of $B$. 
In particular, the Whittaker model of the center of $U({\frak g})$ is identified 
with
$Hk^0(U({\frak g}), U({{\frak n}_+}), \chi)^{opp}$.

The goal of this thesis is to generalize the Kostant's construction
to quantum groups. An obvious obstruction is the fact that the
subalgebra in $U_h({\frak g})$ 
generated by positive root generators (subject to the quantum Serre relations)
does not have non-singular characters. In order to overcome
this difficulty we introduce a family of new realizations of
quantum groups, one for each Coxeter element of the corresponding
Weyl group (see also \cite{S1}). The modified quantum Serre relations allow for 
non--singular
characters, and we are able to construct the Whittaker model
of the center of $U_h({\frak g})$.

The new Whittaker model is applied to the deformed quantum Toda
lattice recently studied by Etingof (see \cite{Et}). We give new proofs of his 
results
which resemble the original Kostant's proofs for the quantum Toda lattice.

Finally, we study the ``quasi-classical'' limit of the Whittaker model
for $U_h({\frak g})$. A remarkable new result is a cross-section theorem
for the action of a complex simple  Lie group on itself by conjugations.
We are able to prove this theorem for all such Lie groups except for
the case of $E_6$! This theorem is a group counterpart of the cross-section
theorem of Kostant (Theorem C, Section \ref{geomappr}). Using the cross--section
theorem we 
establish a relation between the
Whittaker model and the set of conjugacy classes of regular elements in
the corresponding Lie group $G$.

The thesis is organized as follows. Chapter \ref{Witt} contains a review of
Kostant's results on the Whittaker model \cite{K}, \cite{K1}. 
In order to create a pattern for 
proofs in the  quantum group case we recall most of the Kostant's proofs.
Chapter \ref{Hom} is devoted to Hecke algebras. It contains the definition of
the algebra $Hk^*(A,B,\chi)$ and the interpretation of the Whittaker model
as $Hk^0(U({\frak g}), U({{\frak n}_+}), \chi)^{opp}$.
The central part of the thesis is Chapter \ref{qWitt}. There we describe
new realizations of finite-dimensional quantum groups and present
the Whittaker model of the center of $U_h({\frak g})$. Chapter \ref{qWitt}
also contains a discussion of the deformed quantum Toda lattice.
In Chapter \ref{GWitt} we establish a relation between the Whittaker model and
regular elements in algebraic groups. The main result of this Chapter
is a cross-section theorem for the action of a complex simple Lie
group on itself by conjugations (see also \cite{SS} where we prove a 
modification of this 
theorem for loop groups). 

Many results presented in this thesis for finite-dimensional
quantum groups have natural counterparts for affine quantum groups.
In order to simplify the presentation we treat only the finite-dimensional
case, and refer the reader to the papers \cite{S1}, \cite{FRS}, \cite{SS}
for further details on the affine case.

\thispagestyle{plain}

\tableofcontents

%%%%%%%%%%%%%%%%%%%%%%%%%%%%%%%%%%%%%%%%%%%%%%%%%%%%%%%%%%%%%%%%%%%%%%%%%%%%%%%%
%%%%%%%%%%%%%%%
%%%%%%%%%%%%%%%%%%%%%%%%%%%%%%%%%%%%%%%%%%%%%%%%%%%%%%%%%%%%%%%%%%%%%%%%%%%%%%%%
%%%%%%%%%%%%%%
%%%%%%%%%%%%%%%%%%%%%%%%%%%%%%%%%%%%%%%%%%%%%%%%%%%%%%%%%%%%%%%%%%%%%%%%%%%%%%%%
%%%%%%%%%%%%%%%%

\chapter{Whittaker model}\label{Witt}

In this chapter we recall the Whittaker model of the center
of the universal enveloping algebra $U({\frak g})$, where ${\frak g}$ is a 
complex simple Lie algebra.

%%%%%%%%%%%%%%%%%%%%%%%%%%%%%%%%%%%%%%%%%%%%%%%%%%%%%%%%%%%%%%%%%%%%%%%%%%%%%%%%
%%%%%%%%%%

\section{Notation}\label{notation}

Fix the notation used throughout of the text.
Let $G$ be a
connected simply connected finite--dimensional complex simple Lie group, $%
{\frak g}$ its Lie algebra. Fix a Cartan subalgebra ${\frak h}\subset {\frak %
g}\ $and let $\Delta $ be the set of roots of $\left( {\frak g},{\frak h}%
\right) .$ Choose an ordering in the root system. Let $\alpha_i,~i=1,\ldots 
l,~~l=rank({\frak g})$ be the 
simple roots, $\Delta_+=\{ \beta_1, \ldots ,\beta_N \}$  
the set of positive roots. Denote by $\rho$ a half of the sum of positive roots,
$\rho=\frac 12 \sum_{i=1}^N\beta_i$.
Let $H_1,\ldots ,H_l$ be the set of simple root generators of $\frak h$. 

Let $a_{ij}$ be the corresponding Cartan matrix.
Let $d_1,\ldots , d_l$ be coprime positive integers such that the matrix 
$b_{ij}=d_ia_{ij}$ is symmetric. There exists a unique non--degenerate invariant
symmetric bilinear form $\left( ,\right) $ on ${\frak g}$ such that 
$(H_i , H_j)=d_j^{-1}a_{ij}$. It induces an isomorphism of vector spaces 
${\frak h}\simeq {\frak h}^*$ under which $\alpha_i \in {\frak h}^*$ corresponds 
to $d_iH_i \in {\frak h}$. We denote by $\alpha^\vee$ the element of $\frak h$ 
that 
corresponds to $\alpha \in {\frak h}^*$ under this isomorphism.
The induced bilinear form on ${\frak h}^*$ is given by
$(\alpha_i , \alpha_j)=b_{ij}$.

Let $W$ be the Weyl group of the root system $\Delta$. $W$ is the subgroup of 
$GL({\frak h})$ 
generated by the fundamental reflections $s_1,\ldots ,s_l$,
$$
s_i(h)=h-\alpha_i(h)H_i,~~h\in{\frak h}.
$$
The action of $W$ preserves the bilinear form $(,)$ on $\frak h$. 
We denote a representative of $w\in W$ in $G$ by
the same letter. For $w\in W, g\in G$ we write $w(g)=wgw^{-1}$.

Let ${{\frak b}_+}$ be the positive Borel subalgebra and ${\frak b}_-$
the opposite Borel subalgebra; let ${\frak n}_+=[{{\frak b}_+},{{\frak b}_+}]$ 
and $%
{\frak n}_-=[{\frak b}_-,{\frak b}_-]$ be their 
nil-radicals. Let $H=\exp {\frak h},N_+=\exp {{\frak n}_+},
N_-=\exp {\frak n}_-,B_+=HN_+,B_-=HN_-$ be
the Cartan subgroup, the maximal unipotent subgroups and the Borel subgroups
of $G$ which correspond to the Lie subalgebras ${\frak h},{{\frak n}_+},%
{\frak n}_-,{\frak b}_+$ and ${\frak b}_-,$ respectively.

We identify $\frak g$ and its dual by means of the canonical invariant bilinear 
form. 
Then the coadjoint 
action of $G$ on ${\frak g}^*$ is naturally identified with the adjoint one. We 
also identify 
${{\frak n}_+}^*\cong {\frak n}_-,~{{\frak b}_+}^*\cong {\frak b}_-$. 

Let ${\frak g}_\beta$ be the root subspace corresponding to a root $\beta \in 
\Delta$, 
${\frak g}_\beta=\{ x\in {\frak g}| [h,x]=\beta(h)x \mbox{ for every }h\in 
{\frak h}\}$.
${\frak g}_\beta\subset {\frak g}$ is a one--dimensional subspace. 
It is well--known that for $\alpha\neq -\beta$ the root subspaces ${\frak 
g}_\alpha$ and ${\frak g}_\beta$ are orthogonal with respect 
to the canonical invariant bilinear form. Moreover ${\frak g}_\alpha$ and 
${\frak g}_{-\alpha}$
are non--degenerately paired by this form.

Root vectors $X_{\alpha}\in {\frak g}_\alpha$ satisfy the following relations:
$$
[X_\alpha,X_{-\alpha}]=(X_\alpha,X_{-\alpha})\alpha^\vee.
$$

If $V$ is a finite--dimensional complex vector space, $S(V)$ will denote the 
symmetric algebra over $V$ and
$S_k(V)$ denotes the homogeneous subspace of degree k. If $V^*$ is the dual 
space to $V$ then $S(V^*)$ is
regarded as the algebra of polynomial functions on $V$.

Let $U({\frak g})$ be the universal enveloping algebra of $\frak g$, and
$U_k({\frak g})$ the standard filtration in $U({\frak g})$. From the 
Poincar\'{e}--Birkhoff--Witt 
theorem it follows that the associated graded algebra $GrU({\frak g})$ is 
isomorphic to the symmetric algebra 
$S({\frak g})$ of the linear space $\frak g$. 

Equip  $S({\frak g})$ with a Poisson structure as follows. For each 
$s_k\in S_k({\frak g})$ choose a representative $u_k\in U_k({\frak g})$ such 
that 
$u_k/ U_{k-1}({\frak g})=s_k$. We shall denote $s_k=Gru_k$.
Given two such elements $s_i$ and $s_j$ with chosen representatives 
$u_i$ and $u_j$, the commutativity of $S({\frak g})$ implies that 
$$
[u_i,u_j]\in U_{i+j-1}({\frak g}).
$$
Define 
\begin{equation}\label{KK}
\{ s_i,s_j\} =[u_i,u_j]/ U_{i+j-2}({\frak g}).
\end{equation}
It is easy to see that this bracket is independent of the choice of 
representatives $u_i,~u_j$ and equips
$S({\frak g})$ with the structure of a Poisson algebra, i.e. it is a derivation 
of the multiplication in 
$S({\frak g})$. We refer to the procedure described above as the graded limit.

%%%%%%%%%%%%%%%%%%%%%%%%%%%%%%%%%%%%%%%%%%%%%%%%%%%%%%%%%%%%%%%%%%%%%%%%%%%%%%%%
%%%%%%%%%%%%%%%%%%%%

\section{The Whittaker model}\label{whitt}

\setcounter{equation}{0}
\setcounter{theorem}{0}

In this section we introduce the Whittaker model of the center of the universal 
enveloping
algebra $U({\frak g})$. We start by recalling the classical result of Chevalley
which describes the structure of the center.

Let $Z({\frak g})$ be the center of $U({\frak g})$. The standard filtration 
$U_k({\frak g})$ 
in $U({\frak g})$ induces a filtration $Z_k({\frak g})$ in $Z({\frak g})$. The 
following
important theorem may be found for instance in \cite{Bur1}, Ch.8, \S 8, no. 3, 
Corollary 1 and no.5, Theorem 2.
\vskip 0.3cm
\noindent
{\bf Theorem (Chevalley)}
{\em One can choose
elements $I_k\in Z_{m_k+1}({\frak g}),~~k=1,\ldots l$, where $m_k$ are called 
the exponents of $\frak g$,
such that
$Z({\frak g})={\Bbb C}[I_1,\ldots , I_l]$ is a polynomial algebra in $l$ 
generators.}
\vskip 0.3cm
The adjoint action of $G$ on $\frak g$ naturally extends to $S({\frak g})$. 
Let $S({\frak g})^G$ be the algebra of $G$--invariants in $S({\frak g})$.
Clearly, $GrZ({\frak g})\cong S({\frak g})^G$. In particular $S({\frak 
g})^G\cong 
{\Bbb C}[\widehat I_1,\ldots , \widehat I_l]$, where $\widehat I_i=Gr 
I_i,~i=1,\ldots ,l$. The elements
$\widehat I_i,~i=1,\ldots ,l$ are called fundamental invariants.

Following Kostant we shall realize the center $Z({\frak g})$ of the universal 
enveloping 
algebra $U({\frak g})$ as a subalgebra in $U({\frak b}_-)$.
Let 
$$
\chi :{{\frak n}_+} \rightarrow {\Bbb C}
$$
be a character of ${{\frak n}_+}$. 
Since ${{\frak n}_+}=\sum_{i=1}^l{\Bbb C}X_{\alpha_i} \oplus[{{\frak 
n}_+},{{\frak n}_+}]$
it is clear that $\chi$ is completely determined by the constants $c_i=\chi 
(X_{\alpha_i}),~i=1,\ldots ,l$ and 
$c_i$ are arbitrary. In \cite{K} $\chi$ is called non--singular if $c_i\neq 0$ 
for all $i$.

Let $f=\sum_{i=1}^l X_{-\alpha_i}\in {\frak n}_-$ be a regular nilpotent 
element.
From the properties of the invariant bilinear form (see Section \ref{notation}) 
it follows that
$(f,[{{\frak n}_+},{{\frak 
n}_+}])=0,~~(f,X_{\alpha_i})=(X_{-\alpha_i},X_{\alpha_i})$, 
and hence the map
$x\mapsto (f,x),~~x\in {{\frak n}_+}$ is a non--singular character of ${{\frak 
n}_+}$.

Recall that in our choice of root vectors no normalization was made. But now 
given a non--singular
character $\chi :{{\frak n}_+}\rightarrow {\Bbb C}$ we will say that $f$ 
corresponds to $\chi$ in case
$$
\chi (X_{\alpha_i}) =(X_{-\alpha_i},X_{\alpha_i}).
$$
Conversely if $\chi$ is non--singular there is a unique choice of $f$ so that 
$f$ corresponds to $\chi$. 
In this case $\chi (x)=(f,x)$ for every $x\in {{\frak n}_+}$.

Naturally, the character $\chi$ extends to a character of 
the universal enveloping algebra $U({{\frak n}_+})$. 
Let $U_\chi ({{\frak n}_+})$ be the kernel of this extension so that 
one has a direct sum
$$
U({{\frak n}_+})={\Bbb C}\oplus U_\chi ({{\frak n}_+}).
$$

 Since ${\frak g}={\frak b}_-\oplus {{\frak n}_+}$ we have a linear 
isomorphism $U({\frak g})=U({\frak b}_-)\otimes U({{\frak n}_+})$ and hence 
the direct sum
\begin{equation}\label{maindec}
U({\frak g})=U({\frak b}_-) \oplus I_\chi,
\end{equation}
where $I_\chi=U({\frak g})U_\chi ({{\frak n}_+})$ is the left--sided ideal 
generated by 
$U_\chi ({{\frak n}_+})$.

For any $u\in U({\frak g})$ let $u^\chi\in U({\frak b}_-)$ be its component in 
$U({\frak b}_-)$ relative to the decomposition (\ref{maindec}). Denote by 
$\rho_\chi$
the linear map 
$$
\rho_\chi : U({\frak g}) \rightarrow U({\frak b}_-)
$$
given by $\rho_\chi (u)=u^\chi$.
Let $W({\frak b}_-)=\rho_\chi (Z({\frak g}))$. 
\vskip 0.3cm
\noindent
{\bf Theorem A (\cite{K}, Theorem 2.4.2)}
{\em The map
\begin{equation}\label{map}
\rho_\chi : Z({\frak g}) \rightarrow W({\frak b}_-)
\end{equation}
is an isomorphism of algebras. In particular 
$$
W({\frak b}_-)={\Bbb C}[I_1^\chi ,\ldots , I_l^\chi ], 
~~I_i^\chi=\rho_\chi(I_i),~~i=1,\ldots ,l
$$
is a polynomial algebra in $l$ generators.}
\vskip 0.3cm
\noindent
{\em Proof.}
First, we show that the map (\ref{map}) is an algebra homomorphism.
If $u,v\in Z({\frak g})$ then $u^\chi v^\chi \in U({\frak b}_-)$ and
$$
uv-u^\chi v^\chi =(u-u^\chi )v+u^\chi (v-v^\chi ).
$$
Since $(u-u^\chi )v=v(u-u^\chi )$ the r.h.s. of the last equality is an element 
of $I_\chi$.
This proves $u^\chi v^\chi =(uv)^\chi$.

By definition the map (\ref{map}) is surjective. We have to prove that it is 
injective.
Let $U({\frak g})^{\frak h}$ be the centralizer of $\frak h$ in $U({\frak g})$. 
Clearly 
$Z({\frak g})\subseteq U({\frak g})^{\frak h}$. From the 
Poincar\'{e}--Birkhoff--Witt 
theorem it follows that every element $z\in U({\frak g})^{\frak h}$ may be 
uniquely
written as
$$
z=\sum_{p,q\in{\Bbb N}^N,<p>=<q>}X_{-\beta_1}^{p_1}\ldots 
X_{-\beta_N}^{p_N}\varphi_{p,q}
X_{\beta_1}^{q_1}\ldots X_{\beta_N}^{q_N},
$$
where $<p>=\sum_{i=1}^r p_i \beta_i \in {\frak h}^*$ and $\varphi_{p,q} \in 
U({\frak h})$. 

Now recall that $\chi (X_{\beta_i})=0$ if $\beta_i$ is not a simple root, and we 
easily obtain
$$
\rho_\chi (z)=\sum_{p,q\in{\Bbb N}^l,<p>=<q>\neq 
0}X_{-\alpha_{k_1}}^{p_{j_1}}\ldots X_{-\alpha_{k_l}}^{p_{j_l}}\varphi_{p,q}
\prod_{i=1}^lc_{k_i}^{q_{j_i}}+\varphi_{0,0}.
$$

Let $z\in Z({\frak g})$. One knows that the map 
$$
Z({\frak g})\rightarrow U({\frak h}),~~z\mapsto \varphi_{0,0},
$$
called the Harich-Chandra homomorphism, is injective (see (c), p. 232 in 
\cite{Dix}). It follows that 
the map (\ref{map}) is also injective. 
\vskip 0.3cm
\noindent
{\bf Definition A}
{\em The algebra $W({\frak b}_-)$ is called the Whittaker model of $Z({\frak 
g})$.}  
\vskip 0.3cm

Next we equip $U({\frak b}_-)$ with a structure of a left $U({{\frak n}_+})$ 
module in such a
way that $W({\frak b}_-)$ is realized as the space of invariants with respect to 
this action.

Let $Y_\chi$ be the left $U({\frak g})$ module defined by 
$$
Y_\chi =U({\frak g})\otimes_{U({{\frak n}_+})}{\Bbb C}_\chi ,
$$
where ${\Bbb C}_\chi$ denotes the 1--dimensional $U({{\frak n}_+})$--module 
defined by $\chi$.
Obviously $Y_\chi$ is just the quotient module $U({\frak g})/I_\chi$. From 
(\ref{maindec}) it 
follows that the map
\begin{equation}\label{iso2}
U({\frak b}_-)\rightarrow Y_\chi;~~v\mapsto v\otimes 1
\end{equation}
is a linear isomorphism. 

It is convenient
to carry the module structure of $Y_\chi$ to $U({\frak b}_-)$. For $u\in 
U({\frak g}),~~
v\in U({\frak b}_-)$ the induced action $u\circ v$ has the form
\begin{equation}\label{indact}
u\circ v=(uv)^\chi.
\end{equation}
The restriction of this action to $U({{\frak n}_+})$ may be changed by tensoring 
with 1--dimensional
$U({{\frak n}_+})$--module defined by $-\chi$. That is $U({\frak b}_-)$ becomes 
an $U({{\frak n}_+})$
module where if $x\in {U({\frak n}_+)},~~v\in U({\frak b}_-)$ one puts
\begin{equation}\label{mainact}
x\cdot v=x\circ v-\chi (x)v.
\end{equation}
\vskip 0.3cm
\noindent
{\bf Lemma A (\cite{K}, Lemma 2.6.1.)}
{\em Let $v\in U({\frak b}_-)$ and $x\in {U({\frak n}_+)}$. Then}
$$
x\cdot v =[x,v]^\chi.
$$
\vskip 0.3cm
\noindent
{\em Proof.} 
By definition $x\cdot v=(xv)^\chi -\chi (x)v$. Then we have $xv=[x,v]+vx$ and 
hence
$x\cdot v=([x,v])^\chi +(vx)^\chi -\chi (x)v$. But clearly $(vx)^\chi =v\chi 
(x)$. Thus
$x\cdot v=([x,v])^\chi$.

The action (\ref{mainact}) may be lifted to an action of the unipotent group 
$N_+$.
Consider the space $U({\frak b}_-)^{N_+}$ of $N_+$ invariants 
in $U({\frak b}_-)$ with respect to this
action. Clearly, 
$W({\frak b}_-)\subseteq U({\frak b}_-)^{N_+}$ . 
\vskip 0.3cm
\noindent
{\bf Theorem B (\cite{K}, Theorems 2.4.1, 2.6)}
{\em Suppose that the character $\chi$ is non--singular. Then the space of $N_+$ 
invariants 
in $U({\frak b}_-)$ with respect to the
action (\ref{mainact}) is isomorphic to $W({\frak b}_-)$, i.e.}
\begin{equation}\label{inv}
U({\frak b}_-)^{N_+}\cong W({\frak b}_-).
\end{equation}
\vskip 0.3cm
We shall prove Theorem B in the next section.

%%%%%%%%%%%%%%%%%%%%%%%%%%%%%%%%%%%%%%%%%%%%%%%%%%%%%%%%%%%%%%%%%%%%%%%%%%%%%%%%
%%%%%%%%%%%%%%%%%

\section{Geometric approach to the Whittaker model}\label{geomappr}

\setcounter{equation}{0}
\setcounter{theorem}{0}

In this section we establish a relation between the Whittaker model
and the geometry of the adjoint action of
the corresponding Lie group. 

Denote the character of $S({{\frak n}_+})$ that equals to $\chi(x)$ for every 
$x\in {\frak n}_+$ 
by the same letter. Similarly to (\ref{maindec}) we have the following 
decomposition for $S({\frak g})$:
$$
S({\frak g})=S({\frak b}_-)\oplus I_\chi^0,
$$
where $I_\chi^0$ is the ideal in $S({\frak g})$ generated by the kernel of 
$\chi$.
For any $s\in S({\frak g})$ let $s^\chi$ be its component in $S({\frak b}_-)$ 
relative to 
this decomposition.
 
Now using Lemma A we define the graded limit of the action (\ref{mainact}). For 
$x\in {\frak n}_+$
and $s\in S({\frak b}_-)$ we put
\begin{equation}\label{classact}
x\cdot s=(\{ x,s\} )^\chi.
\end{equation}

This action may be lifted to an action of the unipotent group $N_+$ on $S({\frak 
b}_-)$. 
For $a\in N_+,~s\in S({\frak b}_-)$ this action is given by 
$$
a\cdot s=({\rm Ad}(a)(s))^\chi.
$$ 
Observe that $S({\frak b}_-)$ is naturally identified with the algebra of 
polynomial functions on
${\frak b}_+$. We shall describe the space of invariants $S({\frak b}_-)^{N_+}$ 
using the induced action of 
$N_+$ on ${\frak b}_+$.

To calculate this action it suffices to consider the restriction of the action 
(\ref{classact}) to linear
functions. Let $s\in {\frak b}_-$ be such a function. Then for $a\in N_+$ 
$$
a\cdot s=({\rm Ad}(a)(s))^\chi =P_{{\frak b}_-}({\rm Ad}(a)(s))+\chi(P_{{\frak 
n}_+}({\rm Ad}(a)(s))),
$$
where $P_{{\frak b}_-}$ and $P_{{\frak n}_+}$ are the projection operators onto 
${\frak b}_-$ and ${\frak n}_+$, 
respectively,  
in the direct sum 
${\frak g}={\frak b}_- +{{\frak n}_+}$. By the definition of the induced action 
we have
$$
a\cdot s(s')=s(a^{-1}\cdot s'),\mbox{ for every }s'\in {{\frak b}_+}.
$$

On the other hand
$$
\begin{array}{l}
a\cdot s(s')=(P_{{\frak b}_-}({\rm Ad}(a)(s)),s')+\chi(P_{{\frak n}_+}({\rm 
Ad}(a)(s)))=\\
({\rm Ad}(a)(s),s')+({\rm Ad}(a)(s),f),
\end{array}
$$
where $f\in {\frak n}_-$ corresponds to $\chi$. Since the canonical bilinear 
form $(,)$ is Ad-invariant
the last formula may be rewritten as:
$$
a\cdot s(s')=(s,{\rm Ad}(a)^{-1}(s'+f))=s(P_{{\frak b}_+}({\rm 
Ad}(a)^{-1}(s'+f))),
$$
where $P_{{\frak b}_+}$ is the projector onto ${{\frak b}_+}$ in the direct sum 
${\frak g}={{\frak b}_+}+{\frak n}_-$.

Finally observe that the subspace $f+{{\frak b}_+}$ is stable under the adjoint 
action of $N_+$. 
Therefore $P_{{\frak b}_+}({\rm Ad}(a)^{-1}(s'+f))={\rm Ad}(a)^{-1}(s'+f)-f$, 
and the induced action of $N_+$ on ${{\frak b}_+}$ takes 
the form:
\begin{equation}\label{dualact}
a\cdot s'={\rm Ad}(a)(s'+f)-f.
\end{equation}

Now the algebra of invariants $S({\frak b}_-)^{N_+}$ may be identified with a 
certain subalgebra in the algebra of 
functions on the quotient ${{\frak b}_+}/{N_+}$. The space ${{\frak b}_+}/{N_+}$ 
has a nice geometric description.

Observe that $[f,{{\frak n}_+}]\subset {{\frak b}_+}$. Moreover,  $[f,{{\frak 
n}_+}]$ is an ${\rm ad}_\rho$ stable
subspace of ${\frak b}_+$. Since $\rho$ is a semi--simple element there exists 
an
${\rm ad}_\rho$ invariant stable subspace ${\frak s}\subseteq {{\frak b}_+}$ 
such that 
${{\frak b}_+}={\frak s}+[f,{{\frak n}_+}]$ is a direct sum. By Theorem 8 and 
Remark 19' in \cite{K1} $\frak s$ is an
$l$--dimensional subspace in ${{\frak n}_+}$. 
\vskip 0.3cm
\noindent
{\bf Theorem C (\cite{K}, Theorem 1.2 )}
{\em The map 
$$
{N_+}\times {\frak s}\rightarrow {{\frak b}_+}
$$
given by $(a,x)\mapsto a\cdot x$ is an isomorphism of affine varieties.
Therefore the quotient space ${{\frak b}_+}/{N_+}$ is isomorphic to ${\frak 
s}$.}
\vskip 0.3cm
The linear space $\frak s$ naturally appears in the study of regular elements in 
$\frak g$. Recall that an element of 
$\frak g$ is called regular if its centralizer in $\frak g$ is of minimal 
possible dimension. Let $R$ be 
the set of regular elements in $\frak g$. Clearly, $R$ is stable under the 
adjoint action of $G$ and in fact $R$
is the union of all $G$ orbits in $\frak g$ of maximal dimension.
\vskip 0.3cm
\noindent
{\bf Theorem D (\cite{K}, Theorem 1.1; \cite{K1}, Theorem 8)}
{\em The affine space $f+{\frak s}$ is contained in $R$ and is a cross--section 
for the action of $G$ on $R$. 
That is every $G$--orbit in $\frak g$ of maximal dimension intersects $f+{\frak 
s}$ in one and only one point.

Let $\widehat I_1,\ldots , \widehat I_l\in S({\frak g})^G$ be the fundamental 
invariants. 
$\widehat I_1,\ldots , \widehat I_l$ may be viewed as polynomial functions on 
${\frak g}^*
\cong {\frak g}$. The restrictions of these functions to $f+{\frak s}$ define a 
global coordinate system on 
${\frak s}$.}
\vskip 0.3cm
\noindent
{\bf Theorem E (\cite{K}, Theorem 1.3)} 
{\em For any $\widehat I\in S({\frak g})^G$ one has $\widehat I^\chi\in S({\frak 
b}_-)^{N_+}$.
Furthermore the map 
\begin{equation}\label{iso1}
S({\frak g})^G\rightarrow S({\frak b}_-)^{N_+},~~\widehat I\mapsto \widehat 
I^\chi
\end{equation}
is an algebra isomorphism. In particular 
$$
S({\frak b}_-)^N={\Bbb C}[\widehat I_1^\chi ,\ldots , \widehat I_l^\chi ]
$$
is a polynomial algebra in $l$ generators.}
\vskip 0.3cm
\noindent
{\em Proof.}
First observe that elements of $S({\frak g})$ may be viewed as polynomial 
functions on ${\frak g}^*
\cong {\frak g}$. Note also that the ideal $I_\chi^0$ is generated by the 
elements $x-(x,f),~x\in {{\frak n}_+}$. 
Therefore $I_\chi^0$ is the ideal of polynomial functions vanishing on the 
subspace $f+{{\frak b}_+}$ and so for every $\widehat I\in S({\frak g})$ 
$\widehat I^\chi$ may be 
regarded as the restriction of the function $\widehat I$ to the subspace 
$f+{{\frak b}_+}$. 

For $\widehat I\in S({\frak g})^G,~s'\in {{\frak b}_+}$ and $a\in N_+$ one has 
$\widehat I({\rm Ad}(a)(f+s'))=
\widehat I^\chi(a\cdot s')$. Since $\widehat I({\rm Ad}(a)(f+s'))=\widehat 
I(f+s')$ it follows that 
$\widehat I^\chi \in S({\frak b}_-)^{N_+}$.

By Theorem C the map 
$$
S({\frak b}_-)^{N_+} \rightarrow S({\frak s}^*) 
$$
given by the restriction $v\mapsto v|_{\frak s}$ is an algebra isomorphism.
Now by Theorem D the restrictions of the functions $\widehat I_i^\chi 
,~i=1,\ldots ,l$ to ${\frak s}$
are a coordinate system. Therefore (\ref{iso1}) is an isomorphism.
\vskip 0.3cm
\noindent
{\em Proof of Theorem B.}
First observe that elements $\widehat I_i^\chi,~i=1,\ldots ,l$
are the graded limits of 
the elements $I_i^\chi \in U({\frak b}_-)^{N_+},~i=1,\ldots ,l$. Therefore 
$GrW({\frak b}_-)=S({\frak b}_-)^{N_+}$. Recall that
$W({\frak b}_-)\subseteq U({\frak b}_-)^{N_+}$ is a linear subspace. 

Let $J\in U({\frak b}_-)^{N_+}\cap U_{k_1}({\frak b}_-)$ be an invariant 
element. 
Clearly, $GrJ\in S({\frak b}_-)^{N_+}$.
Since $GrW({\frak b}_-)=S({\frak b}_-)^{N_+}$ one can find elements 
$I_1\in W({\frak b}_-)\cap U_{k_1}({\frak b}_-)$ and
$J_1\in U({\frak b}_-)^{N_+}\cap U_{k_2}({\frak b}_-),
~~k_2<k_1$ such that
$$
J-I_1=J_1.
$$

Applying the same procedure to $J_1$ we obtain elements 
$J_2\in U({\frak b}_-)^{N_+}\cap U_{k_3}({\frak b}_-),
~k_3<k_2,~~I_2\in W({\frak b}_-)\cap U_{k_2}({\frak b}_-)$ such that
$$
J_1-I_2=J_2.
$$

We can continue this process. Since the standard filtration in $U({\frak b}_-)$ 
is
bounded below we finally obtain that for some $i~~ J_i-I_i=c\in {\Bbb C}$. By 
construction
the element $J$ is represented as $J=\sum_{j=1}^i I_j+c,~I_j\in W({\frak 
b}_-)\cap U_{k_j}({\frak b}_-)$.
Therefore $J\in W({\frak b}_-)$.  
This concludes the proof. 
 
Now we make an important remark.
\vskip 0.3cm
{\bf Remark A}
Observe that the space $U({\frak b}_-)^{N_+}$ may be interpreted as the zeroth 
cohomology space of the 
$U({{\frak n}_+})$ module $Y_\chi$, where $U({{\frak n}_+})$ is augmented by 
$\chi$. 
Indeed, for every associative algebra $B$ equipped with character $\chi$ and for 
every left 
$B$--module $V$ the cohomology module $H^*(V)$ is defined as the cohomology 
space of the complex (see \cite{carteil})

\begin{equation}\label{cohomol1}
{\rm Hom}_B(X,V),
\end{equation}
where $X$ is a projective resolution of the one--dimensional $B$--module ${\Bbb 
C}_\chi$ defined by 
$\chi$. In homological algebra $\chi$ is called an augmentation of $B$.
It is well--known that the graded vector space $H^*(V)$ does not depend on the 
resolution $X$ and
the zeroth cohomology
space $H^0(V)$ is isomorphic to the space of invariants ${\rm Hom}_B({\Bbb 
C}_\chi ,V)$ (see 
\cite{carteil}). Using the map
$$
{\rm Hom}_B({\Bbb C}_\chi ,V)\rightarrow V;~~ \hat{v}\mapsto \hat{v}(1)=v 
$$
this space may be identified with the subspace in $V$ spanned by elements
$v\in V$ such that $bv=\chi (b)v$ for every $b\in B$, i.e.
$$
H^0(V)={\rm Hom}_B({\Bbb C}_\chi ,V)=\{ v\in V: bv=\chi (b)v \mbox{ for every } 
b\in B \}.
$$

Now for $B=U({{\frak n}_+})$, $\chi$ as in Theorem B and $V=Y_\chi$ we have
$H^0(Y_\chi)=\{ v\in Y_\chi: xv=\chi (x)v \mbox{ for every } x\in U({{\frak 
n}_+}) \}$.
From (\ref{mainact}) and (\ref{iso2}) it follows that $H^0(Y_\chi)=U({\frak 
b}_-)^{N_+}$.

Now recall that by 
Theorem B there exists a linear isomorphism $W({\frak b}_-)\cong U({\frak 
b}_-)^{N_+}$. 
Therefore the associative algebra $W({\frak b}_-)$ is isomorphic to 
$H^0(Y_\chi)$ as a
linear space. In Section \ref{whitthom} we show that the multiplicative 
structure of $W({\frak b}_-)$ 
naturally appears in the context of homological algebra.

%%%%%%%%%%%%%%%%%%%%%%%%%%%%%%%%%%%%%%%%%%%%%%%%%%%%%%%%%%%%%%%%%%%%%%
%%%%%%%%%%%%%%%%%%%%%%%%%%%%%%%%%%%%%%%%%%%%%%%%%%%%%%%%%%%%%%%%%%%%%%%%
%%%%%%%%%%%%%%%%%%%%%%%%%%%%%%%%%%%%%%%%%%%%%%%%%%%%%%%%%%%%%%%%%%%%%

\chapter{Hecke algebras}\label{Hom}

In this section we give a homological definition of Hecke algebras (see Section 
\ref{Hecke}). 
Let $K$ be a ring with unit, $A$ an associative algebra over $K$, and $B$ a 
subalgebra of $A$ with
augmentation, that is, a $K$--algebra homomorphism $\varepsilon : B\rightarrow 
K$. 
The Hecke algebra $Hk^*(A,B,\varepsilon)$ of the triple $(A,B,\varepsilon )$ is 
a natural generalization of 
the algebra ${\rm Hom}_A(A\otimes_BK,A\otimes_BK)$. For 
every left $A$ module $V$ and every right $A$ module $W$ the algebra 
$Hk^*(A,B,\varepsilon)$ acts in both the cohomology 
space $H^{*}(B,V)$ and the homology space $H_*(B,W)$ of $V$ and $W$ as 
$B$--modules. 
Hecke algebras are
also closely related to the quantum BRST cohomology (see \cite{KSt}).

To define Hecke algebras we study complexes of $A$--endomorphisms of graded left 
$A$ modules.
Let $X$ be such a complex, ${\rm End}_{A}(X)$ be the corresponding complex of 
endomorphisms. Our main
observation is that the natural multiplication in ${\rm End}_{A}(X)$ given by 
composition of endomorphisms
induces a multiplicative structure on the cohomology space $H^*({\rm 
End}_{A}(X))$. Furthermore, the
associative algebra $H^*({\rm End}_{A}(X))$ only depends on the homotopy class 
of the complex $X$. 

As an application of our construction we show that the Whittaker model 
$W({\frak b}_-)$ is the zeroth graded component of the 
Hecke algebra of the triple $(U({\frak g}),U({\frak n}),\chi )$.

The exposition in this chapter follows \cite{S2}.
%%%%%%%%%%%%%%%%%%%%%%%%%%%%%%%%%%%%%%%%%%%%%%%%%%%%%%%%%%%%%%%%%

\section{Endomorphisms of complexes}\label{end}

\setcounter{equation}{0}
\setcounter{theorem}{0}

Let $A$ be an associative ring with unit, $X$ a graded complex of left 
$A$~modules equipped 
with a differential $d$ of degree $-1$. Recall the definition of the complex 
$Y={\rm End}_{A}(X)$ \cite{MacLane}.

By definition $Y$ is a $\Bbb Z$--graded complex
$$
Y=\bigoplus_{n=-\infty}^{\infty} Y^n
$$
with graded components defined as
$$
Y^n=\prod_{p+q=n} Y^{p,q},
$$
where 
$$
Y^{p,q}={\rm Hom}_{A}(X^{p},X^{-q}).
$$

Clearly $Y$ is
closed with respect to the multiplication given by composition of endomorphisms. 
Thus it is a graded 
associative algebra.

We introduce a differential on $Y$ of degree +1 as follows:
$$
\begin{array}{c}
({\bf d}f)^{p,q}=(-1)^{p+q}f^{p-1,q}\circ d +d\circ f^{p,q-1}, \\
f=\{ f^{p,q} \}, f^{p,q} \in Y^{p,q},
\end{array}
$$
where $d$ is the differential of $X$. If $f$ is homogeneous then
\begin{equation}\label{diff}
{\bf d}f=d\circ f -(-1)^{{\rm deg} (f)}f\circ d .
\end{equation}
So that ${\bf d}$ is the supercommutator by $d$.

We shall consider also the partial differentials $d'$ and $d''$ on Y :
\begin{equation}\label{part}
\begin{array}{cc}
\mbox{for } f\in Y^{p,q}& \\
(d'f)(x)=(-1)^{p+q+1}f(dx),& x\in X^{p+1}; \\
(d''f)(x)=df(x),& x\in X^{p}.
\end{array}
\end{equation}
It is easy to check that
$$
d'^2=d''^2=d'd''+d''d'=0
$$
These conditions ensure that ${\bf d}^2=0$.

The following property of ${\bf d}$ is crucial for the subsequent 
considerations.
\begin{lemma}
${\bf d}$ is a superderivation of $Y$.
\end{lemma}
{\em Proof.} Let $f$ and $g$ be homogeneous elements of $Y$. Then ${\rm deg} 
(fg) = {\rm deg} (f) + {\rm deg} (g)$ and 
(\ref{diff}) yields:
$$
\begin{array}{l}
{\bf d}(fg)=d\circ fg - (-1)^{{\rm deg} (f) +{\rm deg} (g)} fg \circ d= \\ 
d\circ fg -(-1)^{{\rm deg} (f)} f \circ d \circ g +(-1)^{{\rm deg} (f)} f \circ 
d \circ g - (-1)^{{\rm deg} (f) +{\rm deg} (g)} fg \circ d= \\
({\bf d}f)g+(-1)^{{\rm deg} (f)}f({\bf d}g). 
\end{array}
$$
This completes the proof.

The most important consequence of the lemma is
\begin{proposition}\label{alg}
The homology space $H^*(Y)$ inherits a multiplicative structure from Y. 
Thus $H^*(Y)$ is a graded associative algebra. 
\end{proposition}
{\em Proof.} First, the product of two cocycles is a cocycle. 
For if $f$ and $g$ are homogeneous and ${\bf d}f={\bf d}g=0$ then
$$
{\bf d}(fg)=({\bf d}f)g+(-1)^{{\rm deg} (f)}f({\bf d}g)=0.
$$

Now we have to show that the product of homology classes is well--defined. 
It suffices to verify that the product of a homogeneous cocycle with a 
homogeneous coboundary is cohomologous to zero.
For instance consider the product $f{\bf d}h$. Equation (\ref{diff}) gives
\begin{eqnarray}
f{\bf d}h=f\circ (d\circ h -(-1)^{{\rm deg} (h)}h\circ d)= \\
(-1)^{{\rm deg} (f)}d\circ f \circ h - (-1)^{{\rm deg} (h)}f\circ h\circ d = 
(-1)^{{\rm deg}(f)}{\bf d}(fh). \nonumber
\end{eqnarray}
This completes the proof.

One of the principal statements of homological algebra says that homotopically 
equivalent complexes have
the same homology. In particular the vector space $H^*(Y)$ depends only on the  
homotopy
class of the complex $X$. It turns out that the same is true for the algebraic 
structure of $H^*(Y)$.
Indeed we have the following
\begin{theorem}\label{equiv}
Let $X, X'$ be two homotopically equivalent graded complexes of left 
$A$--modules. Then
$$
 H^*(Y)\cong H^*(Y')
$$
as graded associative algebras.
\end{theorem}
{\em Proof.} Let $F:X \rightarrow X' , F':X' \rightarrow X$ be two maps between 
the complexes such that
$$
\begin{array}{lll}
F'F-{\rm id}_X=d_Xs+sd_X,& s:X\rightarrow X ,& s \in Y^{-1},\\
FF'-{\rm id}_{X'}=d_{X'}s'+s'd_{X'},& s':X'\rightarrow X' ,& s' \in Y'^{-1}.
\end{array}
$$

Consider the induced mappings of the complexes $Y,~~Y'$:
$$
\begin{array}{c}
FF'^{*}:Y \rightarrow Y' ,\\
FF'^{*}f= F \circ f \circ F' , f\in Y ;\\
F'F^{*}:Y' \rightarrow Y ,\\
F'F^{*}g= F' \circ g \circ F , f\in Y' .\\
\end{array}
$$
Their compositions are homotopic to the identity maps of $Y$ and $Y'$ (see 
Chap.~4, \cite{carteil} for a 
general statement about equivalences of functors). But this means that $FF'^{*}$ 
is inverse to $F'F^{*}$ 
when restricted to homology. Thus $H^*(Y)$ is isomorphic to $H^*(Y')$ as a 
vector space. We have to
show that the restrictions of $FF'^{*}$ and $F'F^{*}$ to the homologies are 
homomorphisms of algebras.

Let $f$ and $g$ be homogeneous elements of $Y$ and ${\bf d}_Xf={\bf d}_Xg=0$. By 
the definition of the
induced maps we have
$$
FF'^{*}(fg)=F \circ fg \circ F'.
$$

On the other hand
\begin{eqnarray}\label{hom}
FF'^{*}(f)FF'^{*}(g)=F \circ f \circ F'F \circ g \circ F'=   \\ 
F \circ f({\rm id}_X+d_Xs+sd_X)g \circ F'.\nonumber
\end{eqnarray}

Now recall that $f$ and $g$ are cocycles in $Y$. By (\ref{diff}) they 
supercommute with $d_X$:
\begin{equation}\label{cocycle}
d_X\circ f =(-1)^{{\rm deg} (f)}f\circ d_X.
\end{equation}

Using (\ref{cocycle}) and the fact that $F$ and $F'$ are morphisms of complexes 
we can rewrite (\ref{hom})
as follows:
\begin{eqnarray}\label{hom1}
F \circ f({\rm id}_X+d_Xs+sd_X)g \circ F'= F \circ fg \circ F' + \nonumber \\
+(-1)^{{\rm deg} (f)}d_{X'} \circ F \circ fsg \circ F' + 
(-1)^{{\rm deg} (g)}F \circ fsg \circ F' \circ d_{X'} = \\
=F \circ fg \circ F'+(-1)^{{\rm deg} (f)}{\bf d}_{X'}(F \circ fsg \circ F') . 
\nonumber
\end{eqnarray}

Finally observe that by (\ref{hom1}), $FF'^{*}(fg)$ and $FF'^{*}(f)FF'^{*}(g)$ 
belong to the same 
homology class in $H^{*}(Y')$. This completes the proof.

%%%%%%%%%%%%%%%%%%%%%%%%%%%%%%%%%%%%%%%%%%%%%%%%%%%%%%%%%%%%%%%%%

\section{Hecke algebras}\label{Hecke}

\setcounter{equation}{0}
\setcounter{theorem}{0}

Let $A$ be an associative algebra over a ring $K$ with unit, and $B$ a 
subalgebra of $A$ with
augmentation, that is, a $K$--algebra homomorphism $\varepsilon : B\rightarrow 
K$. 

Let $X$ be a projective resolution of the left $B$--module $K$ defined by 
$\varepsilon$. 
Since $X$ is a complex of left
$B$--modules, the space $A \otimes_B X$ is also a differential complex.   
Observe that this complex has the natural structure of a left $A$--module. 
Therefore we can apply Proposition \ref{alg} to define a 
graded associative 
algebra
$$
Hk^*(A,B,\varepsilon)=H^*({\rm End}_A(A \otimes_B X)).
$$

Note that all $B$--projective resolutions of $K$ are homotopically equivalent 
and so the complexes 
$A \otimes_B X$ are homotopically equivalent for different resolutions $X$. 
Hence 
by Theorem \ref{equiv} the associative algebra $Hk^*(A,B,\varepsilon)$ does not 
depend on the resolution $X$. 
We shall call it the {\it Hecke 
algebra} of the the triple $(A,B,\varepsilon )$.

Now consider $A$ as a left $A$--module and a right $B$--module via 
multiplication. In this way $A$ becomes
a left $A\otimes B^{opp}$--module. Let $X'$ be a projective resolution of this 
module. The complex $X' \otimes_B K$,
where the $B$ module structure on $K$ is defined by $\varepsilon$, is a left 
$A$--module. 
Therefore one can define an associative algebra
$$
\widehat{Hk}^*(A,B,\varepsilon)=H^*({\rm End}_A(X' \otimes_B K))
$$
independent of the resolution $X'$.
\begin{proposition}\label{iso}
$Hk^*(A,B,\varepsilon)$ is isomorphic to $\widehat{Hk}^*(A,B,\varepsilon)$ as a 
graded associative algebra.
\end{proposition}
{\em Proof.} We shall use the standard bar resolutions for computing 
$\widehat{Hk}^*(A,B,\varepsilon)$ and 
${Hk}^*(A,B,\varepsilon)$ \cite{MacLane}, \cite{carteil}. 
Consider the complex $B\otimes T(I(B)) \otimes B$, where $I(B)=B/K$ and $T$ 
denotes the tensor algebra 
of the vector space. Elements of $B\otimes T(I(B)) \otimes B$ are usually 
written as
$a[a_1,\ldots ,a_s]a'$. The differential is given by
\begin{eqnarray}
da[a_1,\ldots ,a_s]a'=aa_1[a_2,\ldots ,a_s]a'+ \\
\sum_{k=1}^{s-1}(-1)^{k}a[a_1,\ldots ,a_ka_{k+1},\ldots  ,a_s]a' + 
(-1)^sa[a_1,\ldots ,a_{s-1}]a_sa'. \nonumber
\end{eqnarray}

Then $B\otimes T(I(B)) \otimes B \otimes_B K = B\otimes T(I(B)) \otimes K$ is a 
free resolution of
the left $B$--module $K$. And $A \otimes_B B\otimes T(I(B)) \otimes B =A\otimes 
T(I(B)) \otimes B$ is 
a free resolution of $A$ as a right $B$--module. The complex $A\otimes T(I(B)) 
\otimes B$ is also a free
left $A$--module via left multiplication by elements of $A$. Hence this is an 
$A\otimes B^{opp}$--
free resolution of $A$.  

Thus the complex ${\rm End}_A(A\otimes_B B\otimes T(I(B)) \otimes K)={\rm 
End}_A(A\otimes T(I(B)) \otimes K)$ 
for the computation of $Hk^*(A,B,\varepsilon)$ is canonically isomorphic to the 
complex
${\rm End}_A(A \otimes T(I(B)) \otimes B \otimes_B K)={\rm End}_A(A \otimes 
T(I(B)) \otimes K)$ for the computation
of $\widehat{Hk}^*(A,B,\varepsilon)$. This establishes the isomorphism of the 
algebras.

%%%%%%%%%%%%%%%%%%%%%%%%%%%%%%%%%%%%%%%%%%%%%%%%%%%%%%%%%%%%%%%%%%%%%%

\section{Action in homology and cohomology spaces}

\setcounter{equation}{0}
\setcounter{theorem}{0}

Recall that for every left $B$--module $V$ the cohomology modules are defined to 
be
\begin{equation}\label{cohomol}
H^*(B,V)={\rm Ext}_B^*(K,V)=H^*({\rm Hom}_B(X,V)),
\end{equation}
where $X$ is a projective resolution of $K$. On the other hand for every right 
$B$--module $W$ one can define
the homology modules
\begin{equation}\label{homol}
H_*(B,W)={\rm Tor}_*^B(W,K)=H_*(W\otimes_B X).
\end{equation}

Now observe that for every left $A$--module $V$ 
the complex in (\ref{cohomol}) for calculating its cohomology as a right 
$B$--module 
may be represented as follows:
\begin{equation}\label{cohomcompl}
{\rm Hom}_B(X,V)={\rm Hom}_A(A\otimes_B X ,V).
\end{equation}

Endow the space ${\rm Hom}_A(A\otimes_B X ,V)$ with a right ${\rm 
End}_A(A\otimes_B X)$--action:
\begin{equation}\label{act1}
\begin{array}{ll}
{\rm Hom}_A(A\otimes_B X ,V) \times {\rm End}_A(A\otimes_B X) \rightarrow {\rm 
Hom}_A(A\otimes_B X ,V) , &\\
\varphi \times f \mapsto  \varphi \circ f, &\\
\varphi \in {\rm Hom}_A(A\otimes_B X ,V), f \in {\rm End}_A(A\otimes_B X).&
\end{array}
\end{equation}

This action is well--defined since $f$ commutes with the left $A$--action. 
Clearly this action respects the
gradings, i.e., it is an action of the graded associative algebra on the graded 
module.
\begin{proposition}\label{cohomact}
For every left $A$ module $V$ 
the action (\ref{act1}) gives rise to a right action
\begin{eqnarray}\label{cohomact1}
H^*(B,V)\times {Hk}^*(A,B,\varepsilon) \rightarrow H^*(B,V), \\
H^n(B,V)\times {Hk}^m(A,B,\varepsilon) \rightarrow H^{n+m}(B,V).\nonumber
\end{eqnarray}
\end{proposition}
{\em Proof.} Let $\varphi \in {\rm Hom}_A(A\otimes_B X ,V)$ and $d\varphi 
=\varphi \circ d =0$. Let also 
$f \in {\rm End}_A(A\otimes_B X)$ be a homogeneous cocycle. By (\ref{cocycle}) 
$\varphi \circ f$ is a cocycle
in ${\rm Hom}_A(A\otimes_B X ,V)$. Indeed
$$
d(\varphi \circ f)=\varphi \circ f \circ d =(-1)^{{\rm deg} (f)} \varphi \circ d 
\circ f =0.
$$

Next we need to show that the action does not depend on the choice of the 
representative $f$ in the
homology class $[f]$, that is $\varphi \circ {\bf d}g$ is homologous to zero for 
every homogeneous 
$g\in {\rm End}_A(A\otimes_B X)$. This is a direct consequence of the 
definitions:
$$
\varphi \circ {\bf d}g= \varphi \circ (d\circ g - (-1)^{{\rm deg}(g)} g\circ d)= 
-(-1)^{{\rm deg}(g)}d(\varphi \circ g),
$$
since $\varphi \circ d =0$. 

Finally let us check that the action is independent of the representative
in the homology class $[\varphi]$. For $\psi \in {\rm Hom}_A(A\otimes_B X ,V)$ 
$d\psi \circ f$ is always 
homologous to zero:
$$
d\psi \circ f= \psi \circ d \circ f = (-1)^{{\rm deg}(f)}\psi \circ f \circ d 
=(-1)^{{\rm deg}(f)}d(\psi \circ f).
$$
This concludes the proof.

Similarly for every right $A$--module $W$ one can equip the homology module 
$H_*(B,W)$ with a 
structure of a left ${Hk}^*(A,B,\varepsilon)$--module. First the complex 
$W\otimes_B X=W\otimes_A A\otimes_B X$ has
the natural structure of a left ${\rm End}_A(A\otimes_B X)$--module:
\begin{equation}\label{act2}
\begin{array}{ll}
{\rm End}_A(A\otimes_B X)\times W\otimes_A A\otimes_B X  \rightarrow W\otimes_A 
A\otimes_B X , &\\
f \times w\otimes x\mapsto  w\otimes f(x), &\\
w\otimes x \in W\otimes_A (A\otimes_B X), f \in {\rm End}_A(A\otimes_B X).&
\end{array}
\end{equation}

Observe that according to the convention of Section \ref{end} elements of ${\rm 
End}_A^n(A\otimes_B X)$ have
degree -n as operators in the graded space $W\otimes_A A\otimes_B X$:
$$
{\rm End}_A^n(A\otimes_B X)\times W\otimes_A A\otimes_B X_m \rightarrow 
W\otimes_A A\otimes_B X_{m-n}.
$$

The following assertion is an analogue of Proposition \ref{cohomact} for 
homology.
\begin{proposition}
For every right $A$ module $W$ the action (\ref{act2}) gives rise to a left 
action
\begin{eqnarray} 
Hk(A,B,\varepsilon)^* \times H_*(B,W) \rightarrow H_*(B,W),\\
Hk(A,B,\varepsilon)^n \times H_m(B,W) \rightarrow H_{m-n}(B,W).\nonumber
\end{eqnarray}
\end{proposition}

%%%%%%%%%%%%%%%%%%%%%%%%%%%%%%%%%%%%%%%%%%%%%%%%%%%%%%%%%%%%%%%%%%%%%%%%

\section{Structure of the Hecke algebras}

\setcounter{equation}{0}
\setcounter{theorem}{0}

In this section we investigate the Hecke algebras under some technical 
assumptions. The main theorem here is
\begin{theorem}\label{struct}
Assume that
$$
{\rm Tor}_n^B(A,K)=0 \mbox{  for } n>0.
$$
Then
$$
Hk^n(A,B,\varepsilon)={\rm Ext}^n_A(A\otimes_BK,A\otimes_BK)={\rm 
Ext}^n_B(K,A\otimes_BK).
$$
In particular
$$
Hk^n(A,B,\varepsilon)=0 ,~ n<0,
$$
and
$$
Hk^0(A,B,\varepsilon )={\rm Hom}_A(A\otimes_BK,A\otimes_BK)
$$
as an associative algebra.
\end{theorem}
{\em Proof.} Equip the complex $Y={\rm End}_A(A\otimes T(I(B)) \otimes K)$, 
which we used in Proposition \ref{iso} for the 
computation of $Hk^*(A,B,\varepsilon)$, with the first filtration as follows:
$$
F^kY=\sum_{n=-\infty}^{\infty}\prod_{p+q=n , p\geq k}Y^{p,q}.
$$

The associated graded complex with respect to this filtration is the double 
direct sum
$$
{\rm Gr}Y=\sum_{p,q=-\infty}^{\infty}Y^{p,q}.
$$
One can show that the filtration is regular and the second term of the 
corresponding spectral sequence is
\begin{equation}\label{spec}
E_2^{p,q}=H^p_{d'}(H^q_{d''}({\rm Gr}Y)),
\end{equation}
where $H^*_{d'}$ and $H^*_{d''}$ denote the homologies of the complex with 
respect to the partial differentials (\ref{part}).

Now observe that at the same time the complex $A\otimes T(I(B)) \otimes K$ is a 
complex 
for the calculation
of ${\rm Tor}_n^B(A,K)$ because $A\otimes T(I(B)) \otimes B$ is a free 
resolution of $A$ as a right $B$--module. It is
also free as a left $A$--module. Therefore the functor ${\rm Hom}_A(A\otimes 
T(I(B)) \otimes K, \cdot )$ is exact.
By assumption $H^*(A\otimes T(I(B)) \otimes K)={\rm Tor}_0^B(A,K)=A\otimes_BK$. 
Using the last two
observations we can calculate the cohomology of the complex ${\rm Gr}Y$ with 
respect to the differential $d''$ :
\begin{equation}\label{degener}
\begin{array}{l}
H^*_{d''}({\rm Gr}Y)=H^*_{d''}({\rm Hom}_A(A\otimes T(I(B)) \otimes K,A\otimes 
T(I(B)) \otimes K))=\\
{\rm Hom}_A(A\otimes T(I(B)) \otimes K,A\otimes_BK).
\end{array}
\end{equation}
Here ${\rm Hom}_A$ should be thought of as a direct sum of the double graded 
components. Now (\ref{degener}) provides that 
the spectral sequence (\ref{spec}) degenerates at the second term. Moreover,
$$
E_2^{p,*}=H^p_{d'}(H^0_{d''}({\rm Gr}Y))=H^p_{d'}({\rm Hom}_A (A\otimes T(I(B)) 
\otimes K,A\otimes_BK)).
$$

But the complex $A\otimes T(I(B)) \otimes K$ may be regarded as a free 
resolution of the left 
$A$--module $A\otimes_BK$. Therefore
$$
E_2^{p,*}={\rm Ext}^p_A(A\otimes_BK,A\otimes_BK).
$$

Finally by Theorem 5.12, \cite{carteil} we have:
$$
Hk^n(A,B,\varepsilon)=H^n(Y)=E_2^{n,0}={\rm Ext}^n_A(A\otimes_BK,A\otimes_BK).
$$

Since ${\rm Tor}_n^B(A,K)=0 \mbox{  for } n>0$ we can apply the Shapiro lemma 
(see Proposition 
4.1.3 in \cite{carteil}) to simplify the last expression:
$$
{\rm Ext}^n_A(A\otimes_BK,A\otimes_BK)={\rm Ext}^n_B(K,A\otimes_BK).
$$
Clearly, $Hk^0(A,B,\varepsilon)={\rm Hom}_A(A\otimes_BK,A\otimes_BK)$
as an associative algebra.
This completes the proof.

\begin{remark}\label{multhk}
In particular the conditions of the theorem are satisfied if $A$ is projective 
as a right $B$--module. For instance
suppose that there exists a subspace $N \subset A$ such that multiplication in 
$A$ provides an isomorphism of vector spaces $A \cong N\otimes B$. Then $A$ is a 
free right $B$--module.
\end{remark}

%%%%%%%%%%%%%%%%%%%%%%%%%%%%%%%%%%%%%%%%%%%%%%%%%%%%%%%%%%%%%%%%%%%%

\section{Comparison with the BRST complex}

\setcounter{equation}{0}
\setcounter{theorem}{0}

Let $\frak g$ be a Lie algebra over a field $K$. For simplicity we suppose that 
$\frak g$ is finite--dimensional. However 
the arguments presented below remain true, with some technical modifications, 
for an arbitrary Lie algebra.
We shall apply the construction of Section \ref{Hecke} in the following 
situation. 

Let $B=U({\frak g})$ and let $A$ be an associative 
algebra over $K$ containing $B$ as a subalgebra. Note that $U({\frak g})$ is 
naturally augmented.
Consider the $U({\frak g})$--free resolution of the left $U({\frak g})$--module 
$K$ as follows:
$$
\begin{array}{l}
X=U({\frak g})\otimes \Lambda ({\frak g}),\\
d (u\otimes x_1 \wedge \ldots \wedge x_n)= 
\sum_{i=1}^n (-1)^{i+1} ux_i\otimes x_1 \wedge \ldots \wedge \widehat{x_i} 
\wedge \ldots \wedge x_n +\\
\sum_{1\leq i< j \leq n}(-1)^{i+j} u\otimes [x_i,x_j] \wedge 
x_1 \wedge \ldots \wedge \widehat{x_i} \wedge \ldots \wedge \widehat{x_j} \wedge 
\ldots \wedge x_n,
\end{array}
$$
where the symbol $\widehat{x_i}$ indicates that $x_i$ is to be omitted.

Introduce operators of exterior and inner multiplication on $\Lambda ({\frak 
g})$ as follows. For every 
$x\in {\frak g}$ and $x^*\in {\frak g}^*$ we define 
$$
\begin{array}{l}
\overline{x} x_1 \wedge \ldots \wedge x_n= x\wedge x_1 \wedge \ldots \wedge 
x_n,\\
\\
\overline{x^*}x_1 \wedge \ldots \wedge x_n=
\sum_{i=1}^n (-1)^{i+1} x^*(x_i) x_1\wedge \ldots \wedge \widehat{x_i} \wedge 
\ldots \wedge x_n.
\end{array}
$$

Equip the linear space ${\frak g}+{\frak g}^*$ with a scalar product given by 
the canonical paring between 
${\frak g}$ and ${\frak g}^*$. Using this scalar product we can construct the 
Clifford algebra
$C({\frak g}+{\frak g}^*)$. The operators $\overline{x}, \overline{y^*}, x\in 
{\frak g},y^*\in{\frak g}^*$ satisfy
the defining relations of this algebra,
$$
\overline{x}\overline{y^*}+\overline{y^*}\overline{x}=y^*(x).
$$
Therefore the algebra $C({\frak g}+{\frak g}^*)$ naturally acts in the space 
$\Lambda ({\frak g})$. 
Moreover, it is well--known that ${\rm End}_K(\Lambda ({\frak g}))=C({\frak 
g}+{\frak g}^*)$.

Now the differential of the complex $A\otimes_{U({\frak g})} X = A\otimes 
\Lambda ({\frak g})$ may be explicitly
described using the operators of exterior and inner multiplications,
\begin{equation}\label{different}
d=\sum_i e_i\otimes \overline{e_i^*}-\sum_{i,j} 1\otimes 
\overline{[e_i,e_j]}\overline{e_i^*}\overline{e_j^*}.
\end{equation}
Here $e_i$ is a linear basis of ${\frak g}$, $e_i^*$ is the dual basis, 
$e_i\otimes 1$ is regarded as the 
operator of right multiplication in $A$, $e_i\otimes 1\cdot u\otimes 
1=ue_i\otimes 1$.

Now consider the complex ${\rm End}_A(A\otimes_{U({\frak g})} X)={\rm 
End}_A(A\otimes \Lambda ({\frak g}))$ for the 
computation of the algebra $Hk(A,B,\varepsilon)$.
Observe that
$$
{\rm End}_A(A\otimes \Lambda ({\frak g})) = A^{opp}\otimes {\rm End}_K(\Lambda 
({\frak g}))=A^{opp} \otimes C({\frak g}+{\frak g}^*).
$$
Under this 
identification $A^{opp}$ acts on $A\otimes \Lambda ({\frak g})$ by 
multiplication in $A$ on the right and the Clifford 
algebra acts by the exterior and inner multiplication in $\Lambda ({\frak g})$. 
This allows to consider the 
differential (\ref{different}) as an element of the complex $A^{opp} \otimes 
C({\frak g}+{\frak g}^*)$. 

It is easy to see that the canonical $\Bbb Z$--grading of the complex $A^{opp} 
\otimes C({\frak g}+{\frak g}^*)$ coincides 
mod 2 with the ${\Bbb Z}_2$--grading inherited from the Clifford algebra. 
Therefore according to (\ref{diff}) 
the differential $\bf d$ is given by the supercommutator in $A^{opp} \otimes 
C({\frak g}+{\frak g}^*)$ by element 
(\ref{different}).

Now recall that the complex $A^{opp} \otimes C({\frak g}+{\frak g}^*)$ with the 
differential given by 
the supercommutator by the element (\ref{different}) is the quantum BRST complex 
proposed in \cite{KSt}. 
This establishes
\begin{proposition}
The complex $({\rm End}_A(A\otimes_{U({\frak g})} X) , {\bf d})$ is isomorphic 
to the BRST one 
$A^{opp} \otimes C({\frak g}+{\frak g}^*)$ with the differential being the 
supercommutator by the element (\ref{different}).
\end{proposition}

%%%%%%%%%%%%%%%%%%%%%%%%%%%%%%%%%%%%%%%%%%%%%%%%%%%%%%%%%%%%%%%%%%%%%%%%%%%%

\section{Whittaker model as a Hecke algebra}\label{whitthom}

\setcounter{equation}{0}
\setcounter{theorem}{0}

In this section we use the notation introduced in Section \ref{notation}. Let 
$\frak g$ be a complex simple 
Lie algebra,
${\frak n}_+\subset {\frak g}$ the maximal nilpotent subalgebra, $\chi :{\frak 
n}_+\rightarrow {\Bbb C}$ a  
character. Let $W({\frak b}_-)$ be the Whittaker model of the center $Z({\frak 
g})$ of the 
universal enveloping algebra $U({\frak g})$.
\begin{proposition}\label{hkhom}
Suppose that the character $\chi$ is non--singular. Then $W({\frak b}_-)$ is 
isomorphic to
$Hk^0(U({\frak g}),U({\frak n}_+),\chi )^{opp}$ as an associative algebra.
\end{proposition}
{\em Proof.}
First observe that since ${\frak g}={\frak b}_-\oplus {\frak n}_+$ we have a 
linear 
isomorphism $U({\frak g})=U({\frak b}_-)\otimes U({\frak n}_+)$. Therefore from 
Remark 
\ref{multhk} and Theorem \ref{struct} it follows that $Hk^0(U({\frak 
g}),U({\frak n}_+),\chi )=
{\rm Hom}_{U({\frak g})}(Y_\chi ,Y_\chi )$, where $Y_\chi =U({\frak 
g})\otimes_{U({\frak n}_+)}{\Bbb C}_\chi$.

Now observe that the map
$$
{\rm Hom}_{U({\frak g})}(Y_\chi ,Y_\chi )\rightarrow {\rm Hom}_{U({\frak 
n}_+)}({\Bbb C}_\chi ,Y_\chi);~~
\tilde{v}\mapsto \hat{v}, 
$$
where $\hat{v}$ is given by $\hat{v}(z)=\tilde{v}(1\otimes z)$ for every $z\in 
{\Bbb C}_\chi$,
is a linear isomorphism.

Note also that by Remark A and Theorem B there exists a linear isomorphism 
$$
{\rm Hom}_{U({\frak n}_+)}({\Bbb C}_\chi ,Y_\chi)\rightarrow W({\frak 
b}_-);~~\hat{v} \mapsto v,
\mbox{ where }v\otimes 1=\hat{v}(1).
$$
Therefore we have a linear isomorphism
\begin{equation}\label{be}
{\rm Hom}_{U({\frak g})}(Y_\chi ,Y_\chi )\rightarrow W({\frak b}_-);~~
\tilde{v}\mapsto v.
\end{equation}
We have to prove that (\ref{be}) is an antihomomorphism.

Let $\tilde{v},\tilde{w}\in {\rm Hom}_{U({\frak g})}(Y_\chi ,Y_\chi )$ be two 
elements such that
$\tilde{v}(1\otimes 1)=v\otimes 1,~\tilde{w}(1\otimes 1)=w\otimes 1$. Then 
$\tilde{v}(\tilde{w}(1\otimes 1))=\tilde{v}(w\otimes 1)$. Since $\tilde{v}$ is 
an $U({\frak g})$ endomorphism
of $Y_\chi$ we have $\tilde{v}(\tilde{w}(1\otimes 1))=w\tilde{v}(1\otimes 
1)=wv\otimes 1$
This completes the proof.

%%%%%%%%%%%%%%%%%%%%%%%%%%%%%%%%%%%%%%%%%%%%%%%%%%%%%%%%%%%%%%%%%%%%%%%%%%
%%%%%%%%%%%%%%%%%%%%%%%%%%%%%%%%%%%%%%%%%%%%%%%%%%%%%%%%%%%%%%%%%%%%%%%%%%
%%%%%%%%%%%%%%%%%%%%%%%%%%%%%%%%%%%%%%%%%%%%%%%%%%%%%%%%%%%%%%%%%%%%%%%%%%%%

\chapter{Quantum deformation of the Whittaker model}\label{qWitt}

Let $\frak g$ be a complex simple Lie algebra, $U_h({\frak g})$ the standard 
quantum group associated with 
${\frak g}$.
In this section we construct a generalization of the Whittaker model $W({\frak 
b}_-)$ for $U_h({\frak g})$.

Let $U_h({\frak n}_+)$
be the subalgebra of $U_h({\frak g})$ corresponding to the nilpotent
Lie subalgebra ${\frak n}_+$. $U_h({\frak n}_+)$ is
generated by simple positive root generators of $U_h({\frak g})$
subject to the quantum Serre relations. It is easy to show that $U_h({\frak 
n}_+)$
has no non--singular characters (taking nonvanishing values
on all simple root generators). Our first main result
is a family of new realizations of the
quantum group $U_h({\frak g})$, one for each Coxeter element
in the corresponding Weyl group (see also \cite{S1}). The counterparts of 
$U({\frak n}_+)$,
which naturally arise in these new realizations of $U_h({\frak g})$,
do have non--singular characters. 

Using these new realizations we can immediately formulate a quantum group 
version
of Definition A. We also prove counterparts of Theorems A and B for $U_h({\frak 
g})$.

Finally we define quantum group generalizations of the Toda Hamiltonians. In the 
spirit of quantum harmonic
analysis these new Hamiltonians are difference operators. An alternative 
definition of these Hamiltonians has been
recently given in \cite{Et}. 

%%%%%%%%%%%%%%%%%%%%%%%%%%%%%%%%%%%%%%%%%%%%%%%%%%%%%%%%%%%

\section{Quantum groups}

\setcounter{equation}{0}
\setcounter{theorem}{0}

In this section we recall some basic facts about quantum groups. 
We follow the notation of \cite{ChP}. 

Let $h$ be an indeterminate, ${\Bbb C}[[h]]$ the ring of formal power series in 
$h$.
We shall consider ${\Bbb C}[[h]]$--modules equipped with the so--called 
$h$--adic 
topology. For every such module $V$ this topology is characterized by requiring 
that 
$\{ h^nV ~|~n\geq 0\}$ is a base of the neighbourhoods of $0$ in $V$, and that 
translations 
in $V$ are continuous. It is easy to see that, for modules equipped with this 
topology, every 
${\Bbb C}[[h]]$--module map is automatically continuous.

A topological Hopf algebra over ${\Bbb C}[[h]]$ is a complete ${\Bbb 
C}[[h]]$--module $A$
equipped with a structure of ${\Bbb C}[[h]]$--Hopf algebra (see \cite{ChP}, 
Definition 4.3.1),
the algebraic tensor products entering the axioms of the Hopf algebra are 
replaced by their 
completions in the $h$--adic topology.
We denote by $\mu , \imath , \Delta , \varepsilon , S$ the multiplication, the 
unit, the comultiplication,
the counit and the antipode of $A$, respectively.

The standard quantum group $U_h({\frak g})$ associated to a complex 
finite--dimensional simple Lie algebra
$\frak g$ is the algebra over ${\Bbb C}[[h]]$ topologically generated by 
elements
$H_i,~X_i^+,~X_i^-,~i=1,\ldots ,l$, and with the following defining relations:
\begin{equation}\label{qgrh}
\begin{array}{l}
[H_i,H_j]=0,~~ [H_i,X_j^\pm]=\pm a_{ij}X_j^\pm,\\
\\
X_i^+X_j^- -X_j^-X_i^+ = \delta _{i,j}{K_i -K_i^{-1} \over q_i -q_i^{-1}} , \\
\\
\mbox{where }K_i=e^{d_ihH_i},~~e^h=q,~~q_i=q^{d_i}=e^{d_ih},
\end{array}
\end{equation}
and the quantum Serre relations:
$$
\begin{array}{l}
\sum_{r=0}^{1-a_{ij}}(-1)^r 
\left[ \begin{array}{c} 1-a_{ij} \\ r \end{array} \right]_{q_i} 
(X_i^\pm )^{1-a_{ij}-r}X_j^\pm(X_i^\pm)^r =0 ,~ i \neq j ,\\ \\
\mbox{ where }\\
 \\
\left[ \begin{array}{c} m \\ n \end{array} \right]_q={[m]_q! \over 
[n]_q![n-m]_q!} ,~ 
[n]_q!=[n]_q\ldots [1]_q ,~ [n]_q={q^n - q^{-n} \over q-q^{-1} }.
\end{array}
$$
$U_h({\frak g})$ is a topological Hopf algebra over ${\Bbb C}[[h]]$ with 
comultiplication 
defined by 
$$
\begin{array}{l}
\Delta_h(H_i)=H_i\otimes 1+1\otimes H_i,\\
\\
\Delta_h(X_i^+)=X_i^+\otimes K_i+1\otimes X_i^+,
\end{array}
$$
$$
\Delta_h(X_i^-)=X_i^-\otimes 1 +K_i^{-1}\otimes X_i^-,
$$
antipode defined by
$$
S_h(H_i)=-H_i,~~S_h(X_i^+)=-X_i^+K_i^{-1},~~S_h(X_i^-)=-K_iX_i^-,
$$
and counit defined by
$$
\varepsilon_h(H_i)=\varepsilon_h(X_i^\pm)=0.
$$

We shall also use the weight--type generators defined by
$$
Y_i=\sum_{j=1}^l d_i(a^{-1})_{ij}H_j,
$$ 
and the elements $L_i=e^{hY_i}$. They commute with the root vectors $X_i^\pm$ as 
follows:
\begin{equation}\label{weight-root}
L_iX_j^\pm L_i^{-1}=q_i^{\pm \delta_{ij}}X_j^\pm .
\end{equation}

The Hopf algebra $U_h({\frak g})$ is a quantization of the standard bialgebra 
structure on $\frak g$, i.e. 
$U_h({\frak g})/hU_h({\frak g})=U({\frak g}),~~ \Delta_h=\Delta~(\mbox{mod }h)$, 
where $\Delta$ is 
the standard comultiplication on $U({\frak g})$, and 
$$
{\Delta_h -\Delta_h^{opp} \over h}~(\mbox{mod }h)=\delta ,
$$
where  
$\delta: {\frak g}\rightarrow {\frak g}\otimes {\frak g}$ is the standard 
cocycle on $\frak g$.
Recall that
$$
\delta (x)=({\rm ad}_x\otimes 1+1\otimes {\rm ad}_x)2r_+,~~ r_+\in {\frak 
g}\otimes {\frak g},
$$
\begin{equation}\label{rcl}
r_+=\frac 12 \sum_{i=1}^lY_i \otimes X_i + \sum_{\beta \in 
\Delta_+}(X_{\beta},X_{-\beta})^{-1} X_{\beta}\otimes X_{-\beta}.
\end{equation}
Here $X_{\pm \beta}\in {\frak g}_{\pm \beta}$ are root vectors of $\frak g$. 
The element $r_+\in {\frak g}\otimes {\frak g}$ is called a classical r--matrix.

The following proposition describes the algebraic structure of $U_h({\frak g})$. 
\begin{proposition}{\bf (\cite{ChP}, Proposition 6.5.5)}\label{algq}
Let $\frak g$ be a finite--dimensional complex simple Lie algebra, let 
$U_h({\frak h})$ be 
the subalgebra of $U_h({\frak g})$ topologically generated by the $H_i, 
i=1,\ldots l$.
Then, there is an isomorphism of algebras $\varphi :U_h({\frak g})\rightarrow 
U({\frak g})[[h]]$
over ${\Bbb C}[[h]]$ such that $\varphi =id$ (mod $h$) and $\varphi|_{U_h({\frak 
h})}=id$.
\end{proposition}

\begin{proposition}{\bf (\cite{ChP}, Proposition 6.5.7)}\label{zq}
If $\frak g$ is a finite--dimensional complex simple Lie algebra, the center 
$Z_h({\frak g})$ of 
$U_h({\frak g})$ is canonically isomorphic to $Z({\frak g})[[h]]$, where 
$Z({\frak g})$ is 
the center of $U({\frak g})$.
\end{proposition}

\begin{corollary}{\bf (\cite{ChP}, Corollary 6.5.6)}\label{rep}
If $\frak g$ be a finite--dimensional complex simple Lie algebra, then the 
assignment
$V\mapsto V[[h]]$ is a one--to--one correspondence between the 
finite--dimensional irreducible
representations of $\frak g$ and indecomposable representations of $U_h({\frak 
g})$
which are free and of finite rank as ${\Bbb C}[[h]]$--modules. Furthermore for
every such $V$ the action of the generators $H_i \in U_h({\frak g}),~~ 
i=1,\ldots l$ on 
$V[[h]]$ coincides with the action of the root generators $H_i \in {\frak h},~~ 
i=1,\ldots l$.
\end{corollary}

The representations of $U_h({\frak g})$ defined in the previous corollary are 
called
finite--dimensional representations of $U_h({\frak g})$. For every 
finite--dimensional representation
$\pi_V:{\frak g}\rightarrow {\rm End}V$ we denote the corresponding 
representation of
$U_h({\frak g})$ in the space $V[[h]]$ by the same letter.

$U_h({\frak g})$ is a quasitriangular Hopf algebra, i.e. there exists an 
invertible element
${\cal R}\in U_h({\frak g})\otimes U_h({\frak g})$, called a universal 
R--matrix, such that
\begin{equation}\label{quasitr}
\Delta^{opp}_h(a)={\cal R}\Delta_h(a){\cal R}^{-1}\mbox{ for all } a\in 
U_h({\frak g}),
\end{equation}
where $\Delta^{opp}=\sigma \Delta$, $\sigma$ is the permutation in $U_h({\frak 
g})^{\otimes 2}$,
$\sigma (x\otimes y)=y\otimes x$, and
\begin{equation}\label{rmprop}
\begin{array}{l}
(\Delta_h \otimes id){\cal R}={\cal R}_{13}{\cal R}_{23},\\
\\
(id \otimes \Delta_h){\cal R}={\cal R}_{13}{\cal R}_{12},
\end{array}
\end{equation}
where ${\cal R}_{12}={\cal R}\otimes 1,~{\cal R}_{23}=1\otimes {\cal R},
~{\cal R}_{13}=(\sigma \otimes id){\cal R}_{23}$.

From (\ref{quasitr}) and (\ref{rmprop}) it follows that $\cal R$ satisfies the 
quantum Yang--Baxter
equation:
\begin{equation}\label{YB}
{\cal R}_{12}{\cal R}_{13}{\cal R}_{23}={\cal R}_{23}{\cal R}_{13}{\cal R}_{12}.
\end{equation}

For every quasitriangular Hopf algebra we also have (see Proposition 4.2.7 in 
\cite{ChP}):
$$
(S\otimes id){\cal R}={\cal R}^{-1},
$$
and
\begin{equation}\label{S}
(S\otimes S){\cal R}={\cal R}.
\end{equation} 

We shall explicitly describe the element ${\cal R}$. 
First following \cite{kh-t} we recall the construction of root vectors of 
$U_h({\frak g})$. 
We shall use the so--called normal ordering in the root system  
$\Delta_+=\{\beta_1,\ldots ,\beta_N\}$ (see \cite{Z1}).
\begin{definition}\label{normord}
An ordering of the root system $\Delta_+$ is called normal if all simple roots 
are written in an arbitrary
order, and
for any theree roots $\alpha,~\beta,~\gamma$ such that
$\gamma=\alpha+\beta$ we have either $\alpha<\gamma<\beta$ or 
$\beta<\gamma<\alpha$.
\end{definition} 
To construct root vectors we shall apply the following  
inductive algorithm. Let $\alpha , \beta , \gamma \in \Delta_+$ be positive 
roots such that
$\gamma=\alpha+\beta,~\alpha<\beta$ and $[\alpha,\beta]$ is the minimal segment 
including
$\gamma$, i.e. the segment has no other roots $\alpha',\beta'$ such that 
$\gamma=\alpha'+\beta'$.
Suppose that $X_{\alpha}^\pm ,~X_{\beta}^\pm$ have
already been constructed. Then we define
\begin{equation}\label{rootvect}
\begin{array}{l}
X_{\gamma}^+=X_{\alpha}^+X_{\beta}^+ - 
q^{(\alpha,\beta)}X_{\beta}^+X_{\alpha}^+,\\
\\
X_{\gamma}^-= X_{\beta}^-X_{\alpha}^- - 
q^{-(\alpha,\beta)}X_{\alpha}^-X_{\beta}^-.
\end{array}
\end{equation} 

\begin{proposition}\label{rootprop}
For $\beta =\sum_{i=1}^lm_i\alpha_i,~m_i\in {\Bbb N}$  $X_{\beta}^\pm $ is a 
polynomial in 
the noncommutative variables $X_i^\pm$ homogeneous in each $X_i^\pm$ of degree 
$m_i$.
\end{proposition}

The root vectors $X_{\beta}$ satisfy the following relations:
$$
[X_\alpha^+,X_{\alpha}^-]=a(\alpha){e^{h\alpha^\vee}-e^{-h\alpha^\vee}\over 
q-q^{-1}}.
$$
where $a(\alpha)\in {\Bbb C}[[h]]$.
They commute with elements of the subalgebra $U_h({\frak h})$ as follows:

\begin{equation}\label{roots-cart}
[H_i,X_{\beta}^\pm]=\pm \beta(H_i)X_{\beta}^\pm,~i=1,\ldots ,l.
\end{equation}

Note that by construction
$$
\begin{array}{l}
X_\beta^+~(\mbox{mod }h)=X_\beta \in {\frak g}_\beta,\\
\\
X_\beta^-~(\mbox{mod }h)=X_{-\beta} \in {\frak g}_{-\beta}
\end{array}
$$
are root vectors of $\frak g$. This implies that $a(\alpha)~(\mbox{mod 
}h)=(X_\alpha,X_{-\alpha})$.

Let $U_h({\frak n}_+),U_h({\frak n}_-)$ be the ${\Bbb C}[[h]]$--subalgebras of 
$U_h({\frak g})$ topologically 
generated by the 
$X_i^+$ and by the $X_i^-$, respectively.

Now using the root vectors $X_{\beta}^\pm$ we can construct a topological basis 
of 
$U_h({\frak g})$.
Define for ${\bf r}=(r_1,\ldots ,r_N)\in {\Bbb N}^N$,
$$
(X^+)^{\bf r}=(X_{\beta_1}^+)^{r_1}\ldots (X_{\beta_N}^+)^{r_N},
$$
$$
(X^-)^{\bf r}=(X_{\beta_1}^-)^{r_1}\ldots (X_{\beta_N}^-)^{r_N},
$$
and for ${\bf s}=(s_1,\ldots s_l)\in {\Bbb N}^{~l}$,
$$
H^{\bf s}=H_1^{s_1}\ldots H_l^{s_l}.
$$
\begin{proposition}{\bf (\cite{kh-t}, Proposition 3.3)}\label{PBW}
The elements $(X^+)^{\bf r}$, $(X^-)^{\bf t}$ and $H^{\bf s}$, for ${\bf 
r},~{\bf t}\in {\Bbb N}^N$, 
${\bf s}\in {\Bbb N}^l$, form topological bases of $U_h({\frak n}_+),U_h({\frak 
n}_-)$ and $U_h({\frak h})$,
respectively, and the products $(X^+)^{\bf r}H^{\bf s}(X^-)^{\bf t}$ form a 
topological basis of 
$U_h({\frak g})$. In particular, multiplication defines an isomorphism of ${\Bbb 
C}[[h]]$ modules:
$$
U_h({\frak n}_-)\otimes U_h({\frak h}) \otimes U_h({\frak n}_+)\rightarrow 
U_h({\frak g}).
$$
\end{proposition}

An explicit expression for $\cal R$ may be written by making use of the 
q-exponential
$$
exp_q(x)=\sum_{k=0}^\infty {x^k \over (k)_q!},
$$ 
where
$$
(k)_q!=(1)_q\ldots (k)_q,~~(n)_q={q^n -1 \over q-1}.
$$

Now the element $\cal R$ may be written as (see Theorem 8.1 in \cite{kh-t}):
\begin{equation}\label{univr}
{\cal R}=exp\left[ h\sum_{i=1}^l(Y_i\otimes H_i)\right]\prod_{\beta}
exp_{q_{\beta}^{-1}}[(q-q^{-1})a(\beta)^{-1}X_{\beta}^+\otimes X_{\beta}^-],
\end{equation}
where $q_\beta =q^{(\beta,\beta)}$;
the product is over all the positive roots of $\frak g$, and the order of the 
terms is such that 
the $\alpha$--term appears to the left of the $\beta$--term if $\alpha <\beta$ 
with respect to the normal
ordering of $\Delta_+$.

\begin{remark}
The r--matrix $r_+=\frac 12 h^{-1}({\cal R}-1\otimes 1)~~(\mbox{mod }h)$, which 
is the classical limit of $\cal R$,
coincides with the classical r--matrix (\ref{rcl}).
\end{remark}

%%%%%%%%%%%%%%%%%%%%%%%%%%%%%%%%%%%%%%%%%%%%%%%%%%%%%%%%%%%%%%%%%%%%%%%%%%%%%%%%
%%%%%%%%%%%%%%%%%%%%%

\section{Non--singular characters and quantum groups}

\setcounter{equation}{0}
\setcounter{theorem}{0}

In this section we construct quantum counterparts of the principal nilpotent 
Lie subalgebras of complex simple Lie algebras and of their non--singular 
characters. We mainly follow the exposition presented in \cite{S1}.

First we would like to show that the algebra $U_h({\frak n}_+)$ spanned by 
$X_i^+ , i=1, \ldots , l$ does not admit characters
which take nonvanishing values on all generators $X_i^+$,
except for the case of $U_h(sl(2))$ when the quantum Serre relations
do not appear.
 
Suppose, $\chi_h$ is such a character, and 
$\chi_h(X_i^+)=c_i\in {\Bbb C}[[h]],~c_i\neq 0,~i=1,\ldots l$. 
By applying the character $\chi_h$ to the quantum Serre relations
one obtains a family of identities,
\begin{equation} \label{false}
\sum_{r=0}^{1-a_{ij}}(-1)^r 
\left[ \begin{array}{c} 1-a_{ij} \\ r \end{array} \right]_{q_i} =0 , 
\, i \neq j.
\end{equation}

We claim that some of these relations fail for  the quantized universal 
enveloping algebra $U_h({\frak g})$ 
of any simple Lie algebra $\frak g$ , with the exception of ${\frak g}=sl(2)$. 
In a more general setting, relations (\ref{false}) 
are analysed in the following lemma.
\begin{lemma}\label{qbinom} 
The only solutions of equation
\begin{equation}\label{c1}
\sum_{k=0}^{m}(-1)^k 
\left[ \begin{array}{c} m \\ k \end{array} \right]_{t}
t^{kc}=0 , 
\end{equation}
where $t$ is an indeterminate, are of the form
\begin{equation}\label{c2}
c=-m+1,-m+2, \ldots ,m-2,m-1.
\end{equation}
\end{lemma}
{\em Proof.}
According to the q--binomial theorem \cite{GR}, 
\begin{equation}\label{z}
\sum_{k=0}^{m}(-z)^k 
\left[ \begin{array}{c} m \\ k \end{array} \right]_{t} 
=\prod_{p=0}^{m-1}(1-t^{m-1-2p}z).
\end{equation}

Put $z=t^c$ in this relation. Then the l.h.s of (\ref{z}) coincides with 
the l.h.s. of (\ref{c1}). 
Now (\ref{z}) implies that $c=m-1-2p , p=0, \ldots ,m-1$
are the only solutions of (\ref{c1}).

Now we return to identities (\ref{false}).
Any Cartan matrix contains at least one off-diagonal element
equal to $-1$. Then, $m= 1- a_{ij} = 2$ and 
$c=\pm 1$, and Lemma \ref{qbinom} implies
that some of identities  (\ref{false}) are false for
any simple Lie algebra, except for $sl(2)$. Hence, subalgebras of $U_h({\frak 
g})$ generated by 
$X_i^+$ do not possess non--singular characters.

It is our goal to construct subalgebras of $U_h({\frak g})$ which resemble the 
subalgebra $U({\frak n}_+) 
\subset U({\frak g})$ and possess non--singular
characters. 
Denote by $S_l$ the symmetric group of $l$ elements.
To any element $\pi \in S_l$ we associate a Coxeter element $s_{\pi}$ by the 
formula
$s_\pi =s_{\pi (1)}\ldots s_{\pi (l)}$.
For each Coxeter element $s_\pi$ we define an associative algebra 
$U_h^{s_\pi}({\frak n}_+) $ generated by elements $e_i ,~ i=1, \ldots l$ subject 
to the relations :
\begin{equation}\label{fqpi}
\sum_{r=0}^{1-a_{ij}}(-1)^r q^{r c_{ij}^{\pi}}
\left[ \begin{array}{c} 1-a_{ij} \\ r \end{array} \right]_{q_i} 
(e_i )^{1-a_{ij}-r}e_j (e_i)^r =0 ,~ i \neq j , 
\end{equation}
where $c_{ij}^{\pi}=\left( {1+s_\pi \over 1-s_\pi }\alpha_i , \alpha_j \right)$ 
are matrix elements of the Caley transform of $s_\pi$ 
in the basis of simple roots.
\begin{proposition}\label{charf}
The map $\chi_h^{s_\pi}:U_h^{s_\pi}({\frak n}_+) \rightarrow {\Bbb C}[[h]]$ 
defined on generators by 
$\chi_h^{s_\pi}(e_i)=c_i,~c_i\in {\Bbb C}[[h]],~c_i\neq 0$ is a character of the 
algebra $U_h^{s_\pi}({\frak n}_+) $.
\end{proposition}

To show that  $\chi_h^{s_\pi}$ is a character of $U_h^{s_\pi}({\frak n}_+) $ it 
suffices to check that the  defining 
relations (\ref{fqpi}) belong to the kernel of $\chi_h^{s_\pi}$, i.e.
\begin{equation}\label{chifqpi}
\sum_{r=0}^{1-a_{ij}}(-1)^r q^{r c_{ij}^{\pi}}
\left[ \begin{array}{c} 1-a_{ij} \\ r \end{array} \right]_{q_i}=0 ,~ i \neq j . 
\end{equation}

As a preparation for the proof of Proposition \ref{charf} we study the matrix 
elements of the Caley transform of 
$s_\pi$ which enter the definition of $U_h^{s_\pi}({\frak n}_+) $.
\begin{lemma}\label{tmatrel}
The matrix elements of  ${1+s_\pi \over 1-s_\pi }$ are of the form :
\begin{equation}\label{matrel}
c_{ij}^{\pi}=\left( {1+s_\pi \over 1-s_\pi }\alpha_i , \alpha_j \right)=
\varepsilon_{ij}^\pi b_{ij},
\end{equation}
where
$$
\varepsilon_{ij}^\pi =\left\{ \begin{array}{ll}
-1 & \pi^{-1}(i) <\pi^{-1}(j) \\
0 & i=j \\
1 & \pi^{-1}(i) >\pi^{-1}(j) 
\end{array}
\right  .
$$
\end{lemma}
{\em Proof.} (compare \cite{Bur}, Ch. V, \S 6, Ex. 3).
First we calculate the matrix of the Coxeter element $s_\pi$ with respect to the 
basis of simple roots. We obtain this matrix in the form of the Gauss 
decomposition of the operator $s_\pi$.

Let $z_{\pi (i)}=s_\pi \alpha_{\pi (i)}$. Recall that 
$s_i(\alpha_j)=\alpha_j-a_{ji}\alpha_i$. 
Using this definition the elements $z_{\pi (i)}$ may be represented as:
$$
z_{\pi (i)}=y_{\pi (i)} -\sum_{k \geq i} a_{\pi (k) \pi (i)}y_{\pi (k)},
$$
where
\begin{equation}\label{y}
y_{\pi (i)}=s_{\pi (1)}\ldots s_{\pi (i-1)}\alpha_{\pi (i)}.
\end{equation}
Using the matrix notation we can rewrite the last formula as follows:
\begin{equation}\label{2*}
\begin{array}{l}
z_{\pi (i)}=
(I+V)_{\pi (k) \pi (i)}y_{\pi (k)} , \\  \\ \mbox{ where } V_{\pi (k) \pi (i)}=
\left\{ \begin{array}{ll}
a_{\pi (k) \pi (i)} & k\geq i \\
0 & k < i
\end{array}
\right  .
\end{array}
\end{equation}
 
To calculate the matrix of the operator $s_\pi$ with respect to the basis of 
simple roots we have to express 
the elements $y_{\pi (i)}$ via the simple roots.
Applying the definition of  simple reflections to (\ref{y}) we can pull out the 
element $\alpha_{\pi (i)}$ to the right:
\[
y_{\pi (i)}=\alpha_{\pi (i)}-\sum_{k<i}a_{\pi (k) \pi (i)}y_{\pi (k)}.
\]
Therefore
\[
\alpha_{\pi (i)}=(I+U)_{\pi (k) \pi (i)}y_{\pi (k)} ~, \mbox{ where } U_{\pi (k) 
\pi (i)}=
\left\{ \begin{array}{ll}
a_{\pi (k) \pi (i)} & k<i \\
0 & k \geq i
\end{array}
\right .
\]
Thus
\begin{equation}\label{1*}
y_{\pi (k)}=(I+U)^{-1}_{\pi (j) \pi (k)}\alpha_{\pi (j)}.
\end{equation}

Summarizing (\ref{1*}) and (\ref{2*}) we obtain:
\begin{equation}\label{**}
s_\pi \alpha_i=\left( (I+U)^{-1}(I-V) \right)_{ki}\alpha_k .
\end{equation}
This implies:
\begin{equation}\label{3*}
{1+s_\pi \over 1-s_\pi}\alpha_i=\left( {2I+U-V \over U+V}\right)_{ki}\alpha_k .
\end{equation}

Observe that $(U+V)_{ki}=a_{ki}$ and 
$(2I+U-V)_{ij}=-a_{ij}\varepsilon_{ij}^\pi$. 
Substituting these expressions into (\ref{3*}) we get :
\begin{eqnarray}
\left( {1+s_\pi \over 1-s_\pi }\alpha_i , \alpha_j \right) = 
-(a^{-1})_{kp}\varepsilon_{pi}^\pi a_{pi}b_{jk}=\\
-d_ja_{jk}(a^{-1})_{kp}\varepsilon_{pi}^\pi a_{pi}  = 
\varepsilon_{ij}^\pi b_{ij}.
\end{eqnarray}
This concludes the proof of the lemma.

\noindent
{\em Proof of Proposition \ref{charf} } Identities (\ref{chifqpi}) follow from 
Lemma \ref{qbinom} for $t=q_i,~~m=1-a_{ij},~~ c=\varepsilon_{ij}^\pi a_{ij}$ 
since set of solutions (\ref{c2}) always contains $\pm (m-1)$.

Motivated by relations (\ref{fqpi}) we suggest new realizations of the quantum 
group $U_h({\frak g})$,
one for each Coxeter element $s_\pi$.
Let 
$U_h^{s_\pi}({\frak g})$ be the associative algebra over ${\Bbb C}[[h]]$ with 
generators 
$e_i , f_i , H_i,~i=1, \ldots l$ subject to the relations:
\begin{equation}\label{sqgr}
\begin{array}{l}
[H_i,H_j]=0,~~ [H_i,e_j]=a_{ij}e_j, ~~ [H_i,f_j]=-a_{ij}f_j,\\
\\
e_i f_j -q^{ c^\pi _{ij}} f_j e_i = \delta _{i,j}{K_i -K_i^{-1} \over q_i 
-q_i^{-1}} , \\
 \\
K_i=e^{d_ihH_i}, \\
 \\
\sum_{r=0}^{1-a_{ij}}(-1)^r q^{r c_{ij}^\pi}
\left[ \begin{array}{c} 1-a_{ij} \\ r \end{array} \right]_{q_i} 
(e_i )^{1-a_{ij}-r}e_j (e_i)^r =0 ,~ i \neq j , \\
 \\
\sum_{r=0}^{1-a_{ij}}(-1)^r q^{r c_{ij}^\pi}
\left[ \begin{array}{c} 1-a_{ij} \\ r \end{array} \right]_{q_i} 
(f_i )^{1-a_{ij}-r}f_j (f_i)^r =0 ,~ i \neq j .
\end{array}
\end{equation}

It follows that the map $\tau_h^\pi :U_h^{s_\pi}({\frak n}_+) \rightarrow 
U_h^{s_\pi}({\frak g});~~~e_i\mapsto e_i$ is a {\em natural} embedding of 
$U_h^{s_\pi}({\frak n}_+) $ into $U_h^{s_\pi}({\frak g})$.  
From now on we identify $U_h^{s_\pi}({\frak n}_+) $ with the subalgebra in 
$U_h^{s_\pi}({\frak g})$ generated by $e_i ,i=1, \ldots l$.
\begin{theorem} \label{newreal}
For every solution $n_{ij}\in {\Bbb C},~i,j=1,\ldots ,l$ of equations 
\begin{equation}\label{eqpi}
d_jn_{ij}-d_in_{ji}=c^\pi_{ij}
\end{equation}
there exists an algebra
isomorphism $\psi_{\{ n\}} : U_h^{s_\pi}({\frak g}) \rightarrow 
U_h({\frak g})$ defined  by the formulas:
$$
\begin{array}{l}
\psi_{\{ n\}}(e_i)=X_i^+ \prod_{p=1}^lL_p^{n_{ip}},\\
 \\
\psi_{\{ n\}}(f_i)=\prod_{p=1}^lL_p^{-n_{ip}}X_i^- , \\
 \\
\psi_{\{ n\}}(H_i)=H_i .
\end{array}
$$
\end{theorem}
{\em Proof} is provided by direct verification of defining relations 
(\ref{sqgr}). The most nontrivial part is to verify deformed quantum Serre 
relations (\ref{fqpi}). 
The defining relations of $U_h({\frak g})$  imply the following relations for 
$\psi_{\{ n\}}(e_i)$,
$$
\sum_{k=0}^{1-a_{ij}}(-1)^k 
\left[ \begin{array}{c} 1-a_{ij} \\ k \end{array} \right]_{q_i}
q^{k({d_j}n_{ij}-d_in_{ji})}\psi_{\{ n\}}(e_i)^{1-a_{ij}-k}\psi_{\{ 
n\}}(e_j)\psi_{\{ n\}}(e_i)^k =0 ,
$$
for any $i\neq j$.
Now using equation (\ref{eqpi}) we arrive to relations (\ref{fqpi}).
\begin{remark} 
The general solution of equation (\ref{eqpi}) 
is given by
\begin{equation}\label{eq3}
n_{ji}=\frac 12 (\varepsilon_{ij}^\pi a_{ij} + \frac{s_{ij}}{d_i}),
\end{equation}
where $s_{ij}=s_{ji}$. 
\end{remark}

We call the algebra $U_h^{s_\pi}({\frak g})$ the Coxeter realization of the 
quantum group $U_h({\frak g})$ corresponding to the Coxeter element $s_\pi$.
\begin{remark}
Let $n_{ij}$ be a solution of the homogeneous system that corresponds to 
(\ref{eqpi}),
$$
d_in_{ji}-d_jn_{ij}=0.
$$
Then the map defined by
\begin{equation}
\begin{array}{l}
X_i^+ \mapsto X_i^+ \prod_{p=1}^lL_p^{n_{ip}},\\
 \\
X_i^- \mapsto \prod_{p=1}^lL_p^{-n_{ip}}X_i^- , \\
 \\
H_i \mapsto H_i 
\end{array}
\end{equation}
is an automorphism of $U_h({\frak g})$. Therefore for given Coxeter element the 
isomorphism $\psi_{\{ n\}}$ 
is defined uniquely up to automorphisms of $U_h({\frak g})$.
\end{remark}

Now we shall study the algebraic structure of $U_h^{s_\pi}({\frak g})$.
Denote by $U_h^{s_\pi}({\frak n}_-) $ the subalgebra in $U_h^{s_\pi}({\frak g})$ 
generated by
$f_i ,i=1, \ldots l$. From defining relations (\ref{sqgr}) it follows that the 
map 
$\overline \chi_h^{s_\pi}:U_h^{s_\pi}({\frak n}_-) \rightarrow {\Bbb C}[[h]]$ 
defined on generators by 
$\overline \chi_h^{s_\pi}(f_i)=c_i, c_i\in {\Bbb C}[[h]], c_i\neq 0$ is a 
character of the algebra $U_h^{s_\pi}({\frak n}_-)$.

Let $U_h^{s_\pi}({\frak h})$ be the subalgebra in $U_h^{s_\pi}({\frak g})$ 
generated by $H_i,~i=1,\ldots ,l$.
Define $U_h^{s_\pi}({\frak b}_\pm)=U_h^{s_\pi}({\frak n}_\pm)U_h^{s_\pi}({\frak 
h})$.

We shall construct a Poincar\'{e}--Birkhoff-Witt basis for $U_h^{s_\pi}({\frak 
g})$.
It is convenient to introduce an operator $K\in {\rm End}~{\frak h}$ such that
\begin{equation}\label{Kdef}
KH_i=\sum_{j=1}^l{n_{ij} \over d_i}Y_j.
\end{equation}
In particular, we have 
$$
{n_{ji} \over d_j}=(KH_j,H_i).
$$

Equation (\ref{eqpi}) is equivalent to the following equation for the operator 
$K$:
$$
K-K^* = {1+s_\pi \over 1-s_\pi}.
$$

\begin{proposition}\label{rootss}
(i)For any solution of equation (\ref{eqpi}) and any normal ordering of the root 
system $\Delta_+$
the elements $e_{\beta}=\psi_{\{ n\}}^{-1}(X_{\beta}^+e^{hK\beta^\vee})$ and 
$f_{\beta}=\psi_{\{ n\}}^{-1}(e^{-hK\beta^\vee}X_{\beta}^-),~\beta \in \Delta_+$
lie in the subalgebras $U_h^{s_\pi}({\frak n}_+)$ and $U_h^{s_\pi}({\frak 
n}_-)$, respectively.

(ii)Moreover, the elements 
$e^{\bf r}=e_{\beta_1}^{r_1}\ldots e_{\beta_N}^{r_N},~~f^{\bf 
t}=e_{\beta_1}^{t_1}\ldots e_{\beta_N}^{t_N}$
and $H^{\bf s}=H_1^{s_1}\ldots H_l^{s_l}$
for ${\bf r},~{\bf t},~{\bf s}\in {\Bbb N}^N$, form 
topological bases of $U_h^{s_\pi}({\frak n}_+),~U_h^{s_\pi}({\frak n}_-)$ and 
$U_h^{s_\pi}({\frak h})$, 
and the products $f^{\bf t}H^{\bf s}e^{\bf r}$ form a topological basis of 
$U_h^{s_\pi}({\frak g})$. In particular, multiplication defines an isomorphism 
of ${\Bbb C}[[h]]$ modules$$
U_h^{s_\pi}({\frak n}_-)\otimes U_h^{s_\pi}({\frak h})\otimes U_h^{s_\pi}({\frak 
n}_+)\rightarrow U_h^{s_\pi}({\frak g}).
$$
\end{proposition}
{\em Proof.} Let $\beta=\sum_{i=1}^l m_i\alpha_i \in \Delta_+$ be a positive 
root, 
$X_{\beta}^+\in U_h({\frak g})$ the corresponding root vector. Then 
$\beta^\vee=\sum_{i=1}^l m_id_iH_i$, and so
$K\beta^\vee=\sum_{i,j=1}^l m_in_{ij}Y_j$. Now the proof of the first statement 
follows immediately from 
Proposition \ref{rootprop}, commutation relations (\ref{weight-root}) and the 
definition of the isomorphism
$\psi_{\{ n\}}$. The second assertion is a consequence of Proposition \ref{PBW}.

Now we would like to choose a normal ordering of the root system $\Delta_+$ in 
such a way that 
$\chi_h^{s_\pi}(e_{\beta})=0$ and $\overline \chi_h^{s_\pi}(f_{\beta})=0$ if 
$\beta$ is not a simple root.  
\begin{proposition}\label{rootsh}
Choose a normal ordering of the root system $\Delta_+$ such that the simple 
roots are written
in the following order: $\alpha_{\pi (1)},\ldots ,\alpha_{\pi (l)}$.  
Then $\chi_h^{s_\pi}(e_{\beta})=0$ and $\overline \chi_h^{s_\pi}(f_{\beta})=0$ 
if $\beta$ is not a simple root.
\end{proposition} 
{\em Proof.} We shall consider the case of positive root generators.
The proof for 
negative root generators is similar to that for the positive ones.
 
The root vectors $X_{\beta}^+$ are defined in terms of iterated q-commutators
(see (\ref{rootvect})). Therefore it suffices to verify that for $i<j$
$$
\begin{array}{l}
\chi_h^{s_\pi}(e_{\alpha_{\pi(i)}+\alpha_{\pi(j)}})=\\
\\
\chi_h^{s_\pi}(\psi_{\{ n\}}^{-1}( (X_{\pi(i)}^+X_{\pi(j)}^+ - 
q^{(\alpha_{\pi(i)},\alpha_{\pi(j)})}X_{\pi(j)}^+X_{\pi(i)}^+)
e^{hK(d_{\pi(i)}H_{\pi(i)}+d_{\pi(j)}H_{\pi(j)})}))=0.
\end{array}
$$

From (\ref{Kdef}) and commutation relations (\ref{weight-root}) we obtain that
\begin{equation}\label{bebe}
\begin{array}{l}
\psi_{\{ n\}}^{-1}((X_{\pi(i)}^+X_{\pi(j)}^+ - 
q^{(\alpha_{\pi(i)},\alpha_{\pi(j)})}X_{\pi(j)}^+X_{\pi(i)}^+)
e^{hK(d_{\pi(i)}H_{\pi(i)}+d_{\pi(j)}H_{\pi(j)})})= \\
\\
q^{-d_{\pi(j)}n_{\pi(i)\pi(j)}}(e_{\pi(i)}e_{\pi(j)} - 
q^{b_{\pi(i)\pi(j)}+d_{\pi(j)}n_{\pi(i)\pi(j)}-d_{\pi(i)}n_{\pi(j)\pi(i)}}e_{\pi
(j)}e_{\pi(i)})
\end{array}
\end{equation}

Now using equation (\ref{eqpi}) and Lemma \ref{tmatrel} the combination 
$b_{\pi(i)\pi(j)}+d_{\pi(j)}n_{\pi(i)\pi(j)}-d_{\pi(i)}n_{\pi(j)\pi(i)}$ may be 
represented as
$b_{\pi(i)\pi(j)}+\varepsilon_{\pi(i)\pi(j)}^\pi b_{\pi(i)\pi(j)}$. But 
$\varepsilon_{\pi(i)\pi(j)}^\pi =-1$ for $i<j$ and therefore the r.h.s. of 
(\ref{bebe}) takes the form
$$
q^{-d_{\pi(j)}n_{\pi(i)\pi(j)}}[e_{\pi(i)},e_{\pi(j)}].
$$
Clearly,
$$
\chi_h^{s_\pi}(e_{\alpha_{\pi(i)}+\alpha_{\pi(j)}})=
q^{-d_{\pi(j)}n_{\pi(i)\pi(j)}}\chi_h^{s_\pi}([e_{\pi(i)},e_{\pi(j)}])=0.
$$

%%%%%%%%%%%%%%%%%%%%%%%%%%%%%%%%%%%%%%%%%%%%%%%%%%%%%%%%%%%%%%%%%%%%%%%%%%%%%%%%
%%%%%%%%%%%%%%%%%%%%%%%%%%%%

\section{Quantum deformation of the Whittaker model}

\setcounter{equation}{0}
\setcounter{theorem}{0}

In this section we define a quantum deformation of the Whittaker model $W({\frak 
b}_-)$.
Our construction is similar the one described in Section \ref{whitt}, the 
quantum group $U_h^{s_\pi}({\frak g})$,
the subalgebra $U_h^{s_\pi}({\frak n}_+)$ and characters 
${\chi_h^{s_\pi}}:U_h^{s_\pi}({\frak n}_+) \rightarrow {\Bbb C}[[h]]$ serve as 
natural
counterparts of the universal enveloping algebra $U({\frak g})$,
of the subalgebra $U({\frak n}_+)$ and of non--singular characters
$\chi:U({\frak n}_+) \rightarrow {\Bbb C}$.

Let ${U_h^{s_\pi}({\frak n}_+)}_{\chi_h^{s_\pi}}$ be the kernel of the character
$\chi_h^{s_\pi}:U_h^{s_\pi}({\frak n}_+) \rightarrow {\Bbb C}[[h]]$ so that
one has a direct sum
$$
U_h^{s_\pi}({\frak n}_+)={\Bbb C}[[h]]\oplus {U_h^{s_\pi}({\frak 
n}_+)}_{\chi_h^{s_\pi}}.
$$

From Proposition \ref{rootss} we have a linear 
isomorphism $U_h^{s_\pi}({\frak g})=U_h^{s_\pi}({\frak b}_-)\otimes 
U_h^{s_\pi}({\frak n}_+)$ and hence 
the direct sum
\begin{equation}\label{maindecq}
U_h^{s_\pi}({\frak g})=U_h^{s_\pi}({\frak b}_-) \oplus I_{\chi_h^{s_\pi}},
\end{equation}
where $I_{\chi_h^{s_\pi}}=U_h^{s_\pi}({\frak g}){U_h^{s_\pi}({\frak 
n}_+)}_{\chi_h^{s_\pi}}$ is the left--sided ideal generated by 
${U_h^{s_\pi}({\frak n}_+)}_{\chi_h^{s_\pi}}$.

For any $u\in U_h^{s_\pi}({\frak g})$ let $u^{\chi_h^{s_\pi}}\in 
U_h^{s_\pi}({\frak b}_-)$ be its component in 
$U_h^{s_\pi}({\frak b}_-)$ relative to the decomposition (\ref{maindecq}). 
Denote by $\rho_{\chi_h^{s_\pi}}$
the linear map 
$$
\rho_{\chi_h^{s_\pi}} : U_h^{s_\pi}({\frak g}) \rightarrow U_h^{s_\pi}({\frak 
b}_-)
$$
given by $\rho_{\chi_h^{s_\pi}} (u)=u^{\chi_h^{s_\pi}}$.

Denote by $Z_h^{s_\pi}({\frak g})$ the center of $U_h^{s_\pi}({\frak g})$.
From Proposition \ref{zq} and Theorem \ref{newreal} we obtain that 
$Z_h^{s_\pi}({\frak g})\cong Z({\frak g})[[h]]$. In particular, 
$Z_h^{s_\pi}({\frak g})$
is freely generated as a commutative topological algebra over ${\Bbb C}[[h]]$ by 
$l$ elements $I_1,\ldots , I_l$.

Let $W_h({\frak b}_-)=\rho_{\chi_h^{s_\pi}} (Z_h^{s_\pi}({\frak g}))$. 
\vskip 0.3cm
\noindent
{\bf Theorem $\bf A_h$ }
{\em The map
\begin{equation}\label{mapq}
\rho_{\chi_h^{s_\pi}} : Z_h^{s_\pi}({\frak g}) \rightarrow W_h({\frak b}_-)
\end{equation}
is an isomorphism of algebras. In particular, $W_h({\frak b}_-)$  
is freely generated as a commutative topological algebra over ${\Bbb C}[[h]]$ by 
$l$ elements 
$I_i^{\chi_h^{s_\pi}}=\rho_{\chi_h^{s_\pi}}(I_i),~~i=1,\ldots ,l$.}
\vskip 0.3cm
\noindent
{\em Proof} is similar to that of Theorem A in the classical case.
\vskip 0.3cm
\noindent
{\bf Definition $\bf A_h$ }
{\em The algebra $W_h({\frak b}_-)$ is called the Whittaker model of 
$Z_h^{s_\pi}({\frak g})$.}  
\vskip 0.3cm

Next we equip $U_h^{s_\pi}({\frak b}_-)$ with a structure of a left 
$U_h^{s_\pi}({\frak n}_+)$ module in such a
way that $W_h({\frak b}_-)$ is identified with the space of invariants with 
respect to this action.
Following Lemma A in the classical case we define this action by
\begin{equation}\label{mainactq}
x\cdot v =[x,v]^{\chi_h^{s_\pi}},
\end{equation}
where $v\in U_h^{s_\pi}({\frak b}_-)$ and $x\in U_h^{s_\pi}({\frak n}_+)$.

Consider the space $U_h^{s_\pi}({\frak b}_-)^{U_h^{s_\pi}({\frak n}_+)}$ of 
$U_h^{s_\pi}({\frak n}_+)$ invariants 
in $U_h^{s_\pi}({\frak b}_-)$ with respect to this
action. Clearly, 
$W_h({\frak b}_-)\subseteq U_h^{s_\pi}({\frak b}_-)^{U_h^{s_\pi}({\frak n}_+)}$. 
\vskip 0.3cm
\noindent
{\bf Theorem $\bf B_h$ }
{\em Suppose that $\chi_h^{s_\pi}(e_i)\neq 0~(\mbox{mod }h)$ for $i=1,\ldots l$. 
Then the space of $U_h^{s_\pi}({\frak n}_+)$ invariants 
in $U_h^{s_\pi}({\frak b}_-)$ with respect to the
action (\ref{mainactq}) is isomorphic to $W_h({\frak b}_-)$, i.e.}
\begin{equation}\label{invq}
U_h^{s_\pi}({\frak b}_-)^{U_h^{s_\pi}({\frak n}_+)}\cong W_h({\frak b}_-).
\end{equation}
\vskip 0.3cm
\noindent
{\em Proof.}
Let $p: U_h^{s_\pi}({\frak g})\rightarrow U_h^{s_\pi}({\frak 
g})/hU_h^{s_\pi}({\frak g})=U({\frak g})$ 
be the canonical projection. Note that $p(U_h^{s_\pi}({\frak n}_+))=U({\frak 
n}_+),~
p(U_h^{s_\pi}({\frak b}_-))=U({\frak b}_-)$ and for every $x\in 
U_h^{s_\pi}({\frak n}_+)$  
$\chi_h^{s_\pi}(x)~(\mbox{mod }h)=\chi(p(x))$ for some non--singular character 
$\chi: U({\frak n}_+)\rightarrow {\Bbb C}$. 
Therefore $p(\rho_{\chi_h^{s_\pi}}(x))=\rho_\chi(p(x))$ for every $x\in 
U_h^{s_\pi}({\frak g})$, 
and hence by Theorem ${\rm A}_q$ $p(W_h({\frak b}_-))=W({\frak b}_-)$. Using 
Lemma A and 
the definition of action (\ref{mainactq}) we also obtain that 
$p(U_h^{s_\pi}({\frak b}_-)^{U_h^{s_\pi}({\frak n}_+)})=
U({\frak b}_-)^{N_+} = W({\frak b}_-)$.

Now let $I\in U_h^{s_\pi}({\frak b}_-)^{U_h^{s_\pi}({\frak n}_+)}$ be an 
invariant element.
Then $p(I)\in W({\frak b}_-)$, and hence one can find an element $K_0\in 
W_h({\frak b}_-)$ such
that $I-K_0=hI_1,~I_1\in U_h^{s_\pi}({\frak b}_-)^{U_h^{s_\pi}({\frak n}_+)}$. 
Applying the same procedure
to $I_1$ one can find elements $K_1\in W_h({\frak b}_-),
~I_2\in U_h^{s_\pi}({\frak b}_-)^{U_h^{s_\pi}({\frak n}_+)}$ such that 
$I_1-K_1=hI_2$, i.e. 
$I-K_0-hK_1=0~(\mbox{mod }h^2)$. 
We can continue this process. Finally we obtain an infinite sequence of elements 
$K_i\in W_h({\frak b}_-)$ such that $I-\sum_{i=0}^p h^pK_p=0~(\mbox{mod 
}h^{p+1})$. Since the space 
$U_h^{s_\pi}({\frak b}_-)$ is complete in the $h$--adic topology the series 
$\sum_{i=0}^\infty h^pK_p\in W_h({\frak b}_-)$ 
converges to $I$. Therefore $I\in W_h({\frak b}_-)$. This completes the proof.

Similarly to Proposition \ref{hkhom} we have
\begin{proposition}
The algebra $W_h({\frak b}_-)$ is isomorphic to the zeroth graded component of 
the Hecke algebra
of the triple $(U_h^{s_\pi}({\frak g}),U_h^{s_\pi}({\frak n}_+),\chi_h^{s_\pi})$ 
with the opposite
multiplication,
$$
W_h({\frak b}_-)=Hk^0(U_h^{s_\pi}({\frak g}),U_h^{s_\pi}({\frak 
n}_+),\chi_h^{s_\pi})^{opp}.
$$
\end{proposition}

%%%%%%%%%%%%%%%%%%%%%%%%%%%%%%%%%%%%%%%%%%%%%%%%%%%%%%%%%%%%%%%%%%%%%%%%%%%%%%%%
%%%%%%%%%%%%%%%%%%%%%%%%%

\section{Coxeter realizations of quantum groups and Drinfeld twist}

\setcounter{equation}{0}
\setcounter{theorem}{0}

In this section we show that the Coxeter realizations $U_h^{s_\pi}({\frak g})$ 
of the quantum group $U_h({\frak g})$
are connected with quantizations of some nonstandard bialgebra structures on 
$\frak g$. At the quantum level 
changing bialgebra structure corresponds to the so--called Drinfeld twist. We 
shall consider a particular class
of such twists described in the following proposition.
\begin{proposition}{\bf (\cite{ChP}, Proposition 4.2.13)}\label{twdef}
Let $(A,\mu , \imath , \Delta , \varepsilon , S)$ be a Hopf algebra over a 
commutative ring. Let $\cal F$ be an invertible element of $A\otimes A$
such that 
\begin{equation}\label{twist}
\begin{array}{l}
{\cal F}_{12}(\Delta \otimes id)({\cal F})={\cal F}_{23}(id \otimes 
\Delta)({\cal F}),\\
\\
(\varepsilon \otimes id)({\cal F})=(id \otimes \varepsilon )({\cal F})=1.
\end{array}
\end{equation}
Then, $v=\mu (id\otimes S)({\cal F})$ is an invertible element of $A$ with
$$
v^{-1}=\mu (S\otimes id)({\cal F}^{-1}).
$$

Moreover , if we define $\Delta^{\cal F}:A\rightarrow A\otimes A$ and $S^{\cal 
F}:A\rightarrow A$ by
$$
\Delta^{\cal F}(a)={\cal F}\Delta(a){\cal F}^{-1},~~S^{\cal F}(a)=vS(a)v^{-1},
$$
then $(A,\mu , \imath , \Delta^{\cal F} , \varepsilon , S^{\cal F})$ is a Hopf 
algebra denoted by $A^{\cal F}$
and called the twist of $A$ by ${\cal F}$.
\end{proposition}

\begin{corollary}{\bf (\cite{ChP}, Corollary 4.2.15)}
Suppose that $A$ and ${\cal F}$ as in Proposition \ref{twdef}, but assume in 
addition that $A$ is quasitriangular
with universal R--matrix $\cal R$. Then $A^{\cal F}$ is quasitriangular with 
universal R--matrix
\begin{equation}\label{rf}
{\cal R}^{\cal F}={\cal F}_{21}{\cal R}{\cal F}^{-1},
\end{equation}
where ${\cal F}_{21}=\sigma {\cal F}$.
\end{corollary}

Fix a Coxeter element $s_\pi\in W$, $s_\pi=s_{\pi (1)}\ldots s_{\pi (l)}$.
Consider the twist of the Hopf algebra $U_h({\frak g})$ by the element
\begin{equation}\label{Ftw}
{\cal F}=exp(-h\sum_{i,j=1}^l {n_{ji} \over d_j}Y_i\otimes Y_j) \in U_h({\frak 
h})\otimes U_h({\frak h}),
\end{equation}
where $n_{ij}$ is a solution of the corresponding equation (\ref{eqpi}).

This element satisfies conditions (\ref{twist}), and so $U_h({\frak g})^{\cal 
F}$ is a quasitriangular 
Hopf algebra with the universal R--matrix ${\cal R}^{\cal F}={\cal F}_{21}{\cal 
R}{\cal F}^{-1}$, 
where $\cal R$ is given by (\ref{univr}). We shall explicitly calculate the 
element ${\cal R}^{\cal F}$.
Substituting (\ref{univr}) and (\ref{Ftw}) into (\ref{rf}) and using 
(\ref{roots-cart}) we obtain
$$
\begin{array}{l}
{\cal R}^{\cal F}=exp\left[ h(\sum_{i=1}^l(Y_i\otimes H_i)+
\sum_{i,j=1}^l (-{n_{ij} \over d_i}+{n_{ji} \over d_j})Y_i\otimes Y_j) 
\right]\times \\
\prod_{\beta}
exp_{q_{\beta}^{-1}}[(q-q^{-1})a(\beta)^{-1}X_{\beta}^+e^{hK\beta^\vee} \otimes 
e^{-hK^*\beta^\vee}X_{\beta}^-],
\end{array}
$$
where $K$ is defined by (\ref{Kdef}).

Equip $U_h^{s_\pi}({\frak g})$ with the comultiplication given by :
$\Delta_{s_\pi}(x)=(\psi_{\{ n\}}^{-1}\otimes \psi_{\{ n\}}^{-1})\Delta_h^{\cal 
F}(\psi_{\{ n\}}(x))$.
Then $U_h^{s_\pi}({\frak g})$ becomes a quasitriangular Hopf algebra with the 
universal R--matrix
${\cal R}^{s_\pi}=\psi_{\{ n\}}^{-1}\otimes \psi_{\{ n\}}^{-1}{\cal R}^{\cal 
F}$. Using equation
(\ref{eqpi}) and Lemma \ref{tmatrel} this R--matrix may be written as follows
\begin{equation}\label{rmatrspi}
\begin{array}{l}
{\cal R}^{s_\pi}=exp\left[ h(\sum_{i=1}^l(Y_i\otimes H_i)+
\sum_{i=1}^l {1+s_\pi \over 1-s_\pi }H_i\otimes Y_i) \right]\times \\
\prod_{\beta}
exp_{q_{\beta}^{-1}}[(q-q^{-1})a(\beta)^{-1}e_{\beta} \otimes 
e^{h{1+s_\pi \over 1-s_\pi} \beta^\vee}f_{\beta}].
\end{array}
\end{equation}

The element ${\cal R}^{s_\pi}$ may be also represented in the form
\begin{equation}\label{rmatrspi'}
\begin{array}{l}
{\cal R}^{s_\pi}=exp\left[ h(\sum_{i=1}^l(Y_i\otimes H_i)\right]\times \\
\prod_{\beta}
exp_{q_{\beta}^{-1}}[(q-q^{-1})a(\beta)^{-1}e_{\beta}e^{-h{1+s_\pi \over 
1-s_\pi}\beta^\vee}\otimes f_\beta]
exp\left[ h(\sum_{i=1}^l {1+s_\pi \over 1-s_\pi }H_i\otimes Y_i)\right] .
\end{array}
\end{equation}

The comultiplication $\Delta_{s_\pi}$ is given on generators by
$$
\begin{array}{l}
\Delta_{s_\pi}(H_i)=H_i\otimes 1+1\otimes H_i,\\
\\
\Delta_{s_\pi}(e_i)=e_i\otimes e^{hd_i{2 \over 1-s_\pi}H_i}+1\otimes e_i,\\
\\
\Delta_{s_\pi}(f_i)=f_i\otimes e^{-hd_i{1+s_\pi \over 
1-s_\pi}H_i}+e^{-hd_iH_i}\otimes f_i.
\end{array}
$$

Note that the Hopf algebra $U_h^{s_\pi}({\frak g})$ is a quantization of the 
bialgebra structure on $\frak g$
defined by the cocycle
\begin{equation}\label{cocycles}
\delta (x)=({\rm ad}_x\otimes 1+1\otimes {\rm ad}_x)2r^{s_\pi}_+,~~ 
r^{s_\pi}_+\in {\frak g}\otimes {\frak g},
\end{equation}
where $r^{s_\pi}_+=r_+ + \frac 12 \sum_{i=1}^l {1+s_\pi \over 1-s_\pi 
}H_i\otimes Y_i$, and $r_+$ is given by (\ref{rcl}).
 
We shall also need the following property of the antipode $S^{s_\pi}$ of 
$U_h^{s_\pi}({\frak g})$.
\begin{proposition}\label{sqant}
The square of the antipode $S^{s_\pi}$ is an inner automorphism of 
$U_h^{s_\pi}({\frak g})$ given by
$$
(S^{s_\pi})^2(x)=e^{2h\rho^\vee}xe^{-2h\rho^\vee},
$$
where $\rho^\vee=\sum_{i=1}^lY_i$.
\end{proposition}
{\em Proof.} 
First observe that by Proposition \ref{twdef} the antipode of 
$U_h^{s_\pi}({\frak g})$ has the form:
$S^{s_\pi}(x)=\psi_{\{ n\}}^{-1}(vS_h(\psi_{\{ n\}}(x))v^{-1})$, where 
$$
v=exp(h\sum_{i,j=1}^l {n_{ji} \over d_j}Y_iY_j).
$$
Therefore $(S^{s_\pi})^2(x)=\psi_{\{ n\}}^{-1}(vS_h(v^{-1})S_h^2(\psi_{\{ 
n\}}(x))S_h(v)v^{-1})$.
Note that $S_h(v)=v$, and hence $(S^{s_\pi})^2(x)=\psi_{\{ 
n\}}^{-1}(S_h^2(\psi_{\{ n\}}(x)))$.

Finally observe that from explicit formulas for the antipode of $U_h({\frak g})$ 
it follows that 
$S_h^2(x)=e^{2h\rho^\vee}xe^{-2h\rho^\vee}$. This completes the proof.

In conclusion we note that using Corollary \ref{rep} and the isomorphism 
$\psi_{\{ n\}}$ one can define finite--dimensional representations of 
$U_h^{s_\pi}({\frak g})$.

%%%%%%%%%%%%%%%%%%%%%%%%%%%%%%%%%%%%%%%%%%%%%%%%%%%%%%%%%%%%%%%%%%%%%%%%%%%%%%%%
%%%%%%%%%%%%%%%%%%%%%%

\section{Quantum deformation of the Toda lattice}\label{toda}

\setcounter{equation}{0}
\setcounter{theorem}{0}

Recall that one of
the main applications of the algebra $W({\frak b}_-)$ is the quantum Toda 
lattice \cite{K'}.
Let $\overline \chi : {\frak n}_- \rightarrow {\Bbb C}$ be a non--singular 
character of
the opposite nilpotent subalgebra ${\frak n}_-$. We denote the character of 
$N_-$ corresponding to 
$\overline \chi$ by the same letter. The algebra $U({\frak b}_-)$ naturally acts 
by differential
operators in the space $C^\infty ({\Bbb C}_{\overline \chi}\otimes_{N_-}{B_-})$. 
This space may be
identified with $C^\infty (H)$.
Let $D_1,\ldots ,D_l$ be the differential operators on $C^\infty (H)$ which 
correspond to the elements 
$I_1^\chi ,\ldots , I_l^\chi\in W({\frak b}_-)$. Denote by $\varphi$ the 
operator of multiplication in 
$C^\infty (H)$ by the function $\varphi (e^h)=e^{\rho(h)}$, where $h\in {\frak 
h}$. The operators $M_i=\varphi D_i\varphi^{-1}, i=1,\ldots l$
are called the quantum Toda Hamiltonians. Clearly, they commute with each other.

In particular if $I$ is the quadratic Casimir element then the corresponding 
operator $M$ is 
the well--known second--order differential operator:
$$
M=\sum_{i=1}^l \partial_i^2
+\sum_{i=1}^l \chi(X_{\alpha_i})\overline 
\chi(X_{-\alpha_i})e^{-\alpha_i(h)}+(\rho,\rho),
$$
where $\partial_i={\partial \over \partial y_i}$, and $y_i,~i=1,\ldots l$ is an 
ortonormal basis of $\frak h$.

Using the algebra $W_h({\frak b}_-)$ we shall construct quantum group analogues 
of the Toda Hamiltonians.
A slightly different approach has been recently proposed in \cite{Et}.

Denote by $A$ the space of linear functions on 
${\Bbb C}[[h]]_{\overline \chi_h^{s_\pi}}\otimes_{U_h^{s_\pi}({\frak 
n}_-)}U_h^{s_\pi}({\frak b}_-)$, where
${\Bbb C}[[h]]_{\overline \chi_h^{s_\pi}}$ is the one--dimensional 
$U_h^{s_\pi}({\frak n}_-)$ module
defined by ${\overline \chi_h^{s_\pi}}$.
Note that 
${\Bbb C}[[h]]_{\overline \chi_h^{s_\pi}}\otimes_{U_h^{s_\pi}({\frak 
n}_-)}U_h^{s_\pi}({\frak b}_-)\cong U_h^{s_\pi}({\frak h})$
as a linear space. Therefore $A=U_h^{s_\pi}({\frak h})^*$.
The algebra $U_h^{s_\pi}({\frak b}_-)$ naturally acts on 
${\Bbb C}[[h]]_{\overline \chi_h^{s_\pi}}\otimes_{U_h^{s_\pi}({\frak 
n}_-)}U_h^{s_\pi}({\frak b}_-)$ 
by multiplications from the right. This action induces an $U_h^{s_\pi}({\frak 
b}_-)$--action in the space $A$.
We denote this action by $L$, $L:U_h^{s_\pi}({\frak b}_-)\rightarrow {\rm 
End}A$. Clearly, this action generates 
an action of the algebra $W_h({\frak b}_-)$ on $A$.

To construct deformed Toda Hamiltonians we shall use certain elements in  
$W_h({\frak b}_-)$.
These elements may be described as follows.
Let $\mu : U_h^{s_\pi}({\frak g}) \rightarrow {\Bbb C}[[h]]$ be a map such that 
$\mu(uv)=\mu(vu)$. By Proposition 
\ref{sqant} $(S^{s_\pi})^2(x)=e^{2h\rho^\vee}xe^{-2h\rho^\vee}$. Hence from 
Remark 1 in \cite{D} it follows that  
$(id\otimes \mu)({\cal R}_{21}^{s_\pi}{\cal R}^{s_\pi}(1\otimes 
e^{2h\rho^\vee}))$, where
${\cal R}_{21}^{s_\pi}=\sigma {\cal R}^{s_\pi}$, is a central element.
In particular, for any finite--dimensional $\frak g$--module $V$ the element
\begin{equation}\label{centrelv}
C_V=(id\otimes tr_V)({\cal R}_{21}^{s_\pi}{\cal R}^{s_\pi}(1\otimes 
e^{2h\rho^\vee})),
\end{equation}
where $tr_V$ is the 
trace in $V[[h]]$, is central in $U_h^{s_\pi}({\frak g})$.

Using formulas (\ref{rmatrspi}) and (\ref{rmatrspi'}) we can easily compute 
elements 
$\rho_{\chi_h^{s_\pi}}(C_V)\in W_h({\frak b}_-)$.
For every finite--dimensional $\frak g$--module $V$ we have
\begin{equation}\label{todah}
\begin{array}{l}
\rho_{\chi_h^{s_\pi}}(C_V)=(id\otimes tr_V)( e^{t_0}\prod_{\beta}
exp_{q_{\beta}^{-1}}[(q-q^{-1})a(\beta)^{-1}f_\beta \otimes 
e_{\beta}e^{-h{1+s_\pi \over 1-s_\pi}\beta^\vee}]\times \\
\\
e^{t_0}\prod_{\beta}
exp_{q_{\beta}^{-1}}[(q-q^{-1})a(\beta)^{-1}{\chi_h^{s_\pi}}(e_{\beta}) \otimes 
e^{h{1+s_\pi \over 1-s_\pi} \beta^\vee}f_{\beta}](1\otimes e^{2h\rho^\vee})),
\end{array}
\end{equation}
where $t_0=h\sum_{i=1}^l(Y_i\otimes H_i)$.

We denote by $W_h^{Rep}({\frak b}_-)$ the subalgebra in $W_h({\frak b}_-)$ 
generated by the elements
$\rho_{\chi_h^{s_\pi}}(C_V)$, where $V$ runs through all finite--dimensional 
representations of $\frak g$.
Note that for every finite--dimensional $\frak g$--module $V$ 
$\rho_{\chi_h^{s_\pi}}(C_V)$ is a 
polynomial in noncommutative elements $f_i,~e^{hx},~x\in {\frak h}$.

Now we shall realize elements of $W_h^{Rep}({\frak b}_-)$ as difference 
operators.
Let $H_h\in U_h^{s_\pi}({\frak h})$ be the subgroup generated by elements 
$e^{hx},~x\in {\frak h}$.
A difference operator on $A$ is an operator $T$ of the form
$T=\sum f_iT_{x_i}$ (a finite sum), where $f_i \in A$, and for every $y\in H_h~$ 
$T_{x}f(y)=(ye^{hx}),~x\in {\frak h}$.
\begin{proposition}{\bf (\cite{Et}, Proposition 3.2)}\label{diffh}
For any $Y \in U_h^{s_\pi}({\frak b}_-)$, 
which is a polynomial in noncommutative elements $f_i,~e^{hx},~x\in {\frak h}$,
the operator $L(Y)$ is a difference operator on $A$.
In particular, the operators $L(I),~I\in W_h^{Rep}({\frak b}_-)$ are mutually 
commuting 
difference operators on $A$.
\end{proposition}
{\em Proof.} It suffices to verify that $L(f_i)$ are difference operators on 
$H_h$.
Indeed,
$$
L(f_i)f(e^{hx})=f(e^{hx}f_i)=e^{-h\alpha_i(x)}f(f_ie^{hx})=\overline 
\chi_h^{s_\pi}(f_i)e^{-h\alpha_i(x)}f(e^{hx}).
$$
This completes the proof.

Let $\jmath : H_h\rightarrow U_h^{s_\pi}({\frak h})$ be the canonical embedding. 
Denote $A_h=\jmath^*(A)$.
Let $T$ be a difference operator on $A$. Then one can define a difference 
operator $\jmath^*(T)$ on the space
$A_h$ by $\jmath^*(T)f(y)=T(\jmath(y))$.  

Let $D_i^h=\jmath^*(L(\rho_{\chi_h^{s_\pi}}(C_{V_i})))$, where $V_i,~i=1,\ldots 
l$ are the fundamental representations of
$\frak g$. 
Denote by $\varphi_h$ the operator of multiplication in 
$A_h$ by the function $\varphi_h (e^{hx})=e^{h\rho(x)}$, where $x\in {\frak h}$. 
The operators $M_i^h=\varphi_h D_i^h\varphi^{-1}_h, i=1,\ldots l$
are called the quantum deformed Toda Hamiltonians.

From now on we suppose that $\pi=id$ and that the ordering of positive roots 
$\Delta_+$ is fixed as in 
Proposition \ref{rootsh}. We denote $s_{id}=s$.
Now using formula (\ref{todah}) we outline computation of the operators $M_i^h$. 
This computation is simplified by
the following lemma.
\begin{lemma}{ \bf (\cite{Et}, Lemma 5.2)} 
Let $X=f_{\gamma_1}...f_{\gamma_n}$. If the roots
$\gamma_1,...,\gamma_n$ are not all simple 
then $L(X)=0$. 
Otherwise, if $\gamma_i=\alpha_{k_i}$, then 
$$
\jmath^*(L(X))f(e^{hy})=e^{-h(\sum\alpha_{k_i},y)}f(e^{hy})\prod_i\overline 
\chi_h^{s}(f_{k_i})
$$ 
\end{lemma}
{\em Proof }  follows immediately from Proposition \ref{rootsh} and the 
arguments used in the proof of 
Proposition \ref{diffh}. 

Using this lemma we obtain that if $\beta$ is not a simple root then the term in 
(\ref{todah}) containing 
root vector $f_\beta$ gives a trivial contribution to the operators
$L(\rho_{\chi_h^{s}}(C_{V_i}))$. Note also that by Proposition \ref{rootsh} 
${\chi_h^{s}}(e_\beta)=0$ if 
$\beta$ is not a simple root. Therefore from formula (\ref{todah}) we have
\begin{equation}
\begin{array}{l}
L(\rho_{\chi_h^{s}}(C_{V_i}))=\\
\\
L(id\otimes tr_V)( e^{t_0}\prod_{i}
exp_{q^{-2d_i}}[(q_i-q_i^{-1})f_i \otimes e_ie^{-hd_i{1+s \over 1-s}H_i}]\times 
\\
\\
e^{t_0}\prod_{i}
exp_{q^{-2d_i}}[(q_i-q_i^{-1}){\chi_h^{s}}(e_i) \otimes 
e^{hd_i{1+s \over 1-s}H_i}f_i](1\otimes e^{2h\rho^\vee})).
\end{array}
\end{equation}

In particular, let ${\frak g}=sl(n)$, $V_1=V$ the fundamental representation of 
$sl(n)$. 
Then direct calculation gives
$$
M_1f(e^{hy})=\left( \sum_{j=1}^n T_j^2-
(q-q^{-1})^2\sum_{i=1}^{n-1}{\chi_h^{s}}(e_i){\overline \chi_h^{s}}(f_i)
e^{-h(y,\alpha_i)}T_{i+1}T_i\right) f(e^{hy}),
$$
where $T_i=T_{\omega_i}$, ${\omega_i}$ are the weights of $V$. 
The last expression coincides with formula (5.7) obtained in \cite{Et}.

%%%%%%%%%%%%%%%%%%%%%%%%%%%%%%%%%%%%%%%%%%%%%%%%%%%%%%%%%%%%%%%%%%%%%%%%%%%%%%%%
%%%%%%%%%%%%%%%%%%%%%%%%%%%
%%%%%%%%%%%%%%%%%%%%%%%%%%%%%%%%%%%%%%%%%%%%%%%%%%%%%%%%%%%%%%%%%%%%%%%%%%%%%%%%
%%%%%%%%%%%%%%%%%%%%%%%%%%%%%
%%%%%%%%%%%%%%%%%%%%%%%%%%%%%%%%%%%%%%%%%%%%%%%%%%%%%%%%%%%%%%%%%%%%%%%%%%%%%%%%
%%%%%%%%%%%%%%%%%%%%%%%%%%%%

\chapter{Poisson--Lie groups and Whittaker model}\label{GWitt}

In this Chapter we introduce another quantum version of the Whittaker model.
We consider quantizations of algebras of regular functions on algebraic 
Poisson--Lie
groups. We define the Whittaker model of the center of these quantum algebras. 
The algebraic structure of this model is related to the structure of the set of 
regular
elements in the corresponding algebraic group. This relation is parallel to the 
one
established by Kostant for Lie algebras (see Section \ref{geomappr}). Our main 
geometric 
result is an analog of Theorem C for algebraic groups. A generalization of this 
theorem
for loop groups is contained in \cite{SS}.

%%%%%%%%%%%%%%%%%%%%%%%%%%%%%%%%%%%%%%%%%%%%%%%%%%%%%%%%%%%%%%%%%%%%%%%%%%%%%%%%
%%%%%%%%%%%%%%%%%%%%%%%%%

\section{Poisson--Lie groups}

\setcounter{equation}{0}
\setcounter{theorem}{0}

Recall some notions concerned with Poisson--Lie groups (see \cite{Dm}, 
\cite{fact}, \cite{dual}, \cite{ChP}). 
Let $G$ be a finite--dimensional Lie group equipped with a Poisson bracket, 
$\frak g$ its Lie algebra. $G$ is called 
a Poisson--Lie group if the multiplication $G\times G \rightarrow G$ is a 
Poisson map. 
A Poisson bracket satisfying this axiom is degenerate and, in particular, is 
identically zero 
at the unit element of the group. Linearizing this bracket at the unit element 
defines the 
structure of a Lie algebra in the space $T^*_eG\simeq {\frak g}^*$. 
The pair (${\frak g},{\frak g}^{*})$ is called the tangent bialgebra of $G$.

Lie brackets in $\frak{g}$ and $\frak{g}^{*}$ satisfy the following
compatibility condition:

{\em Let }$\delta: {\frak g}\rightarrow {\frak g}\wedge {\frak g}$ {\em be
the dual  of the commutator map } $[,]_{*}: {\frak g}^{*}\wedge
{\frak g}^{*}\rightarrow {\frak g}^{*}$. {\em Then } $\delta$ {\em is a
1-cocycle on} $  {\frak g}$ {\em (with respect to the adjoint action
of } $\frak g$ {\em on} ${\frak g}\wedge{\frak g}$).

Let $c_{ij}^{k}, f^{ab}_{c}$ be the structure constants of
${\frak g}, {\frak g}^{*}$ with respect to the dual bases $\{e_{i}\},
\{e^{i}\}$ in ${\frak g},{\frak g}^{*}$. The compatibility condition
means that
     
$$
c_{ab}^{s} f^{ik}_{s} ~-~ c_{as}^{i} f^{sk}_{b} ~+~ c_{as}^{k}
f^{si}_{b} ~-~ c_{bs}^{k} f^{si}_{a} ~+~ c_{bs}^{i} f^{sk}_{a} ~~=
~~0.
$$
This condition is symmetric with respect to exchange of $c$ and
$f$. Thus if $({\frak g},{\frak g}^{*})$ is a Lie bialgebra, then
$({\frak g}^{*}, {\frak g})$ is also a Lie bialgebra.

The following proposition shows that the category of finite--dimensional Lie 
bialgebras is isomorphic to
the category of finite--dimensional connected simply connected Poisson--Lie 
groups.
\begin{proposition}{\bf (\cite{ChP}, Theorem 1.3.2)}
If $G$ is a connected simply connected finite--dimensional Lie group, every 
bialgebra structure on $\frak g$
is the tangent bialgebra of a unique Poisson structure on $G$ which makes $G$ 
into a Poisson--Lie group.
\end{proposition}

Let $G$ be a finite--dimensional Poisson--Lie group, $({\frak g},{\frak g}^{*})$ 
the tangent bialgebra of $G$. 
The connected simply connected finite--dimensional 
Poisson--Lie group corresponding to the Lie bialgebra $({\frak g}^{*}, {\frak 
g})$ is called the dual
Poisson--Lie group and denoted by $G^*$.

$({\frak g},{\frak g}^{*})$ is called a {\em factorizable Lie
bialgebra }if the following conditions are satisfied (see \cite{fact} , 
\cite{Dm}):
\begin{enumerate}
\item
${\frak g}${\em \ is equipped with a non--degenerate invariant
scalar product} $\left( \cdot ,\cdot \right)$.

We shall always identify ${\frak g}^{*}$ and ${\frak g}$ by means of this
scalar product.

\item  {\em The dual Lie bracket on }${\frak g}^{*}\simeq {\frak g}${\em \
is given by} 
\begin{equation}
\left[ X,Y\right] _{*}=\frac 12\left( \left[ rX,Y\right] +\left[ X,rY\right]
\right) ,X,Y\in {\frak g},  \label{rbr}
\end{equation}
{\em where }$r\in {\rm End}\ {\frak g}${\em \ is a skew symmetric linear 
operator
(classical r-matrix).}

\item  $r${\em \ satisfies} {\em the} {\em modified classical Yang-Baxter
identity:} 
\begin{equation}
\left[ rX,rY\right] -r\left( \left[ rX,Y\right] +\left[ X,rY\right] \right)
=-\left[ X,Y\right] ,\;X,Y\in {\frak g}{\bf .}  \label{cybe}
\end{equation}
\end{enumerate}

Define operators $r_\pm \in {\rm End}\ {\frak g}$ by 
\[
r_{\pm }=\frac 12\left( r\pm id\right) . 
\]
We shall need some properties of the operators $r_{\pm }$.
Denote by ${\frak b}_\pm$ and ${\frak n}_\mp$ the image and the kernel of the 
operator 
$r_\pm $: 
\begin{equation}\label{bnpm}
{\frak b}_\pm = Im~r_\pm,~~{\frak n}_\mp = Ker~r_\pm.
\end{equation}
\begin{proposition}{\bf (\cite{BD}, \cite{rmatr})}\label{bpm}
Let $({\frak g}, {\frak g}^*)$ be a factorizable Lie bialgebra. Then

(i) ${\frak b}_\pm \subset {\frak g}$ is a Lie subalgebra, the subspace ${\frak 
n}_\pm$ is a Lie ideal in
${\frak b}_\pm,~{\frak b}_\pm^\perp ={\frak n}_\pm$.

(ii) ${\frak n}_\pm$ is an ideal in ${\frak {g}}^{*}$.

(iii) ${\frak b}_\pm$ is a Lie subalgebra in ${\frak {g}}^{*}$. Moreover ${\frak 
b}_\pm ={\frak {g}}^{*}/ {\frak n}_\pm$.

(iv) $({\frak b}_\pm,{\frak b}_\pm ^*)$ is a subbialgebra of $({\frak 
{g}},{\frak {g}}^{*})$ and
$({\frak b}_\pm,{\frak b}_\pm ^*)\simeq ({\frak b}_\pm,{\frak b}_\mp)$. The 
canonical paring
between ${\frak b}_\mp$ and ${\frak b}_\pm$is given by
\begin{equation}
(X_\mp ,Y_\pm )_\pm=(X_\mp,r_\pm^{-1}Y_\pm ) ,~ X_\mp \in {\frak b}_\mp ;~ Y_\pm 
\in {\frak b}_\pm .
\end{equation}
\end{proposition}
The classical Yang--Baxter equation implies that $r_{\pm }$ , regarded as a 
mapping from 
${\frak g}^{*}$ into ${\frak g}$, is a Lie algebra homomorphism.  
Moreover, $r_{+}^{*}=-r_{-},$\ and $r_{+}-r_{-}=id.$

Put ${\frak {d}}={\frak g\oplus {g}}$ (direct sum of two
copies). The mapping
\begin{eqnarray}\label{imbd}
{\frak {g}}^{*}\rightarrow {\frak {d}}~~~:X\mapsto (X_{+},~X_{-}),~~~X_{\pm
}~=~r_{\pm }X
\end{eqnarray}
is a Lie algebra embedding. Thus we may identify ${\frak g^{*}}$ with a Lie
subalgebra in ${\frak {d}}$.

Naturally, embedding (\ref{imbd}) extends to an embedding 
$$
G^*\rightarrow G\times G,~~L\mapsto (L_+,L_-).
$$
We shall identify $G^*$ with the corresponding subgroup in $G\times G$.

%%%%%%%%%%%%%%%%%%%%%%%%%%%%%%%%%%%%%%%%%%%%%%%%%%%%%%%%%%%%%%%%%%%%%%%%%%%%%%%%
%%%%%%%%%%%%%%%%%%%%%%%%

\section{Poisson reduction}\label{poisred}

\setcounter{equation}{0}
\setcounter{theorem}{0}

In this section we recall basic facts on Poisson reduction (see \cite{W}, 
\cite{RIMS}). 
These facts will be used in the proof of the group counterpart of Theorem E (see 
Section \ref{geomappr}).

Let $M,~B,~B'$ be Poisson manifolds. Two Poisson surjections 
$$
\begin{array}{ccccc}
&  & M &  &  \\ 
& \stackrel{\pi^{\prime } }{\swarrow } &  & \stackrel{\pi }{\searrow } &  \\ 
B^{\prime } &  &  &  & B
\end{array}
$$
form a dual pair if the pullback $\pi^{^{\prime
}*}C^\infty(B^{\prime })$ is the centralizer of $\pi^* C^\infty (B)$ in the 
Poisson algebra 
$C^\infty (M)$. In that case the sets $B^{\prime }_b=\pi^{\prime } \left( \pi
^{-1}(b) \right),~b\in
B$ are Poisson submanifolds in $B^{\prime }$ (see \cite{W}) called reduced 
Poisson manifolds.

Fix an element $b\in B$. Then the algebra of functions $C^\infty (B^{\prime 
}_b)$ may be described as
follows. Let $I_b$ be the ideal in $C^\infty (M)$ generated by elements
${\pi}^*(f),~f\in C^\infty (B),~f(b)=0$. Denote $M_b=\pi^{-1}(b)$. Then the 
algebra $C^\infty (M_b)$ 
is simply the quotient of $C^\infty (M)$ by $I_b$.  Denote by 
$P_b:C^\infty (M)\rightarrow C^\infty (M)/I_b=C^\infty (M_b)$ 
the canonical projection onto the quotient.
\begin{lemma}\label{redspace}
Suppose that the map $f\mapsto f(b)$ is 
a character of the Poisson algebra $C^\infty (B)$. Then one can define an action 
of 
the Poisson algebra $C^\infty (B)$ on the space $C^\infty (M_b)$ by
\begin{equation}\label{redact}
f\cdot \varphi=P_b(\{ {\pi}^*(f), \tilde \varphi \}),
\end{equation}
where $f\in C^\infty (B)$, $\varphi \in C^\infty (M_b)$, $\tilde \varphi \in 
C^\infty (M)$ is a 
representative of $\varphi$ in $C^\infty (M)$ such that $P_b(\tilde 
\varphi)=\varphi$.
Moreover, $C^\infty (B^{\prime }_b)$ is the subspace of invariants in $C^\infty 
(M_b)$
with respect to this action. 
\end{lemma}
{\em Proof.} 
Let $\varphi \in C^\infty (M_b)$. Choose a representative 
$\tilde \varphi \in C^\infty (M)$ such that $P_b(\tilde \varphi)=\varphi$.
Since the map $f\mapsto f(b)$ is 
a character of the Poisson algebra $C^\infty (B)$, Hamiltonian vector fields of 
functions 
${\pi}^*(f),~f\in C^\infty (B)$ are tangent to the surface $M_b$. 
Therefore using the the definition of the dual 
pair we obtain that
$\varphi={\pi^{\prime }}^*(\psi)$ for some $\psi \in C^\infty(B^{\prime }_b)$ if 
and only if
$P_b(\{ {\pi}^*(f), \tilde \varphi\})=0$ for every $f\in C^\infty (B)$.
Note also that the r.h.s. of (\ref{redact}) only depends on $\varphi$ but not on 
the representative
$\tilde \varphi$, and hence
formula (\ref{redact})
defines an action of the Poisson algebra $C^\infty (B)$ on the space $C^\infty 
(M_b)$.
Finally we obtain that $C^\infty (B^{\prime }_b)$ is exactly the subspace of 
invariants in $C^\infty (M_b)$
with respect to this action. 
\begin{definition}
The algebra $C^\infty (B^{\prime }_b)$ is called a reduced Poisson algebra.
We also denote it by $C^\infty (M_b)^{C^\infty (B)}$.
\end{definition}
\begin{remark}\label{redpoisalg}
Note that the description of the algebra $C^\infty (M_b)^{C^\infty (B)}$ 
obtained in Lemma \ref{redspace} 
is independent
of both the manifold $B^{\prime }$ and the projection $\pi^{\prime }$.
Observe also that the reduced space $B^{\prime }_b$ may be identified with a 
cross--section 
of the action of the Poisson algebra $C^\infty (B)$ on $M_b$ by Hamiltonian 
vector fields.
In particular, $B^{\prime }_b$ may be regarded as a submanifold in $M_b$.
\end{remark} 

An important example of dual pairs is provided by Poisson group actions.
Recall that a Poisson group action of a Poisson--Lie group $A$ on a Poisson 
manifold $M$ 
is a group action $A\times M\rightarrow M$ which is also a Poisson map (as
usual, we suppose that $A\times M$ is equipped with the product Poisson
structure). In \cite{RIMS} it is proved that if the space $M/A$ is a smooth 
manifold,
there exists a unique Poisson structure on $M/A$
such that the canonical projection $M\rightarrow M/A$ is a Poisson map.
 
Let $\frak a$ be the Lie algebra of $A$. Denote by $\langle\cdot,\cdot\rangle$ 
the 
canonical paring between ${\frak a}^*$ and $\frak a$.
A map $\mu :M\rightarrow A^*$ is called a moment map for a right Poisson group 
action 
$A\times M\rightarrow M$ if (\cite{Lu})
\begin{equation}
L_{\widehat X} \varphi =\langle \mu^*(\theta_{A^*}) , X \rangle (\xi_\varphi ) ,
\end{equation}
where $\theta_{A^*}$ is the universal right--invariant Maurer--Cartan form on 
$A^*$, $X \in {\frak a}$, 
$\widehat X$ is the corresponding vector field on $M$ and
$\xi_\varphi $ is the Hamiltonian vector field of $\varphi \in C^\infty (M)$.

By Theorem 4.9, \cite{Lu} one can always equip $A^*$ with a Poisson structure in 
such a way
that $\mu$ becomes a Poisson mapping. Then 
from the definition of the moment map it follows that if $M/A$ is a smooth
manifold then the canonical projection $M\rightarrow M/A$ and the moment map 
$\mu:M\rightarrow A^*$
form a dual pair (see \cite{Lu} for details).

The main example of Poisson group actions is the so--called dressing action.
The dressing action may be described as follows (see \cite{Lu}, \cite{RIMS}).
\begin{proposition}\label{dressingact}
Let $G$ be a connected simply connected Poisson--Lie group with factorizable 
tangent Lie bialgebra, 
$G^*$ the dual group. Then there exists a unique right Poisson group action
$$
G^*\times G\rightarrow G^*,~~((L_+,L_-),g)\mapsto g\circ (L_+,L_-),
$$
such that the identity mapping $\mu: G^* \rightarrow G^*$ is the moment map for 
this action.

Moreover, let $q:G^* \rightarrow G$ be the map defined by
$$
q(L_+,L_-)=L_-L_+^{-1}.
$$
Then
$$
q(g\circ (L_+,L_-))=g^{-1}L_-L_+^{-1}g.
$$
\end{proposition}
The notion of Poisson group actions may be generalized as follows.
Let $A\times M \rightarrow M$ be a Poisson group action of a Poisson--Lie
group $A$ on a Poisson manifold $M$.
A subgroup $K\subset A$ is called {\em admissible} if the set $%
C^\infty \left( M\right) ^K$ of $K$-invariants is a Poisson subalgebra in $%
C^\infty \left( M\right)$. If space $M/K$ is a smooth manifold, we may identify 
the algebras 
$C^\infty(M/K)$ and $C^\infty \left( M\right) ^K$.  Hence there exists a Poisson 
structure on $M/K$
such that the canonical projection $M\rightarrow M/K$ is a Poisson map. 
\begin{proposition}\label{admiss}{\bf (\cite{RIMS}, Theorem 6; \cite{Lu}, \S 2)}
Let $\left( {\frak a},{\frak a}^{*}\right) $ be the tangent
Lie bialgebra of $A.$ A connected Lie subgroup $K\subset A$ with Lie algebra 
${\frak k}\subset {\frak a}$ is admissible if ${\frak k}^{\perp }\subset
{\frak a}^{*}$ is a Lie subalgebra.
\end{proposition}
We shall need the following particular example of dual pairs arising from 
Poisson group actions.

Let $A\times M \rightarrow M$ be a right Poisson group action of a Poisson--Lie 
group $A$ on a manifold $M$.
Suppose that this action possesses a moment mapping $\mu : M\rightarrow A^*$.
Let $K$ be an admissible subgroup in $A$. Denote by $\frak k$ the Lie algebra of 
$K$. 
Assume that ${\frak k}^\perp \subset {\frak a}^*$ is a Lie subalgebra in ${\frak 
a}^*$. 
Suppose also that there is a splitting ${\frak a}^*={\frak t}\oplus {\frak 
k}^\perp$, and that 
$\frak t$ is a Lie subalgebra in ${\frak a}^*$. Then the linear space ${\frak 
k}^*$ is naturally
identified with $\frak t$.
Assume that $A^*$ is the semidirect 
product of the Lie subgroups $K^\perp , T$ corresponding to the Lie algebras
${\frak k}^\perp , {\frak t}$ respectively. Suppose that $K^\perp$ is a 
connected subgroup in $A^*$.
Fix the decomposition
$A^*=K^\perp T$ and denote by $\pi_{K^\perp} , \pi_{T}$ the projections onto
$K^\perp$ and $T$ in this decomposition. 
\begin{proposition}\label{QPmoment}
Define a map $\overline{\mu}:M\rightarrow T$ by
$$
\overline{\mu}=\pi_{T}\mu.
$$ 
Then
 
(i) 
$\overline{\mu}^*\left( C^\infty \left( T\right)\right)$ is a Poisson subalgebra 
in $C^\infty \left( M\right)$,
and hence one can equip $T$ with a Poisson structure such that 
$\overline{\mu}:M\rightarrow T$ is 
a Poisson map.

(ii)Moreover, the algebra $C^\infty \left( M\right) ^K$ is the centralizer of 
$\overline{\mu}^*\left( C^\infty \left( T\right)\right)$ in the Poisson algebra 
$C^\infty \left( M\right)$.
In particular, if $M/K$ is a smooth manifold the maps
\begin{equation}\label{dp}
\begin{array}{ccccc}
&  & M &  &  \\ 
& \stackrel{\pi }{\swarrow } &  & \stackrel{\overline{\mu}}{\searrow } &  ,\\ 
M/K &  &  &  & T
\end{array}
\end{equation}
form a dual pair.
\end{proposition}
{\em Proof.} (i)First, by Theorem 4.9 in \cite{Lu} there exists a Poisson 
bracket on $A^*$ such that 
$\mu :M\rightarrow A^*$ is a Poisson map. Moreover, we can choose this bracket 
to be the sum of the standard
Poisson--Lie bracket of $A^*$ and of a left invariant bivector on $A^*$. 
Denote by $A^*_M$ the manifold $A^*$ equipped with this Poisson structure.
Now observe that $T$ is identified with the quotient $K^\perp \setminus A^*_M$, 
where $K^\perp$ acts on
$A^*_M$ by multiplications from the left. Therefore to prove part (i) of the 
proposition it suffices to show
that $K^\perp$--invariant functions on $A^*_M$ form a Poisson subalgebra in 
$C^\infty(A^*_M)$.

Observe that since $A^*$ is a Poisson--Lie group and the Poisson structure of 
$A^*_M$ is obtained from that of
$A^*$ by adding a left--invariant term, the action of $A^*$ on $A^*_M$ by 
multiplications from the left is
a Poisson group action. Note also that $K^\perp$ is a connected subgroup in 
$A^*$ and $({\frak k}^\perp)^\perp 
\cong {\frak k}$ is a Lie subalgebra in $\frak a$. Therefore by Proposition 
\ref{admiss} $K^\perp$ is 
an admissible subgroup in $A^*$. Therefore 
$K^\perp$--invariant functions on $A^*_M$ form a Poisson subalgebra in 
$C^\infty(A^*_M)$, and hence
$\overline{\mu}^*\left( C^\infty \left( T\right)\right)$ is a Poisson subalgebra 
in $C^\infty \left( M\right)$.
This proves part (i).

(ii)By the definition of the moment map we have:
\begin{equation}\label{X5}
L_{\widehat X} \varphi =\langle \mu^*(\theta_{A^*}) , X \rangle (\xi_\varphi ) ,
\end{equation}
where $X \in {\frak a} , \widehat X$ is the corresponding vector field on $M$ 
and
$\xi_\varphi $ is the Hamiltonian vector field of $\varphi \in C^\infty (M)$. 
Since $A^*$ is the semidirect product of $K^\perp$ and $T$ the pullback of the 
right--invariant Maurer--Cartan form $\mu^*(\theta_{A^*})$ may be represented as 
follows:
$$
\mu^*(\theta_{A^*})= {\rm Ad}(\pi_{K^\perp}\mu 
)({\overline{\mu}}^*\theta_{T})+(\pi_{K^\perp}\mu )^*\theta_{K^\perp},
$$
where ${\rm Ad}(\pi_{K^\perp}\mu )({\overline{\mu}}^*\theta_{T})\in {\frak 
t},~(\pi_{K^\perp}\mu )^*\theta_{K^\perp}\in {\frak k}^\perp$.

Now let $X \in {\frak k}$. Then $\langle (\pi_{K^\perp}\mu 
)^*\theta_{K^\perp}),X\rangle =0$ and formula (\ref{X5}) takes the form:
\begin{equation}\label{+}
\begin{array}{l}
L_{\widehat X} \varphi =
\langle {\rm Ad}(\pi_{K^\perp}\mu )({\overline{\mu}}^*\theta_{T}),X \rangle 
(\xi_\varphi )=\\
\\
\langle {\rm Ad}(\pi_{K^\perp}\mu )(\theta_{T}),X \rangle 
({\overline{\mu}}_*(\xi_\varphi )) .
\end{array}
\end{equation}

Since ${\rm Ad}(\pi_{K^\perp}\mu )$ is a non--degenerate transformation, 
$L_{\widehat X} \varphi =0$ for every $X\in {\frak k}$
if and only if ${\overline{\mu}}_*(\xi_\varphi )=0$, i.e. a function $\varphi 
\in C^\infty (M)$ is 
$K$--invariant if and only if $\{ \varphi ,\overline{\mu}^*(\psi) \}=0$ for 
every 
$\psi \in C^\infty (T)$. This completes the proof.
\begin{remark}\label{remred}
Let $t\in T$ be as in Lemma \ref{redspace}.
Assume that $\pi(\overline{\mu}^{-1}(t))$ is a smooth manifold ($M/K$ does not 
need to be smooth). Then  
the algebra $C^\infty(\pi(\overline{\mu}^{-1}(t)))$ is isomorphic to
the reduced Poisson algebra $C^\infty(\overline{\mu}^{-1}(t))^{C^\infty(T)}$.
\end{remark}

%%%%%%%%%%%%%%%%%%%%%%%%%%%%%%%%%%%%%%%%%%%%%%%%%%%%%%%%%%%%%%%%%%%%%%%%%%%%%%%%
%%%%%%%%%%%%%%%%%%%%%%%%

\section{Quantization of Poisson--Lie groups and Whittaker model}

\setcounter{equation}{0}
\setcounter{theorem}{0}

Let $\frak g$ be a finite--dimensional complex simple Lie algebra. Observe that 
cocycle (\ref{cocycles}) equips 
$\frak g$ with
the structure of a factorizable Lie bialgebra. For simplicity we suppose that 
$\pi =id$, and denote $s_{id}=s$.
Using the identification 
${\rm End}~{\frak g}\cong {\frak g}\otimes {\frak g}$ the corresponding 
r--matrix may be represented as
$$
r^{s}=P_+-P_-+{1+s \over 1-s}P_0,
$$
where $P_+,P_-$ and $P_0$ are the projection operators onto ${\frak n}_+,{\frak 
n}_-$ and $\frak h$ in 
the direct sum
$$
{\frak g}={\frak n}_+ +{\frak h} + {\frak n}_-.
$$

Let $G$ be the connected simply connected simple Poisson--Lie group with the 
tangent Lie bialgebra $({\frak g},{\frak g}^*)$, 
$G^*$ the dual group. Observe that $G$ is an algebraic group (see \S 104, 
Theorem 12 in \cite{Z}). 

Note also that 
$$
r^{s}_+=P_+ + {1 \over 1-s}P_0,~~r^{s}_-=-P_- + {s \over 1-s}P_0,
$$
and hence the subspaces ${\frak b}_\pm$ and ${\frak n}_\pm$ defined by 
(\ref{bnpm}) coincide with
the Borel subalgebras in $\frak g$ and their nil--radicals, respectively.
Therefore every element $(L_+,L_-)\in G^*$ may be uniquely written as
\begin{equation}\label{fact}
(L_+,L_-)=(h_+,h_-)(n_+,n_-),
\end{equation}
where $n_\pm \in N_\pm$, $h_+=exp({1 \over 1-s}x),~h_-=exp({s \over 1-s}x),~x\in 
{\frak h}$.
In particular, $G^*$ is a solvable algebraic subgroup in $G\times G$.

For every algebraic variety $V$ we denote by ${\cal F}(V)$ the algebra of 
regular functions on $V$.
Our main object will be the algebra of regular functions on $G^*$, ${\cal 
F}(G^*)$.
This algebra may be explicitly described as follows.
Let $\pi_V$ be a finite--dimensional representation of $G$. Then matrix
elements of $\pi_V(L_\pm)$ are well--defined functions on $G^*$, and ${\cal 
F}(G^*)$ is the subspace
in $C^\infty(G^*)$ generated by matrix elements of $\pi_V(L_\pm)$, where $V$ 
runs through all finite--dimensional
representations of $G$.

The elements $L^{\pm,V}=\pi_V(L_\pm)$ may be viewed as elements of the space 
${\cal F}(G^*)\otimes {\rm End}V$. For every two finite--dimensional ${\frak g}$ 
modules $V$ and $W$
we denote ${r^s_+}^{VW}=(\pi_V\otimes \pi_W)r^s_+$, where $r^s_+$ is regarded as 
an 
element of ${\frak g}\otimes {\frak g}$. 
\begin{proposition}{\bf (\cite{dual}, Section 2)}
${\cal F}(G^*)$ is a Poisson subalgebra in the Poisson algebra $C^\infty(G^*)$, 
the Poisson brackets
of the elements $L^{\pm,V}$ are given by
\begin{equation}\label{pbf}
\begin{array}{l}
\{L^{\pm,W}_{1},L^{\pm,V}_{2}\}~=~
2[{r_+^s}^{VW},L^{\pm,W}_{1}L^{\pm,V}_{2}],\\
\\
\{L^{-,W}_{1},L^{+,V}_{2}\}~=2[{r_{+}^s}^{VW},L^{-,W}_{1}L^{+,V}_{2}],
\end{array}
\end{equation}
where 
$$
L^{\pm,W}_1=L^{\pm,W}\otimes I_V,~~L^{\pm,V}_2=I_W\otimes L^{\pm,V},
$$
and $I_X$ is the unit matrix in $X$.

Moreover, the map $\Delta:{\cal F}(G^*)\rightarrow {\cal F}(G^*)\otimes {\cal 
F}(G^*)$ dual to the multiplication in $G^*$,
\begin{equation}\label{comultcl}
\Delta(L^{\pm,V}_{ij})=\sum_k L^{\pm,V}_{ik}\otimes L^{\pm,V}_{kj},
\end{equation}
is a homomorphism of Poisson algebras. 
\end{proposition}
\begin{remark}
Recall that a Poisson--Hopf algebra is a Poisson algebra which is also a Hopf 
algebra such that the 
comultiplication is a homomorphism of Poisson algebras. According to Proposition 
\ref{pbf}
${\cal F}(G^*)$ is a Poisson--Hopf algebra.
\end{remark}

Now we describe a quantization of the Poisson--Hopf algebra ${\cal F}(G^*)$. 
Let $\tilde U_h^{s}({\frak g})$ be the subalgebra in $U_h^{s}({\frak g})$ 
topologically
generated, in the sense of formal power series over ${\Bbb C}[[h]]$, by elements
$\tilde H_i=hH_i,~i=1,\ldots l,~\tilde e_{\beta}=he_{\beta},~\tilde 
f_{\beta}=hf_{\beta},~\beta\in \Delta_+$. 

In fact $\tilde U_h^{s_\pi}({\frak g})$
is a Hopf subalgebra in $U_h^{s}({\frak g})$, explicit formulas for the 
comultiplication may be
obtained using Proposition 8.3 in \cite{kh-t}.
\begin{proposition} 
$\tilde U_h^{s}({\frak g})$ is a quantum formal series Hopf algebra (or QFSH 
algebra), i.e. 
$\tilde U_h^{s}({\frak g})$ is isomorphic as a ${\Bbb C}[[h]]$--module to 
$Map(I,{\Bbb C}[[h]])$ for
some set $I$, and 
$\tilde U^{s}({\frak g})=\tilde U_h^{s}({\frak g})/h\tilde U_h^{s}({\frak 
g})\cong {\Bbb C}[[\xi_1,\xi_2,\ldots ]]$ as a
topological algebra, for some (possibly infinite) sequence of indeterminates 
$\xi_1,\xi_2,\ldots $.
\end{proposition}
{\em Proof} is similar to the proof of the same result for $U_h({\frak g})$ (see 
Section 8.3 C in
\cite{ChP}).

Note that $\tilde U^{s}({\frak g})$ is naturally a Poisson--Hopf
algebra, the Poisson bracket is given by
\begin{equation}\label{quasipb}
\{x_1,x_2\}={[a_1,a_2] \over h}~(\mbox{mod }h),
\end{equation}
if $a_1,a_2\in \tilde U_h^{s}({\frak g})$ reduce to 
$x_1,x_2\in \tilde U^{s}({\frak g})~(\mbox{mod }h)$.
 
For any finite--dimensional $U_h^{s}({\frak g})$ module $V[[h]]$ we denote by 
${^h{L^{\pm,V}}}$ the following
elements of $U_h^{s}({\frak g})\otimes {\rm End}V[[h]]$ (see \cite{FRT}):
$$
{^h{L^{+,V}}}=(id\otimes \pi_V){{\cal R}_{21}^{s}}^{-1}=(id\otimes 
\pi_VS^{s}){\cal R}_{21}^{s}
,~~ {^h{L^{-,V}}}=(id\otimes \pi_V){\cal R}^{s}.
$$

We also denote $R^{VW}=(\pi_V\otimes \pi_W){\cal R}^{s}$. 
Observe that from formula (\ref{rmatrspi}) it follows that actually 
${^h{L^{\pm,V}}}\in \tilde U_h^{s}({\frak g})\otimes {\rm End}V[[h]]$.
If we fix a basis in $V[[h]]$, ${^h{L^{\pm,V}}}$ may be regarded as matrices 
with matrix 
elements $({^h{L^{\pm,V}}})_{ij}$ being elements of $\tilde U_h^{s}({\frak g})$. 
From the Yang--Baxter equation for $\cal R$ we get 
relations between $L^{\pm,V}$:
\begin{equation}\label{ppcomm}
\begin{array}{l}
R^{VW}{^h{L^{\pm,W}_1}}{^h{L^{\pm,V}_2}}={^h{L^{\pm,V}_2}}{^h{L^{\pm,W}_1}}R^{VW
},
\end{array}
\end{equation}
\begin{equation}\label{pmcomm}
R^{VW}{^h{L^{-,W}_1}}{^h{L^{+,V}_2}}={^h{L^{+,V}_2}}{^h{L^{-,W}_1}}R^{VW}.
\end{equation}
By ${^h{L^{\pm,W}_1}},~{^h{L^{\pm,V}_2}}$ we understand the following matrices 
in $V\otimes W$:
$$
{^h{L^{\pm,W}_1}}={^h{L^{\pm,W}}}\otimes I_V,~~{^h{L^{\pm,V}_2}}=I_W\otimes 
{^h{L^{\pm,V}}},
$$
where $I_X$ is the unit matrix in $X$.

From (\ref{rmprop}) we can obtain the action of the comultiplication on the 
matrices ${^h{L^{\pm,V}}}$:
\begin{equation}\label{comult}
\Delta_s({^h{L^{\pm,V}_{ij}}})=\sum_k {^h{L^{\pm,V}_{ik}}}\otimes 
{^h{L^{\pm,V}_{kj}}}.
\end{equation}

We denote by ${\cal F}_h(G^*)$ the Hopf subalgebra in $\tilde U_h^{s}({\frak 
g})$ generated in the sense of
$h$--adic topology by matrix elements
of ${^h{L^{\pm,V}}}$, where $V$ runs through all finite--dimensional 
representations of $\frak g$.
\begin{proposition}\label{quantreg}
Denote by $p:\tilde U_h^{s}({\frak g})\rightarrow \tilde U^{s}({\frak g})$ the 
canonical projection.
Then $p({\cal F}_h(G^*))$ is isomorphic to ${\cal F}(G^*)$ as a Poisson--Hopf 
algebra.
\end{proposition}
{\em Proof.}
Denote ${\cal F}(G^*)'=p({\cal F}_h(G^*)),~{\tilde 
L^{\pm,V}}=p({^h{L^{\pm,V}}})\in 
{\cal F}(G^*)'\otimes {\rm End}V$. 
First observe that the map 
$$
\imath :{\cal F}(G^*)'\rightarrow {\cal F}(G^*),~~(\imath \otimes id){\tilde 
L^{\pm,V}}={L^{\pm,V}}
$$
is a well--defined linear isomorphism.
Indeed, consider, for instance, element ${\tilde L^{-,V}}$.
From (\ref{rmatrspi}) it follows that
\begin{equation}
\begin{array}{l}
{\tilde L^{-,V}}_{ij}=\{ exp\left[ \sum_{i=1}^l-2p(hY_i)\otimes \pi_V({s \over 
1-s}H_i)\right]\times \\
\prod_{\beta}
exp[2(X_\beta,X_{-\beta})^{-1}p(he_{\beta}) \otimes 
\pi_V(X_{-\beta})]\}_{ij}.
\end{array}
\end{equation}

On the other hand (\ref{fact}) implies that every element $L_-$ may be 
represented in the form
\begin{equation}
\begin{array}{l}
L_- = exp\left[ \sum_{i=1}^lb_i{s \over 1-s}H_i\right]\times \\
\prod_{\beta}
exp[b_{\beta}X_{-\beta}],~b_i,b_\beta\in {\Bbb C},
\end{array}
\end{equation}
and hence
\begin{equation}
\begin{array}{l}
L^{-,V}_{ij}=\{ exp\left[ \sum_{i=1}^lb_i\otimes \pi_V({s \over 
1-s}H_i)\right]\times \\
\prod_{\beta}
exp[b_{\beta} \otimes 
\pi_V(X_{-\beta})]\}_{ij}.
\end{array}
\end{equation}
Therefore $\imath$ is a linear isomorphism. We have to prove that $\imath$ is an 
isomorphism of 
Poisson--Hopf algebras.

Recall that ${\cal R}^{s}=1\otimes 1 +2hr_+^{s}$ (mod $h^2$). Therefore from
commutation relations (\ref{ppcomm}), (\ref{pmcomm}) it follows that ${\cal 
F}(G^*)'$ is a commutative
algebra, and the Poisson brackets of matrix elements ${\tilde L^{\pm,V}}_{ij}$ 
(see (\ref{quasipb}))
are given by (\ref{pbf}), where $L^{\pm,V}$ are replaced by ${\tilde 
L^{\pm,V}}$. From (\ref{comult})
we also obtain that the action of the comultiplication on the matrices ${\tilde 
L^{\pm,V}}$ is given by 
(\ref{comultcl}), where $L^{\pm,V}$ are replaced by ${\tilde L^{\pm,V}}$.
This completes the proof.

We shall call the map $p:{\cal F}_h(G^*) \rightarrow {\cal F}(G^*)$ the 
quasiclassical limit.

Now using the Hopf algebra ${\cal F}_h(G^*)$ we shall define another quantum 
version of 
the Whittaker model $W({\frak b}_-)$.
Let ${\cal F}_h(N_\pm)$ be the subalgebras in ${\cal F}_h(G^*)$ generated by 
matrix elements of the matrices
$N^{-,V}=(id\otimes \pi_V){\cal R}^{s}_\Delta,~N^{+,V}=(id\otimes \pi_V){{\cal 
R}_{\Delta}^{s}}^{-1}_{21}$, where
$$
{\cal R}^{s}_\Delta=
\prod_{\beta}exp_{q_{\beta}^{-1}}[(q-q^{-1})a(\beta)^{-1} e_{\beta} \otimes 
e^{h{1+s \over 1-s} \beta^\vee}f_{\beta}].
$$

Suppose that the ordering of the root system $\Delta_+$ is fixed as in 
Proposition \ref{rootsh}.
Then by Proposition \ref{rootsh} the map $\chi_h^{s}:{\cal F}_h(N_-)\rightarrow 
{\Bbb C}$ 
defined by
\begin{equation}\label{charq}
(\chi_h^{s}\otimes id)(N^{-,V})=
\prod_{i=1}^lexp_{q_{\alpha_i}^{-1}}[{(q_i-q_i^{-1})\over h}c_{i} \otimes 
\pi_V(e^{hd_i{1+s \over 1-s}H_i}f_i)],
c_i\in {\Bbb C}[[h]],~c_i\neq 0
\end{equation}
is a character of ${\cal F}_h(N_-)$.

We also denote by ${\cal F}_h(H)$ the intersection $U_h^{s}({\frak h})\cap {\cal 
F}_h(G^*)$. 
Clearly, ${\cal F}_h(H)$ is a commutative subalgebra in ${\cal F}_h(G^*)$.
From commutation relations (\ref{pmcomm}) one can obtain the following weak 
version of 
the Poincar\'{e}--Birkhoff--Witt theorem for ${\cal F}_h(G^*)$.
\begin{proposition}\label{Pbw}
Multiplication defines an isomorphism of ${\Bbb C}[[h]]$--modules
$$
{\cal F}_h(N_+)\otimes {\cal F}_h(H)\otimes {\cal F}_h(N_-)\rightarrow {\cal 
F}_h(G^*).
$$
\end{proposition}
Define ${\cal F}_h(B_\pm)={\cal F}_h(N_\pm){\cal F}_h(H)$.
Let ${{\cal F}_h(N_-)}_{\chi_h^{s}}$ be the kernel of the character $\chi_h^{s}$ 
so that one has a direct sum
\begin{equation}\label{ker}
{\cal F}_h(N_-)={\Bbb C}[[h]]\oplus {{\cal F}_h(N_-)}_{\chi_h^{s}}.
\end{equation}

From Proposition \ref{Pbw} and formula (\ref{ker}) we obtain also the direct sum
\begin{equation}\label{maindecqg}
{\cal F}_h(G^*)={\cal F}_h(B_+)\oplus I_{\chi_h^{s}},
\end{equation}
where $I_{\chi_h^{s}}={\cal F}_h(G^*){{\cal F}_h(N_-)}_{\chi_h^{s}}$ is the 
left--sided
ideal generated by ${{\cal F}_h(N_-)}_{\chi_h^{s}}$.

Denote by $\rho_{\chi_h^{s}}$ the projection onto ${\cal F}_h(B_+)$ in the 
direct sum 
(\ref{maindecqg}).
Let $Z({\cal F}_h(G^*))$ be the center of ${\cal F}_h(G^*)$.
Similarly to the classical case we define a subspace $W_h(B_+)$ in ${\cal 
F}_h(B_+)$ by
$W_h(B_+)=\rho_{\chi_h^{s}}(Z({\cal F}_h(G^*)))$.

To formulate the quantum version of Theorem A for $W_h(B_+)$ we recall that
for any finite--dimensional $\frak g$--module $V$ the element
$$
C_V=(id\otimes tr_V)((S^{s}\otimes id)(L^{+,V})L^{-,V}(1\otimes 
e^{2h\rho^\vee})),
$$
where $tr_V$ is the 
trace in $V[[h]]$, is central in ${\cal F}_h(G^*)$ (see formulas 
(\ref{centrelv}) and (\ref{S})).
\vskip 0.3cm
\noindent
{\bf Theorem $\bf A_q$}
{\em (i)The map
\begin{equation}\label{qgrho}
\rho_{\chi_h^{s_\pi}}:Z({\cal F}_h(G^*))\rightarrow W_h(B_+)
\end{equation}
is an isomorphism of algebras. 

(ii) The algebra $W_h(B_+)$ is freely generated as a commutative topological 
algebra over ${\Bbb C}[[h]]$ by  
the elements 
$C_{V_i}^{\rho_{\chi_h^{s}}}=\rho_{\chi_h^{s}}(C_{V_i}),~i=1,\ldots ,l$, where 
$V_i,~i=1,\ldots l$ are the fundamental representations of $\frak g$.}
\vskip 0.3cm
\noindent
{\em Proof} of (i) is similar to that of Theorem A in the classical case.
Part (ii) will be proved in Section \ref{cross}.
\begin{corollary}
The algebra $Z({\cal F}_h(G^*))$ is freely generated as a commutative 
topological algebra over ${\Bbb C}[[h]]$ by the 
elements $C_{V_i}$, 
where $V_i,~i=1,\ldots l$ are the fundamental representations of $\frak g$.
\end{corollary}
\vskip 0.3cm
\noindent
{\bf Definition $\bf A_q$}
{\em The algebra 
$$
W_h(B_+)=\rho_{\chi_h^{s}}(Z({\cal F}_h(G^*))).
$$
is called the Whittaker model of the center $Z({\cal F}_h(G^*))$.}
\vskip 0.3cm
Now following Section \ref{whitt} (see Lemma A) we equip ${\cal F}_h(B_+)$ with 
a structure of a left ${\cal F}_h(N_-)$ module in such a
way that $W_h(B_+)$ is realized as the space of invariants with respect to this 
action.
For every $v\in {\cal F}_h(B_+)$ and $x\in {\cal F}_h(N_-)$ we put
\begin{equation}\label{mainactqg}
x\cdot v =\rho_{\chi_h^{s}}([x,v]).
\end{equation}

Consider the space ${\cal F}_h(B_+)^{{\cal F}_h(N_-)}$ of ${\cal F}_h(N_-)$ 
invariants 
in ${\cal F}_h(B_+)$ with respect to this
action. Clearly, 
$W_h(B_+)\subseteq {\cal F}_h(B_+)^{{\cal F}_h(N_-)}$.

To formulate the quantum version of Theorem B for $W_h(B_+)$ we have to impose a 
restriction on the coefficients
$c_i$ in (\ref{charq}).
Define an element $u\in N_-$ by 
\begin{equation}\label{u}
u=\prod_{i=1}^lexp[2d_ic_{i}^0 X_{-\alpha_i}],~c_i^0=c_i~(\mbox{mod }h),
\end{equation}
where the terms in the product are ordered as in (\ref{charq}).
The motivation for this definition will be explained in the next section.
\vskip 0.3cm
\noindent
{\bf Theorem $\bf B_q$ }
{\em 
Suppose that $u\in N_+sN_+\cap N_-$, where $s$ stands for a representative of 
the Coxeter element in $G$.
Then the space of ${\cal F}_h(N_-)$ invariants 
in ${\cal F}_h(B_+)$ with respect to the
action (\ref{mainactqg}) is isomorphic to $W_h(B_+)$, i.e.}
\begin{equation}\label{invqg}
{\cal F}_h(B_+)^{{\cal F}_h(N_-)}\cong W_h(B_+).
\end{equation}
\vskip 0.3cm
The proof of this theorem occupies two next sections.
\begin{remark}
The following lemma shows that the set $N_+sN_+\cap N_-$ is not empty.
\begin{lemma}{\bf (\cite{st}, Lemma 4.5)}\label{f}
Let $w_0\in W$ be the longest element; let $\tau \in Aut$ $\Delta _{+}$ be
the automorphism defined by $\tau \left( \alpha \right) =-w_0\alpha
,\alpha \in \Delta _{+}.$ Let $N_i\subset N_+$ be the 1-parameter subgroup
generated by the root vector $X_{\tau \left( \alpha _i\right) },i=1,\ldots l$. 
Choose an element $u_i\in N_i,u_i\neq 1.$ Then we have 
$w_0u_iw_0^{-1}\in B_+s_iB_+.$
We may fix $u_i$ in such a way that $w_0u_iw_0^{-1}\in N_+s_iN_+.$ Set 
$x=u_1u_{2}...u_l$. Then $f=w_0xw_0^{-1}\in N_+sN_+\cap N_-.$
\end{lemma}
\end{remark}

Similarly to Proposition \ref{hkhom} we also have the following homological 
description of $W_h(B_+)$.
\begin{proposition}
Suppose that the conditions of Theorem $B_q$ are satisfied. Then $W_h(B_+)$ is 
isomorphic to
$Hk^0({\cal F}_h(G^*),{\cal F}_h(N_-),\chi_h^{s})^{opp}$ as an associative 
algebra.
\end{proposition}

%%%%%%%%%%%%%%%%%%%%%%%%%%%%%%%%%%%%%%%%%%%%%%%%%%%%%%%%%%%%%%%%%%%%%%%%%%%%%%%%
%%%%%%%%%%%%%%%%%%%%%%%%%%

\section{Poisson reduction and the Whittaker model}

\setcounter{equation}{0}
\setcounter{theorem}{0}

In this section we start the proof of Theorem ${\rm B}_q$. We shall analyse the 
quasiclassical limit of 
the algebra ${\cal F}_h(B_+)^{{\cal F}_h(N_-)}$. Using results of Section 
\ref{poisred}
we realize this limit algebra as the algebra of functions on a reduced 
Poisson manifold.

Denote ${\cal F}(N_\pm)=p({\cal F}_h(N_\pm)),~{\cal F}(B_\pm)=p({\cal 
F}_h(B_\pm)),
~{\cal F}(H)=p({\cal F}_h(H))$. 
We denote by $\chi_h^{s}$ the character of the Poisson subalgebra ${\cal 
F}(N_-)$ such that
$\chi^{s}(p(x))=\chi_h^{s}(x)~(\mbox{mod }h)$ for every $x\in {\cal F}_h(N_-)$. 
From (\ref{charq}) we have
\begin{equation}\label{charcl}
(\chi^{s}\otimes id)p(N^{-,V})=
\prod_{i=1}^lexp[2d_ic_{i}^0 \otimes 
\pi_V(X_{-\alpha_i})],~~c_i^0=c_i~(\mbox{mod }h).
\end{equation}

Let ${{\cal F}(N_-)}_{\chi^{s}}$ be the kernel of the character $\chi^{s}$ 
so that one has a direct sum
\begin{equation}
{\cal F}(N_-)={\Bbb C}\oplus {{\cal F}(N_-)}_{\chi^{s}}.
\end{equation}

Similarly to (\ref{maindecqg}) we have the direct sum
\begin{equation}\label{maindeccl}
{\cal F}(G^*)={\cal F}(B_+)\oplus I_{\chi^{s}},
\end{equation}
where $I_{\chi^{s}}={\cal F}(G^*){{\cal F}(N_-)}_{\chi^{s}}$ is the left--sided
ideal generated by ${{\cal F}(N_-)}_{\chi^{s}}$.

Denote by $\rho_{\chi^{s}}$ the projection onto ${\cal F}(B_+)$ in the direct 
sum 
(\ref{maindeccl}).
Using Lemma ${\rm A}_q$ we define the quasiclassical limit of action 
(\ref{mainactqg}) by
\begin{equation}\label{mainactcl}
x\cdot v =\rho_{\chi^{s}}(\{x,v\} ),
\end{equation}
where $v\in {\cal F}(B_+)$ and $x\in {\cal F}(N_-)$.
We shall describe the space of invariants ${\cal F}(B_+)^{{\cal F}(N_-)}$ with 
respect to this action by analysing
``dual geometric objects''.

First observe that algebra ${\cal F}(B_+)^{{\cal F}(N_-)}$ is a particular 
example of the 
reduced Poisson algebra introduced in Lemma \ref{redspace}. 
Indeed, define a map $\mu_{N_-}:G^* \rightarrow N_-$ by
\begin{equation}\label{mun}
\mu_{N_+}(L_+,L_-)=n_-,
\end{equation}
where $n_-$ is given by (\ref{fact}). $\mu_{N_-}$ is a morphism of algebraic 
varieties.
We also note that by definition ${\cal F}(N_-)=\{ \varphi\in {\cal 
F}(G^*):\varphi=
\varphi(n_-)\}$. Therefore ${\cal F}(N_-)$ is generated by the pullbacks of 
regular functions on $N_-$.
Since ${\cal F}(N_-)$ is a Poisson subalgebra in ${\cal F}(G^*)$, and  regular 
functions
on $N_-$ are dense in $C^\infty(N_-)$ on every compact subset, we can equip the 
manifold $N_-$ with
the Poisson structure in such a way that $\mu_{N_+}$ becomes a Poisson mapping. 
Let $u$ be the element defined by (\ref{u}),
\begin{equation}
u=\prod_{i=1}^lexp[2d_ic_{i}^0 X_{-\alpha_i}]~ \in N_-.
\end{equation}
Then from (\ref{charcl}) it follows that $\chi^s(\varphi)=\varphi (u)$ 
for every $\varphi \in {\cal F}(N_-)$. $\chi^s$ naturally extends to a character 
of the
Poisson algebra $C^\infty(N_-)$.

Now applying Lemma \ref{redspace} for $M=G^*,~B=N_-,~\pi=\mu_{N_+},~b=u$ we can 
define the 
reduced Poisson algebra $C^\infty(\mu_{N_+}^{-1}(u))^{C^\infty(N_-)}$ (see also 
Remark \ref{redpoisalg}). 
Denote by $I_u$ the ideal in $C^\infty(G^*)$ generated by elements 
$\mu_{N_+}^*\psi,~\psi \in C^\infty(N_-),
~\psi(u)=0$. Let $P_u:C^\infty(G^*)\rightarrow 
C^\infty(G^*)/I_u=C^\infty(\mu_{N_+}^{-1}(u))$ be the 
canonical projection. Then the action (\ref{redact}) of $C^\infty(N_-)$ on 
$C^\infty(\mu_{N_+}^{-1}(u))$
takes the form:
\begin{equation}\label{actred}
\psi\cdot \varphi=P_u(\{ \mu_{N_+}^*\psi, \tilde \varphi\}),
\end{equation}
where $\psi \in C^\infty(N_-),~\varphi \in C^\infty(\mu_{N_+}^{-1}(u))$ and 
$\tilde \varphi \in C^\infty(G^*)$ 
is a representative of $\varphi$ such that $P_u\tilde \varphi=\varphi$.
\begin{lemma}\label{redreg}
$\mu_{N_+}^{-1}(u)$ is a subvariety in $G^*$. Moreover, the algebra 
${\cal F}(B_+)^{{\cal F}(N_-)}$ is isomorphic to the algebra of regular 
functions on  
$\mu_{N_+}^{-1}(u)$ which are invariant with respect to the action 
(\ref{actred}) of
$C^\infty(N_-)$ on $C^\infty(\mu_{N_+}^{-1}(u))$, i.e.
$$
{\cal F}(B_+)^{{\cal F}(N_-)}={\cal F}(\mu_{N_+}^{-1}(u))\cap 
C^\infty(\mu_{N_+}^{-1}(u))^{C^\infty(N_-)}.
$$
\end{lemma}
{\em Proof.}
By definition $\mu_{N_+}^{-1}(u)$ is a subvariety in $G^*$. Next observe that 
$I_{\chi^{s}}=
{\cal F}(G^*)\cap I_u$. Therefore the algebra ${\cal F}(B_+)={\cal 
F}(G^*)/I_{\chi^{s}}$ is
identified with the algebra of regular functions on $\mu_{N_+}^{-1}(u)$.

Since ${\cal F}(N_-)$ is dense in $C^\infty(N_-)$ on every compact subset in 
$N_-$ we have:
$$
C^\infty(\mu_{N_+}^{-1}(u))^{C^\infty(N_-)}\cong 
C^\infty(\mu_{N_+}^{-1}(u))^{{\cal F}(N_-)}.
$$

Finally observe that action (\ref{actred}) coincides with action 
(\ref{mainactcl}) when restricted to
regular functions.

We shall realize the 
algebra $C^\infty(\mu_{N_+}^{-1}(u))^{C^\infty(N_-)}$ as the algebra of 
functions on a reduced 
Poisson manifold. In the spirit of Lemma \ref{redspace} we shall construct a map 
that forms
a dual pair together with the mapping $\mu_{N_+}$. In this construction we use 
the dressing
action of the Poisson--Lie group $G$ on $G^*$ (see Proposition 
\ref{dressingact}).

Consider the restriction of the dressing action $G^*\times G \rightarrow G^*$ to 
the subgroup $N_+\subset G$.
Note that by Proposition \ref{bpm} (i), (iii) and Proposition \ref{admiss} $N_+$ 
is an admissible 
subgroup in $G$.  
Therefore $C^\infty (G^*)^{N_+}$ is a subalgebra in the Poisson algebra 
$C^\infty (G^*)$.
\begin{proposition}
The algebra $C^\infty (G^*)^{N_+}$ is the centralizer of 
$\mu_{N_+}^*\left( C^\infty \left( N_-\right)\right)$ in the Poisson algebra 
$C^\infty (G^*)$.
\end{proposition}
{\em Proof.}
We shall prove the proposition in a few steps. First we restrict the dressing 
action of $G$ on $G^*$ the 
the Borel subgroup $B_+$.
According to part (iii) of Proposition \ref{bpm} $({\frak b}_+,{\frak b}_-)$ is 
a subbialgebra of 
$({\frak g},{\frak g}^*)$. Therefore $B_+$ is a Poisson--Lie subgroup in $G$. 

By Proposition 
\ref{dressingact} for $X \in {\frak b}_+$ we have:
\begin{equation}
L_{\widehat X} \varphi(L_+,L_-) =( \theta_{G^*}(L_+,L_-) , X ) (\xi_\varphi )= 
(r_-^{-1}\mu_{B_+}^*(\theta_{B_-}) , X) (\xi_\varphi ),
\end{equation}
where $\widehat X$ is the corresponding vector field on $G^*$,
$\xi_\varphi $ is the Hamiltonian vector field of $\varphi \in C^\infty (G^*)$, 
and the map 
$\mu_{B_+}:G^*\rightarrow B_-$ is defined by $\mu_{B_+}(L_+,L_-)=L_-$.
Now from Proposition \ref{bpm} (iv) and the definition of the moment map it 
follows that 
$\mu_{B_+}$ is a moment map for the dressing action of the subgroup $B_+$ on 
$G^*$.

Observe that the orthogonal complement of the Lie subalgebra ${\frak n}_+ 
\subset {\frak b}_+$ 
in the dual space ${\frak b}_-$ coincides with the Lie subalgebra ${\frak 
h}\subset {\frak b}_-$. 
Hence by Proposition \ref{admiss} $N_+$ is an admissible subgroup 
in the Lie--Poisson group $B_+$. Moreover the dual group $B_-$ is the semidirect 
product of the Lie groups 
$H$ and $N_-$ corresponding to the Lie algebras ${\frak n}_+^\perp ={\frak h}$ 
and 
${\frak n}_+^*={\frak n}_-$ , respectively.
We conclude that all the conditions of Proposition \ref{QPmoment}
are satisfied with $A=B_+ , K=N_+ , A^*=B_-, 
T=N_- , K^\perp = H, \mu =\mu_{B_+}$. 
It follows that the algebra $C^\infty (G^*)^{N_+}$ is the centralizer of 
$\mu_{N_+}^*\left( C^\infty \left( N_-\right)\right)$ in the Poisson algebra 
$C^\infty (G^*)$.
This completes the proof.

Let $G^*/N_+$ be the quotient of $G^*$ with respect to the dressing action of 
$N_+$, 
$\pi:G^* \rightarrow G^*/N_+$ 
the canonical projection. Note that the space $G^*/N_+$ is not a smooth 
manifold. However,
in the next section we will see that the subspace $\pi(\mu_{N_+}^{-1}(u))\subset 
G^*/N_+$ is 
a smooth manifold. Therefore by Remark \ref{remred} the algebra 
$C^\infty(\pi(\mu_{N_+}^{-1}(u)))$ 
is isomorphic to $C^\infty(\mu_{N_+}^{-1}(u))^{C^\infty(N_-)}$. Moreover we will 
see that 
$\pi(\mu_{N_+}^{-1}(u))$ has a structure of algebraic variety.
Using Lemma \ref{redreg} we will obtain that the algebra ${\cal F}(B_+)^{{\cal 
F}(N_-)}$
is the algebra of regular functions on this variety.

%%%%%%%%%%%%%%%%%%%%%%%%%%%%%%%%%%%%%%%%%%%%%%%%%%%%%%%%%%%%%%%%%%%%%%%%%%%%%%%%
%%%%%%%%%%%%%%%%%%%%

\section{Cross--section theorem}\label{cross}

\setcounter{equation}{0}
\setcounter{theorem}{0}

In this section we describe the reduced space $\pi(\mu_{N_+}^{-1}(u))\subset 
G^*/N_+$ and the algebra
${\cal F}(B_+)^{{\cal F}(N_-)}$. We also complete the proof of Theorem ${\rm 
B}_q$.

First observe that using the embedding $q:G^*\rightarrow G$ (see Proposition 
\ref{dressingact}) one can
reduce the study of the dressing action to the study of the action of $G$ on 
itself by conjugations.
This simplifies many geometric problems. In particular, consider the restriction 
of this action to
the subgroup $N_+$. Denote by $\pi_q:G\rightarrow G/N_+$ 
the canonical projection onto the quotient with respect to this action. Then we 
can identify the reduced 
space $\pi(\mu_{N_+}^{-1}(u))$ with the subspace $\pi_q(q(\mu_{N_+}^{-1}(u)))$ 
in $G/N_+$. 
Using this identification we shall explicitly
describe the reduced space $\pi(\mu_{N_+}^{-1}(u))$.  
We start with description of the image of the ``level surface'' 
$\mu_{N_+}^{-1}(u)$ under the embedding
$q$. 
\begin{proposition}\label{constrt}
Let $q:G^*\rightarrow G$ be the map introduced in Proposition \ref{dressingact},
$$
q(L_+,L_-)=L_-L_+^{-1}.
$$
Then $q(\mu_{N_+}^{-1}(u))$ is a subvariety in $N_+sN_+$. 
\end{proposition}
{\em Proof.} First, using definition (\ref{mun}) of the map $\mu_{N_+}$ we can 
describe the space 
$\mu_{N_+}^{-1}(u)$ as 
follows:
\begin{equation}\label{mun1}
\mu_{N_+}^{-1}(u)=\{(h_+n_+,s(h_+)u) | n_+ \in N_+ , h_+ \in H \},
\end{equation}
since by (\ref{fact}) $h_-=s(h_+)$. Therefore
\begin{equation}\label{dva}
q(\mu_{N_+}^{-1}(u))=
\{ s(h_+)un_+^{-1}h_+^{-1}| n_+ \in N_+ , h_+ \in H \}.
\end{equation}

Now recall that $u \in N_+sN_+\cap N_-$, and hence
\begin{equation}\label{ras} 
un_+^{-1}\in N_+sN_+.
\end{equation}

Next, the space $N_+sN_+$ is invariant with respect to 
the following action of $H$:
\begin{equation}\label{tri}
h\circ L= s(h)Lh^{-1}.
\end{equation}

Indeed, let $L=vsu,~ v,u \in N_+$ be an element of $N_+sN_+$. Then
\begin{equation}
h\circ L=s(h)vs(h)^{-1}s(h)sh^{-1}huh^{-1}=s(h)vs(h)^{-1}shuh^{-1}.
\end{equation}
The r.h.s. of the last equality belongs to $N_+sN_+$ because $H$ normalizes 
$N_+$. 

Comparing action (\ref{tri}) with (\ref{dva}) and adding (\ref{ras}) we obtain 
that
$q(\mu_{N_+}^{-1}(f)) \subset N_+sN_+$. Since $q$ is an embedding, 
$q(\mu_{N_+}^{-1}(f))$ is a subvariety in $N_+sN_+$.
This concludes the proof.

We identify $\mu_{N_+}^{-1}(u)$ with the subvariety in 
$N_+sN_+$ described in the previous proposition. 
As we observed in the beginning of this section
the reduced space $\pi(\mu_{N_+}^{-1}(u))$ is isomorphic to 
$\pi_q(q(\mu_{N_+}^{-1}(u)))$. 
Note that by Proposition \ref{constrt} $q(\mu_{N_+}^{-1}(u))\subset N_+sN_+$. 
But the variety $N_+sN_+$ 
is stable under the action of $N_+$ by conjugations. Therefore to describe
the reduced space $\pi_q(q(\mu_{N_+}^{-1}(u)))$ we have to study the structure 
of the quotient 
$N_+sN_+/N_+$. Our main geometric result is
\vskip 0.3cm
\noindent
{\bf Theorem $\bf C_q$ (\cite{SS}, Theorem 3.1)}
{\em Let $N_+'=\{ v \in N_+|s^{-1}(v)\in N_- \}$.
Then the action of $N_+$ on $N_+sN_+$ by conjugations is free, and 
$N_+'s$ is a cross--section for this action, i.e. 
for each $L\in N_+sN_+$ there exists a unique
element $n\in N_+$ such that $n L n^{-1}\in N_+'s$.
Moreover, the projection $\pi_q: N_+sN_+\rightarrow N_+'s$ is a morphism of 
varieties.}
\begin{lemma}{\bf (\cite {Ch}, Theorem 8.4.3, \cite{st} Lemma 7.2)}
$N_+'\subset N$ is an
abelian subgroup, $\dim N_+'=l.$ Moreover,
every element $L\in N_+sN_+$ may be uniquely represented in the form $L=vsu ,
v\in N_+' , u\in N_+$.
\end{lemma}
\noindent
{\em Proof of Theorem $ C_q$.}
Denote by $h$ the Coxeter number of $\frak g$. By definition $h$ is the order of 
the Coxeter element, 
$s^h=id$. Note that $h={2N \over l}$.

Let $C_s\subset W$ be the cyclic subgroup generated by the Coxeter element. $%
C_s$ has exactly $l$ different orbits in $\Delta$. 
The proof depends on the structure of these orbits. For this
reason we have to distinguish several cases\footnote{ The proofs given below do 
not apply when $\frak g$
is the simple Lie algebra of type $E_6$.}.

1.Let ${\frak g}$ be of type $ A_l.$

The following lemma is checked by straightforward calculation.
\begin{lemma}
(i) Each orbit of $C_s$ in $\Delta$ consists of exactly $h$ elements. 
One can order these orbits in such a way
that $k$-th orbit contains all positive roots of height $k$ and all negative
roots of height $h-k.$
\end{lemma}

Put 
\[
{\frak n}_k=\bigoplus_{\left\{ \alpha \in \Delta _{+},\;ht\,\alpha
=k\right\} }{\frak n}_\alpha ,N_k=\exp {\frak n}_k. 
\]
For each $k$ we can choose $\gamma _k\in \Delta _{+}$ in such a way that 
\[
{\frak n}_k=\bigoplus_{p=0}^{h-k-1}{\frak n}_{s^{-p}\left( \gamma_k\right) }. 
\]
Put ${\frak n}_k^p={\frak n}_{s^{-p}\left( \gamma_k\right) },N_k^p=\exp {\frak %
n}_k^p.$

Let $L=vsu,v\in N_+',u\in N_+$. We must
find $n\in N_+$ such that 
\begin{equation}
nvsu=v_0sn,~v_0\in  N_+'. \label{nx}
\end{equation}
For any $n\in N_+$ there exists a factorization
\[
n=n_1n_2\ldots n_l,\mbox{ where } n_k\in N_k. 
\]
Moreover, each $n_k$ may be factorized as 
\[
n_k=n_k^0n_k^1\ldots n_k^{h-k-1},\;n_k^p\in N_k^p. 
\]

For any $n\in N_+$ the element $nvsu$
admits a representation 
\[
nvsu=\tilde vs\tilde u,\;\tilde v\in N_+',~
\tilde u\in N_+. 
\]
Let 
\[
\tilde u=\overrightarrow{\prod_{k=1}^l}\overrightarrow{\,\prod_{p=0}^{h-k-1}}%
\tilde u_k^p,\;\tilde u_k^p\in N_k^p, 
\]
be the corresponding factorization of $\tilde u.$ 
\begin{lemma}
We have $\tilde u_k^p= s^{-1}\left( n_k^{p-1}\right) V_k^p,$ where the
factors $V_k^p\in N_k^p$ depend only on $u,v$ and on $n_j^q$ with $%
j<k.$
\end{lemma}

Assume now that $n$ satisfies (\ref{nx}). Then we have $\tilde v=v_0,\tilde
u=n.$ This leads to the following relations: 
\begin{equation}
s^{-1}\left( n_k^{p-1}\right) V_k^p=n_k^p,  \label{recur}
\end{equation}
where we set formally $n_k^{-1}=1.$
\begin{lemma}
The system (\ref{recur}) may be solved recursively starting with $k=1,$ $%
p=0. $
\end{lemma}
Clearly, the solution is unique. This concludes the proof for ${\frak g}$ of
type $A_l.$

2. Let now ${\frak g}$ be a simple Lie algebra of type other than $A_l$ and $%
E_6$.

\begin{lemma}
(i) The Coxeter number $h$ is even. 

(ii) Each orbit
of $C_s$ in $\Delta$ consists of exactly $%
h$ elements and contains an equal number of positive and negative roots.

(iii) Put 
\[
\Delta _{+}^p=\{\alpha \in \Delta _{+};s^{p}\alpha \notin \Delta
_{+}\},\;{\frak n}^p=\bigoplus_{\alpha \in \Delta _{+}^p}{\Bbb C}\cdot
X_\alpha ; 
\]
then ${\frak n}^p\subset {\frak n}$ is an abelian subalgebra, $\dim {\frak n}%
^p=l.$
\end{lemma}

If  $\frak g$ is not of type $D_{2k+1}$ this assertion follows from  
Proposition 33, Chap.6, no. 1.11 and Corollary 3, Chap.5, no. 6.2 in \cite{Bur}. 
For $\frak g$ of type $D_{2k+1}$ it may be 
checked directly.

Put $N^p=\exp {\frak n}^p$. Let $N^p$ be the corresponding subgroup
of $G.$ Let $L=vsu,~v\in 
N_+',~u\in N_+$. We must find $n\in N_+$ such that 
\[
vsu=nv_0sn^{-1},~v_0\in N_+'. 
\]

Put 
\begin{equation}
n=n_1n_2\ldots n_{\frac h2},~~n_p\in {N}_p.  \label{n}
\end{equation}
The elements $n_p$ will be determined recursively. 
We have 
\begin{equation}
vs\left( u\right) =\overrightarrow{\prod_p }n_p v_0 s\left( 
\overleftarrow{\prod_p}n_p^{-1}\right).  
\label{nnn}
\end{equation}
We shall say that an element $x\in G$ is in the big cell in $G$
if $x \in B_+N_-\subset G.$
\begin{lemma}
$vs\left( u\right) $ is in the big cell in $G$ and admits
a factorization 
\[
vs\left( u\right) =x_{+}^1 x_{-}^1,\;x_{+}^1\in N_+%
,\;x_{-}^1\in N_-. 
\]
\end{lemma}

Indeed, let $u=u_{h/2}u_{h/2-1}\ldots u_1,~u_p\in N^p$ be a
similar decomposition of $u$. Then we have simply $x_{-}=s\left(
u_1\right) .\ $(It is clear that $x_{+}^1\in B_+$ actually does not
have an $H$-component and so belongs to $N_+$

A comparison of the r.h.s in (\ref{nnn}) with the Bruhat decomposition of
the l.h.s. immediately yields that the first factor in (\ref{n}) is given by 
$n_1=s^{-1}\left( x_{-}\right) ^{-1}.$

Assume that $n_1,n_2,\ldots ,n_{k-1}$ are already computed. Put 
\[
m_k=n_1n_2\ldots n_{k-1} 
\]
and consider the element 
\begin{equation}
L^k:=s^{k-1}\left( m_k^{-1}vs(u)s(m_k)\right) .
\label{Lk}
\end{equation}

\begin{lemma}
$L^k$ is in the big cell in $G$ and admits a factorization 
\begin{equation}
L^k=x_{+}^kx_{-}^k,\;x_{+}^k\in N_+,\;x_{-}^k\in N_-.
\label{fk}
\end{equation}
\end{lemma}

The elements $x_{\pm }^k$ are computed recursively from the known
quantities. By applying a similar transform to the r.h.s. of (\ref{nnn}) we
get 
\begin{eqnarray}
L^k &=&s^{k-1}\left( m_k^{-1}  
\overrightarrow{\prod_p}n_p v_0s\left( 
\overleftarrow{\prod_p}n_p^{-1}\right) s(m_k) \right)
=  \label{Lrhs}
\end{eqnarray}
\[
s^{k-1}\left( \overrightarrow{\prod_{p\geq k}}n_p
v_0\right) s^{k}\left( \overleftarrow{\prod_{p\geq k+1}} 
n_p^{-1}\right)s^{k}\left( n_k^{-1}\right) . 
\]
Comparison of (\ref{Lrhs}) and (\ref{fk}) yields $x_{-}^k=s^{k}\left(
n_k^{-1}\right)$. Hence $n_k=s^{-k}\left( x_{-}^k\right)^{-1}$, which concludes
the induction.

Finally observe that by construction the map $\pi_q: N_+sN_+\rightarrow N_+'s$ 
is a morphism of varieties.

\begin{corollary}\label{var} 
The space
$\pi(\mu_{N_+}^{-1}(u))$ is a subvariety in $N_+'s$. The algebra
${\cal F}(B_+)^{{\cal F}(N_-)}$ is isomorphic to the algebra of regular 
functions on
$\pi(\mu_{N_+}^{-1}(u))$.
\end{corollary}
{\em Proof.}
First observe that by construction
$\pi(\mu_{N_+}^{-1}(u))\cong \pi_q(q(\mu_{N_+}^{-1}(u)))$ is a subvariety in 
$N_+'s$.
In particular, $\pi(\mu_{N_+}^{-1}(u))$ is a smooth manifold. Hence by Remark 
\ref{remred}
the map
$$
C^\infty(\pi(\mu_{N_+}^{-1}(u)))\rightarrow 
C^\infty(\mu_{N_+}^{-1}(u))^{C^\infty(N_-)},~~\psi \mapsto \pi^*\psi
$$
is an isomorphism.

Now observe that by construction the map $\pi: \mu_{N_+}^{-1}(u)\rightarrow 
\pi(\mu_{N_+}^{-1}(u))$
is a morphism of varieties. Therefore if $\psi \in {\cal 
F}(\pi(\mu_{N_+}^{-1}(u)))$ then
$\pi^*\psi$ is a regular function on $\mu_{N_+}^{-1}(u)$. Conversely, suppose 
that $\varphi \in 
{\cal F}(\mu_{N_+}^{-1}(u))\cap C^\infty(\mu_{N_+}^{-1}(u))^{C^\infty(N_-)}$.
Note that  $\pi(\mu_{N_+}^{-1}(u))$ may be regarded as a subvariety in 
$\mu_{N_+}^{-1}(u)$ 
(see Remark \ref{redpoisalg}). Then the restriction of 
$\varphi$ to $\pi(\mu_{N_+}^{-1}(u))\subset \mu_{N_+}^{-1}(u)$ is a regular 
function.
Therefore the map
$$
{\cal F}(\pi(\mu_{N_+}^{-1}(u)))\rightarrow {\cal F}(\mu_{N_+}^{-1}(u))\cap 
C^\infty(\mu_{N_+}^{-1}(u))^{C^\infty(N_-)},~~\psi \mapsto \pi^*\psi
$$
is an isomorphism. 

Finally observe that by Lemma \ref{redreg} the algebra ${\cal 
F}(\mu_{N_+}^{-1}(u))\cap 
C^\infty(\mu_{N_+}^{-1}(u))^{C^\infty(N_-)}$ is isomorphic to ${\cal 
F}(B_+)^{{\cal F}(N_-)}$.
This completes the proof.

Theorem ${\rm C}_q$ is a group counterpart of Theorem C. Moreover the space
$N_+'s$ naturally appears in the study of regular elements in $G$. Recall that 
an element of 
$G$ is called regular if its centralizer in $G$ is of minimal possible 
dimension. Let $R$ be 
the set of regular elements in $G$. Clearly, $R$ is stable under the action of 
$G$ on itself
by conjugations and in fact $R$
is the union of all $G$ orbits in $G$ of maximal dimension.
A function $\psi$ on $G$ is called a class function if $f(x)=f(y)$ whenever $x$ 
and $y$ are conjugate
points of definition of $\psi$.
We denote by ${\cal F}^G(G)$ the algebra of regular class functions on $G$.
\vskip 0.3cm
\noindent
{\bf Theorem $\bf D_q$ (\cite{st}, Theorems 1.4 and 6.1)}
{\em Let $G$ be a complex connected simply connected simple algebraic group. 
Then 
The space $N_+'s$ is contained in $R$ and is a cross--section for the action of 
$G$ on $R$. 
That is every $G$--orbit in $G$ of maximal dimension intersects $N_+'s$ in one 
and only one point.
The algebra of regular class functions on $G$ is freely generated as a 
commutative algebra
over $\Bbb C$ by the characters of fundamental representations of $G$, 
$\chi_1,\ldots , \chi_l$.
Moreover, $N_+'s$ is an algebraic variety, and the algebra of regular functions 
on $N_+'s$ is 
freely generated as a commutative algebra over $\Bbb C$ by the restrictions of 
the characters 
$\chi_1,\ldots , \chi_l$ to $N_+'s$ .}
\vskip 0.3cm
\noindent
{\bf Theorem $\bf E_q$} 
{\em For any $\psi \in {\cal F}^G(G)$ one has $\rho_{\chi^{s}}(p^*\psi) \in 
{\cal F}(B_+)^{{\cal F}(N_-)}$.
Furthermore the map
\begin{equation}\label{isoreg}
{\cal F}^G(G)\rightarrow {\cal F}(B_+)^{{\cal F}(N_-)}, \psi \mapsto 
\rho_{\chi^{s}}(p^*\psi)
\end{equation}
is an algebra isomorphism. In particular, 
$$
{\cal F}(B_+)^{{\cal F}(N_-)}={\Bbb C}[\rho_{\chi^{s}}(p^*\chi_1),\ldots , 
\rho_{\chi^{s}}(p^*\chi_l)]
$$
is a polynomial algebra in $l$ generators.}
\vskip 0.3cm
\noindent
{\em Proof.} 
Let $\psi$ be an element of ${\cal F}^G(G)$. The restriction of $\psi$ to 
the subvariety $\pi(\mu_{N_+}^{-1}(u))\cong \pi_q(q(\mu_{N_+}^{-1}(u)))\subset 
N_+'s \subset G$ 
is a regular function. Using the isomorphism
${\cal F}(\pi(\mu_{N_+}^{-1}(u)))\cong {\cal F}(B_+)^{{\cal F}(N_-)}$ (see 
Corollary \ref{var}) this 
restriction may be identified with
$\rho_{\chi^{s}}(q^*\psi)\in {\cal F}(B_+)^{{\cal F}(N_-)}$.
By Theorem ${\rm D}_q$ the 
algebra ${\cal F}(N_+'s)$ is freely generated as a commutative algebra over 
$\Bbb C$ by the restrictions of 
the fundamental characters
$\chi_1,\ldots , \chi_l$. Since $\pi(\mu_{N_+}^{-1}(u))$ is a subvariety in 
$N_+'s$ the algebra 
${\cal F}(\pi(\mu_{N_+}^{-1}(u)))$ is generated by the restrictions of the 
fundamental characters
$\chi_1,\ldots , \chi_l$. Therefore the map (\ref{isoreg}) is surjective.
We have to prove that it is injective.

Let $\chi_i$ be a fundamental character.
Consider the restriction of the function $\rho_{\chi^{s}}(q^*\chi_i)$ to the 
subspace in 
$\mu_{N_+}^{-1}(u)$ formed by elements (see (\ref{mun1})):
$$
(h_+,s(h_+)u),~h_+\in H.
$$
Then $\rho_{\chi^{s}}(q^*\chi_i)(h_+,s(h_+)u)=\chi_i(s(h_+)uh_+^{-1})$. Since 
$\chi_i$ is a character
we have $\chi_i(s(h_+)uh_+^{-1})=\chi_i(h_+^{-1}s(h_+)u)$. The element $u$ is 
unipotent, and hence
$\chi_i(h_+^{-1}s(h_+)u)=\chi_i(h_+^{-1}s(h_+))$. Now recall that the 
restrictions of the fundamental 
characters to the Cartan subgroup are algebraically independent (they are given 
by the well--known Weyl
formula). Therefore (\ref{isoreg}) is an isomorphism. This completes the proof.
\vskip 0.3cm
\noindent
{\em Proof of Theorem $B_q$.}
Let $p:{\cal F}_h(G^*)\rightarrow {\cal F}(G^*)$ be the map defined in 
Proposition \ref{quantreg}.
Let $W_h^{Rep}(B_+)$ be the subalgebra in 
$W_h(B_+)$ topologically generated by the elements 
$C_{V_i}^{\rho_{\chi_h^{s}}}=\rho_{\chi_h^{s}}(C_{V_i}),~i=1,\ldots ,l$.
From the the definition of the elements $C_{V_i}^{\rho_{\chi_h^{s}}}$ it follows 
that
$p(C_{V_i}^{\rho_{\chi_h^{s}}})=\rho_{\chi^{s}}(p^*\chi_i)$.
Therefore by Theorem ${\rm E}_q$ $p(W_h^{Rep}(B_+))={\cal F}(B_+)^{{\cal 
F}(N_-)}$,
and $W_h^{Rep}(B_+)$ is freely generated as a commutative topological algebra 
over ${\Bbb C}[[h]]$ by the elements 
$C_{V_i}^{\rho_{\chi_h^{s}}}=\rho_{\chi_h^{s}}(C_{V_i}),~i=1,\ldots ,l$.

On the other hand using the definitions of the algebras ${\cal F}_h(B_+)^{{\cal 
F}_h(N_-)}$ and
${\cal F}(B_+)^{{\cal F}(N_-)}$ it is easy to see that $p({\cal F}_h(B_+)^{{\cal 
F}_h(N_-)})=
{\cal F}(B_+)^{{\cal F}(N_-)}$. We shall prove that $W_h^{Rep}(B_+)$ is 
isomorphic to
${\cal F}_h(B_+)^{{\cal F}_h(N_-)}$.

Let $I\in {\cal F}_h(B_+)^{{\cal F}_h(N_-)}$ be an invariant element.
Then $p(I)\in {\cal F}(B_+)^{{\cal F}(N_-)}$, and hence one can find an element 
$K_0\in 
W_h^{Rep}(B_+)$ such
that $I-K_0=hI_1,~I_1\in {\cal F}_h(B_+)^{{\cal F}_h(N_-)}$. Applying the same 
procedure
to $I_1$ one can find elements $K_1\in W_h^{Rep}(B_+),
~I_2\in {\cal F}_h(B_+)^{{\cal F}_h(N_-)}$ such that $I_1-K_1=hI_2$, i.e. 
$I-K_0-hK_1=0~(\mbox{mod }h^2)$. 
We can continue this process. Finally we obtain an infinite sequence of elements 
$K_i\in W_h^{Rep}(B_+)$ such that $I-\sum_{i=0}^p h^pK_p=0~(\mbox{mod 
}h^{p+1})$. Since the space 
${\cal F}_h(B_+)$ is complete in the $h$--adic topology the series 
$\sum_{i=0}^\infty h^pK_p\in W_h^{Rep}(B_+)$ 
converges to $I$. Therefore $I\in W_h^{Rep}(B_+)$, and hence
${\cal F}_h(B_+)^{{\cal F}_h(N_-)}$ is isomorphic to $W_h^{Rep}(B_+)$. 

We also have the following inclusions:
$$
W_h^{Rep}(B_+)\subseteq W_h(B_+)\subseteq {\cal F}_h(B_+)^{{\cal F}_h(N_-)}
\cong W_h^{Rep}(B_+).
$$
Therefore $W_h^{Rep}(B_+)$ coincides with $W_h(B_+)$.
This proves part (ii) of Theorem ${\rm A}_q$ and Theorem ${\rm B}_q$.

\chapter*{Acknowledgements}

The first words of gratitude are due to Professor Michael Semenov--Tian--Shansky 
for guidance into
the world of Modern Mathematical Physics. This thesis is essentially based on 
his ideas.

I am greatly indebted my advisor in Uppsala Doctor Anton Alekseev for his 
continuous support,
for many insightful discussions and for continuous encouragement.

I gratefully acknowledge Prof. Antti Niemi and his research group at the 
Department 
of Theoretical Physics at Uppsala University for providing 
an excellent working atmosphere.

The Department of Theoretical Physics in Uppsala University has provided 
excellent
facilities for research, which are gratefully acknowledged.

I wish to thank MSc Mats Lilja for his patient help in many and various 
practical matters
during my stay in Uppsala.

\pagestyle{plain}

\end{document}